\documentclass[a4paper]{article}[11pts]

\usepackage[english]{babel}
\usepackage[utf8x]{inputenc}
\usepackage[T1]{fontenc}

\normalsize

\usepackage[a4paper,top=1in,bottom=1in,left=1in,right=1in,marginparwidth=1in]{geometry}

\usepackage{setspace}
\usepackage{amsmath}
\usepackage{amssymb} 
\usepackage{amsthm}
\usepackage{amsfonts, mathtools}
\usepackage{algorithm,algpseudocode}
\usepackage{bm}
\usepackage{comment}
\usepackage{graphicx}
\usepackage{multirow, booktabs}
\usepackage{caption}
\usepackage{subcaption}
\usepackage[colorinlistoftodos]{todonotes}
\usepackage[colorlinks=true, allcolors=blue]{hyperref}

\newcommand{\Bb}{\bm{B}}
\newcommand{\Dd}{\bm{b}}
\newcommand{\Aa}{\bm{a}}

\newcommand{\E}{\mathbb{E}}
\newcommand{\R}{\mathbb{R}}
\newcommand{\prob}{\mathbb{P}}

\newtheorem{proposition}{Proposition}
\newtheorem{lemma}{Lemma}
\newtheorem{example}{Example}
\newtheorem{assumption}{Assumption}
\newtheorem{theorem}{Theorem}
\newtheorem{corollary}{Corollary}

\usepackage{natbib}

\title{An Improved Analysis of LP-based Control for Revenue Management}
\author{}
\date{}

\author{Guanting Chen$^\dagger$ \and Xiaocheng Li$^\Diamond$ \and  Yinyu Ye$^\ddagger$}
\date{\small 
$^\dagger$ Institute for Computational and Mathematical Engineering, Stanford University\\
$^\Diamond$ Imperial College Business School, Imperial College London\\
$^\ddagger$Department of Management Science and Engineering, Stanford University\\
$\{$guanting, chengli1, yinyu-ye$\}$@stanford.edu\\
}

\begin{document}
\maketitle

\onehalfspacing

\begin{abstract}
In this paper, we study a class of revenue management problems where the decision maker aims to maximize the total revenue subject to budget constraints on multiple type of resources over a finite horizon. At each time, a new order/customer/bid is revealed with a request of some resource(s) and a reward, and the decision maker needs to either accept or reject the order. Upon the acceptance of the order, the resource request must be satisfied and the associated revenue (reward) can be collected. We consider a stochastic setting where all the orders are i.i.d. sampled, i.e., the reward-request pair at each time is drawn from an unknown distribution with finite support. The formulation contains many classic applications such as the quantity-based network revenue management problem and the Adwords problem. We focus on the classic LP-based adaptive algorithm and consider \textit{regret} as the performance measure defined by the gap between the optimal objective value of the certainty-equivalent linear program (LP) and the expected revenue obtained by the online algorithm. Our contribution is two-fold: (i) when the underlying LP is nondegenerate, the algorithm achieves a problem-dependent regret upper bound that is independent of the horizon/number of time periods $T$; (ii) when the underlying LP is degenerate, the algorithm achieves a regret upper bound that scales on the order of $\sqrt{T}\log T$. To our knowledge, both results are new and improve the best existing bounds for the LP-based adaptive algorithm in the corresponding setting. We conclude with numerical experiments to further demonstrate our findings. 

\end{abstract}

\section{Introduction}
Consider a revenue management problem over a horizon of $T$ time periods. The objective is to maximize the cumulative reward over the horizon subject to the budget constraints on $m$ types of resource. At each time period, a customer order arrives and it requests certain amount of each resource. We need to decide whether to accept or reject the customer order. Upon the acceptance of the order, we need to satisfy the resource request, and we will collect the revenue/reward associated with the order. From the resource viewpoint, we act as a market maker who allocates the resources among all the customer orders. The allocation decisions are made in an online manner and all the past decisions are irrevocable. In this paper, we study a stochastic model where the request-reward pair made by each customer order is assumed to be i.i.d. and follows an unknown distribution with finite support. 

Specifically, the considered problem takes the following LP as its underlying form:
\begin{align}
   \max \ \ & \sum_{t=1}^T \bm{r}^\top_t \bm{x}_t \label{eqn:multiILP}  \\
    \text{s.t. }\ & \sum_{t=1}^T \bm{A}_t \bm{x}_t \le \bm{B} \nonumber  \\ 
    & \bm{1}^\top \bm{x}_t \le 1, \ \ \bm{x}_t \ge \bm{0}, \ \ t=1,...,T \nonumber
\end{align}
where $\bm{r}_t = (r_{1t},...,r_{kt})^\top \in \R^k$, $\bm{A}_t = (\bm{a}_{1t},...,\bm{a}_{kt}) \in \R^{m\times k}$, and $\bm{a}_{st} = (a_{1st},...,a_{mst})^\top\in\mathbb{R}^m,$ for $t=1,...,T$ and $s=1,...,k.$ The right-hand-side vector $\bm{B}=(B_1,...,B_m)^\top$ encapsulates the capacity for each resource. The decision variables are $\bm{x}=\left(\bm{x}_1,...,\bm{x}_T\right)$ where $\bm{x}_t=(x_{1t},...,x_{kt})^\top$ for $t=1,...,T$. In an online setting, the parameters of the optimization problem (\ref{eqn:multiILP}) are revealed in an online fashion and one needs to determine the value of decision variables sequentially. At each time $t,$ the coefficients $(\bm{r}_t, \bm{A}_t)$ are revealed, and we need to decide the value of $\bm{x}_t$ instantly. Different from the offline setting, at time $t$, we do not have the information of the subsequent coefficients to be revealed, i.e., $\{(\bm{r}_{t'}, \bm{A}_{t'})\}_{t'=t+1}^T$. The problem \eqref{eqn:multiILP} in an online setting is often referred to as \textit{online linear programming} \citep{agrawal2014dynamic,kesselheim2014primal}. With different specification of the input of the LP, the problem encompasses a wide range of applications, including secretary problem \citep{ferguson1989solved}, knapsack problem \citep{kellerer2003knapsack},  network routing problem \citep{buchbinder2009online}, matching and Adwords problem \citep{mehta2005adwords}, service reservation and scheduling problem \citep{conforti2014integer, stein2020advance}, network revenue management problem \citep{talluri2006theory}, order fulfillment problem \citep{acimovic2015making,jasin2015lp} and more.

We generally refer to the problem as the online revenue management problem. In this paper, we analyze the performance of the classic LP-based adaptive algorithm under both the cases when the underlying LP is nondegenerate or degenerate.

Our contribution is two-fold:
\begin{itemize}
\item When the underlying LP is nondegenerate, we derive a regret bound independent of the horizon $T$ which improves the existing logarithmic bound in \cite{jasin2015performance}. The result highlights that the online revenue management is different in its information-theoretic nature from other constrained online learning problem such as bandits with knapsacks problem (BwK, \cite{badanidiyuru2013bandits}) and online convex optimization with constraints (OCOwC, \cite{yu2017online}). Specifically, the online revenue management problem concerns a full-information setting where the decision maker first observes the information (customer order) and then makes the decision of acceptance and rejection. In contrast, the problems of BwK and OCOwC are in a partial-information setting where the decision is made prior to the observation, and thus bounded regret is not achievable in such setting. Moreover, the finite support condition on the distribution is critical for the achievability of bounded regret. \cite{bray2019does} shows an $\Omega(\log T)$ lower bound when the underlying distribution is continuously supported. 

\item When the underlying LP is degenerate, we derive a regret bound on the order of $\sqrt{T}\log T$. The LP-based adaptive algorithm that re-optimizes the decision rule at every time period (also known as frequent re-solving) is criticized for performance deterioration when the underlying LP becomes nearly degenerate. The existing bounds for the algorithm \citep{jasin2012re, jasin2015performance, wu2015algorithms} can be arbitrarily large when the underlying LP is still nondegenerate but close to degenerate. Our analysis provides a positive result that the algorithm achieves $O(\sqrt{T}\log T)$ regret upper bound regardless of whether the underlying LP is nondegenerate or not. The analysis is nearly tight in that an $\Omega(\sqrt{T})$ lower bound is established by \cite{bumpensanti2020re}.

\end{itemize}

In the following, we review the related literature and highlight the technical novelty of our result and analysis.

\subsection{Literature Review}

There have been a proliferate literature on algorithm design and analysis for the online revenue management problem. We provide a short summary of our result and the related literature in Table \ref{tab:result}. In the following, we will elaborate on our contribution upon the existing literature. 

\begin{table}[ht!]
    \centering
    \begin{tabular}{c|c|c|c}
    \toprule
         & Benchmark & Regret Bound &  Key Assumption(s)\\
         \midrule
        \cite{jasin2012re} & Fluid & Bounded & Nondegeneracy, distribution knowledge\\ 
        \cite{jasin2015performance} & Fluid &$\tilde{O}(\log T)$ & Nondegeneracy \\ 
         \cite{vera2019bayesian} & Hindsight & Bounded & Distribution knowledge \\
         \cite{bumpensanti2020re} & Hindsight & Bounded & Distribution knowledge \\
         \cite{asadpour2019online} & Full flex. & Bounded  & Long-chain design, $\xi$-Hall condition \\
          \cite{veradynamic} & Hindsight & Bounded & Replenishment, binary entries\\ 
        Ours & Fluid & Bounded & Nondegenarcy \\
        Ours & Fluid &  $\tilde{O}(\sqrt{T})$ & -- \\
         \bottomrule
    \end{tabular}
    \caption{Result comparison against literature: The precise definitions of the fluid benchmark and the hindsight benchmark are given in \eqref{eqn:DLP} and \eqref{eqn:Hind}. We note that all the regret bounds here are problem-dependent bound that mainly focuses on the dependence on horizon length $T$ but will inevitably involve certain parameters related to the underlying distribution/optimization problem.}
    \label{tab:result}
\end{table}

\paragraph{\textbf{Known distribution}}\

One stream of literature investigates the existence of an online algorithm that achieves bounded regret under a known distribution setting. Specifically, a line of works \citep{jasin2012re, wu2015algorithms, bumpensanti2020re} study the canonical quantity-based network revenue management problem and design algorithms that achieve bounded regret under the knowledge of the order arrival distribution. \cite{jasin2012re} show that the LP-based re-solving algorithm achieves bounded regret where the algorithm computes a control policy by periodically solving a linear program specified by the known arrival distribution. A critical condition in \cite{jasin2012re} is that the underlying LP should be nondegenerate. A subsequent work \citep{bumpensanti2020re} fully solves the problem through an infrequent re-solving scheme which updates the control policy only at a few selected points. Analysis-wise, one key difference of our analysis compared to \cite{jasin2012re} and \cite{bumpensanti2020re} is a new regret decomposition (Proposition \ref{prop_general_regret_form}). The regret decomposition helps us to circumvent the complication of dealing with the partially-accepted order types and transforms the analysis of partially-accepted order types to that of the constraint process. In this light, the analyses in \cite{jasin2012re} and \cite{bumpensanti2020re} are purely primal-based while our analysis also utilizes the dual problem.

Another line of works \citep{vera2019online, vera2019bayesian, banerjee2020uniform} devise algorithms that achieve bounded regret for various problems including dynamic pricing, knapsack problem, and bin-packing problem. The authors develop a novel and intuitive approach called ``compensated coupling'' to derive regret upper bounds. The idea is to bound the cumulative expected loss induced by the different decisions made by the online algorithm against the hindsight optimal which has all the future information (formally defined in \eqref{eqn:Hind}). As in the aforementioned works on network revenue management, this line of works are also built upon the knowledge of the underlying distribution and the algorithms are closely related to the idea of approximate dynamic programming.

\paragraph{\textbf{Unknown distribution}}\

Efforts have also been made to relax the assumption on knowing the underlying distribution. \cite{banerjeeconstant} consider one single historical trace of observations in substitute of the knowledge of the true distribution and derive an algorithm that achieves bounded regret for certain constrained online optimization problem. In certain sense, the single historical trace can contribute to $\Omega(n)$ observations (with $n$ being the length of horizon). In a similar spirit, \cite{shivaswamy2012multi} characterize the number of historical observations needed to achieve a bounded regret for the stochastic multi-armed bandits problem. The historical observations required in these two works can be viewed as a warm start for the online procedure, and the empirical distribution constructed from the historical observations provides a moderately good estimation for the true underlying distribution at the very beginning of the online procedure. 

As a subsequent work of \cite{jasin2012re}, \cite{jasin2015performance} studies the same problem as ours and proves an $O(\log T)$ regret for the LP-based adaptive algorithm. Apart from the new regret decomposition mentioned above, we now provide the second explanation for why we manage to improve the regret bound. First, we point out that the key of the analysis under known distribution in \cite{jasin2012re} is a martingale analysis of the constraint process. A recent work \citep{balseiro2021survey} extends the result in \cite{jasin2012re} to more general constrained revenue management problems but still requires the distribution knowledge. The common point of \cite{jasin2012re} and \cite{balseiro2021survey} is that with the knowledge of the underlying distribution, the constraint process is a martingale under the adaptive (re-solving) scheme. When the distribution is unknown, the constraint process is no longer a martingale. \cite{jasin2015performance} leaves the martingale approach and defines a sequence of high probability events to analyze the constraint process. In contrary, we still stick to the martingale approach: different from the known distribution case, we identify conditions under which the constraint process is not a martingale, but it still behaves ``stably'' to suffice for the bounded regret analysis. The conditions need to be carefully tuned so that they can be met with high probability. We remark that the idea applies to both cases of nondegenerate and degenerate. Under both cases, we first start to identify the desirable condition on the constraint process to achieve the corresponding regret bound and then show that such condition can be met with high probability. 

A stream of works \citep{agrawal2014dynamic, bray2019does, li2021online} also consider the problem under unknown distribution but impose no finite support condition. Similar to our work, the primal-dual approach is often used in such setting for both algorithm design and analysis. Both lower and upper bounds of order $\log T$ are established for the case when the distribution support is infinite. In this light, the finite support condition is critical in the achievability of bounded regret. Intuitively, when the LP is nondegenerate, the finite support condition creates a well-separatedness between different order types, and the separation can be learned with high probability using $O(1)$ number of samples. This contrast between the finite and infinite support is analogous to the findings in the (unconstrained) newsvendor problem \citep{besbes2013implications}.

\paragraph{\textbf{Other related works}}\

Another work related to our result is \cite{asadpour2019online}, where the authors derive bounded regret for a resource allocation problem to study the effectiveness of the long-chain design without knowing the true distribution. Technically, the formulation in \cite{asadpour2019online} can be cast in (\ref{eqn:multiILP}) by imposing a binary structure on the constraint matrices $\bm{A}_t$'s, along with certain other conditions. Another recent work \citep{veradynamic} on a revenue management/resource allocation problem considers a setting where there is resource replenishment. Similar to our work, the paper focuses on analyzing the resource consumption process and takes a geometric viewpoint. In addition to the resource replenishment, the paper assumes the entries in the constraint matrix is binary and all the resources are binding, while both assumptions are critical in the regret derivation. Two subsequent works \citep{kerimov2021dynamic, kerimov2021optimality} consider a dynamic matching problem and develop algorithms achieving bounded regret utilizing the structure of the underlying LP. 



\section{Model, Assumption, and Algorithm}
\label{sec_model}

In this section, we formulate the problem and present the LP-based adaptive algorithm. We first focus on the one-dimensional case for notation simplicity where $k=1$ in LP  (\ref{eqn:multiILP}) and will discuss the general problem in the appendix. For $k=1$, the online formulation of LP  (\ref{eqn:multiILP}) reduces to a one-dimensional online LP problem,
\begin{align}
   \max \ \ & \sum_{t=1}^T r_t x_t \label{eqn:OLP}  \\
    \text{s.t. }\ & \sum_{t=1}^T \bm{a}_t x_t \le \bm{B} \nonumber  \\ 
    & 0 \le x_t \le 1,  \ \ t=1,...,T \nonumber
\end{align}
where $\bm{a}_t = (a_{1,t},...,a_{m,t})^\top \in \R^{m}$ and the decision variables are $\bm{x}=\left({x}_1,...,{x}_T\right)^\top\in \R^{T}.$ There are $m$ constraints and $T$ decision variables. Throughout the paper, we use $i$ to index constraints and $t$ to index decision variables. 

Now, we introduce our first group of assumptions on the distribution that governs the generation of $(r_t, \bm{a}_t)$'s. In the next section, we will introduce an additional assumption on the nondegeneracy of the underlying LP.

\begin{assumption}[Distribution] We assume
\begin{itemize}
    \item[(a)] Stochastic:  The column-coefficient pair $(r_t,\bm{a}_t)$'s are i.i.d. sampled from a distribution $\mathcal{P}.$ The distribution $\mathcal{P}$ takes a finite and known support $\{(\mu_j, \bm{c}_j)\}_{j=1}^n$ where $\mu_j \in \mathbb{R}$ and $\bm{c}_j \in \mathbb{R}^m$. Specifically, $$\prob((r_t, \bm{a}_t) = (\mu_j, \bm{c}_j)) = p_j$$ for $j=1,...,n$. The probability vector $\bm{p} = (p_1,...,p_n)^\top$ is unknown. 
    \item[(b)] Positiveness and Boundedness: $0\le\mu_j\le1$, $\bm{c}_j\ge \bm{0}$ and $\|\bm{c}_j\|_\infty\le 1$ for $j=1,...,n.$
    \item[(c)] Linear growth: The right-hand-side $\bm{B}=T\bm{b}$ for some $\bm{b}=(b_1,...,b_m)^\top > \bm{0}$.
\end{itemize}
\label{assp_dist}
\end{assumption}

Assumption \ref{assp_dist} (a) imposes a stochastic assumption for the customer orders. In addition, it states that the support of the distribution is finite and known, but the parameters are unknown. In other words, it means that there are $n$ known order types, and the order type at each time $t$ follows a multinomial distribution with unknown parameters. Assumption \ref{assp_dist} (b) requires all the entries of $(\mu_j,\bm{c}_j)$ between $0$ and $1$. We remark that all the results in this paper still hold (up to a constant) when this part is violated, and the positiveness and boundedness are introduced only for notation simplicity. Lastly, the linear growth condition in Assumption \ref{assp_dist} (c) is commonly assumed in problem setup and regret analysis online resource allocation problems. In our context, the condition is mild in that if $\bm{B} = o(T)$, we can always adjust the time horizon with $T'\ll T$ such that $\bm{B} = T'\bm{b}$, and consequently the linear growth condition holds for $T'$.

\subsection{Performance measure} 

In the literature of revenue management and the more general constrained online learning problems, a commonly considered performance benchmark is the \textit{certainty equivalent} version of the ``offline'' LP \eqref{eqn:OLP},
\begin{equation}
    \begin{aligned}
\text{OPT}_{\text{D}}  \coloneqq  \max \ \ & \sum_{j=1}^n p_j \mu_j y_j \label{eqn:DLP},  \\
    \text{s.t. }\ & \sum_{j=1}^n p_j \bm{c}_j \cdot y_j \le \bm{b}, \\ 
    & 0 \le y_j \le 1,  \ \ j=1,...,n.
    \end{aligned}
\end{equation}

Recall from Assumption \ref{assp_dist} that $\mu_j$ and $\bm{c}_j$ represent the revenue and requested resource consumption of the $j$-th order type, respectively. The right-hand-side $\bm{b}=\bm{B}/T$ represents the average resource capacity per time period, and $p_j$ is the probability of the $j$-th order type. The decision variables $y_j$'s prescribe a ``probabilistic'' decision rule for the orders, and $y_j$ can be interpreted as the proportion of accepted orders (or the probability of accepting orders) for the $j$-th order type. The connection between the LPs \eqref{eqn:OLP} and \eqref{eqn:DLP} can be seen from Assumption \ref{assp_dist} that the probability of the $j$-th order type is $p_j$. In this light, LP \eqref{eqn:DLP} can be viewed as a deterministic version obtained by taking expectation of the objective and the left-hand-side of LP \eqref{eqn:OLP}. For such reason, we refer to \eqref{eqn:DLP} as the deterministic LP (DLP). It is easy to verify that the optimal objective value of \eqref{eqn:DLP}, with a proper scaling of factor $T$, upper bounds the expected offline (hindsight) optimal objective value of \eqref{eqn:OLP}. Specifically, the expected offline (hindsight) optimal objective value is denoted by $\mathbb{E}[\text{OPT}_{\text{Hind}}]$, where $\text{OPT}_{\text{Hind}}$ is a random variable and is defined by
\begin{align}
   \text{OPT}_{\text{Hind}} \coloneqq \max \ \ & \sum_{t=1}^T r_t x_t \label{eqn:Hind}  \\
    \text{s.t. }\ & \sum_{t=1}^T \bm{a}_t x_t \le \bm{B} \nonumber  \\ 
    & 0 \le x_t \le 1,  \ \ t=1,...,T, \nonumber
\end{align}
where $(r_t, \bm{a}_t)$ follows the distribution in Assumption \ref{assp_dist}. It can be shown that $\mathbb{E}[\text{OPT}_{\text{Hind}}] \leq T\cdot\text{OPT}_{\text{D}}$. In this light, the DLP (fluid or certainty-equivalent) benchmark is a stronger one than the expected offline (hindsight) one. 

For the online problem, at each time $t$, we decide the value of $x_t$: $x_t=1$ means that we accept the order and allocate $\bm{a}_t$ amount of resources to this order accordingly; $x_t=0$ means that we reject the order. In this paper, we focus on the case of integer-valued solution, i.e., $x_t=0 \text{ or } 1$, but the analysis can be easily extended to the case where partial acceptance is allowed. Like the offline problem, we need to conform to the constraints throughout the procedure, i.e., no shorting of the resources is allowed. In this paper, we consider \textit{regret} as the performance measure, formally defined as follows:
$$\text{Reg}_{T}^\pi \coloneqq \E\left[T\cdot \text{OPT}_{\text{D}} - \sum_{t=1}^\tau r_tx_t\right]$$
where the quantity $\text{OPT}_{\text{D}}$ represents the optimal objective value of the DLP problem (\ref{eqn:DLP}) and $x_t$'s represent the online solution. Here $\tau$ is the stopping time for an algorithm when one or more types of the resource is depleted. The superscript $\pi$ denotes the online algorithm/policy according to which the online decisions are made. The expectation is taken with respect to $(r_t,\bm{a}_t)$'s and the (possible) randomness introduced by the algorithm.

\subsection{LP-based Adaptive Algorithm}

Before formally describing the algorithm, we first introduce a few additional notations to characterize the constraint consumption process. Define $\Bb_1=\Bb$ and $\Bb_t=(B_{1,t},...,B_{m,t})^\top$ as the remaining resource capacity at the beginning of time $t$, i.e.,
$$\bm{B}_t = \bm{B}_{t-1} - \bm{a}_{t-1} x_{t-1}.$$
Accordingly, we define $\Dd_t = \Bb_t/(T-t+1)$ as the average resource capacity for the remaining time periods.
In addition, we use $\bm{B}_{T+1}$ to denote the remaining constraint at the end of horizon, and the initial $\Dd_1 = (b_{1,1},...,b_{1,m})^\top = \Bb_1/T= \bm{B}/T=\bm{b}.$ We formally define the stopping time $\tau$ based on $\bm{B}_t$
$$\tau \coloneqq \min \{t: B_{i,t}\le 1 \text{ for some } i=1,...,m\} -1.$$
The rationale behind the definition is that when $B_{i,t}\le 1$, there may arrive an order that is profitable but cannot be fulfilled due to the resource constraint. On the opposite, when $B_{i,t} > 1$ for all $i$, we can fulfill any possible arriving order (given $\|\bm{c}_j\|_\infty\le 1$). Let $n_{j}(t)$ denote the counting process of the $j$-th order type, i.e., the number of observations $(\mu_j,\bm{c}_j)$ up to time $t$ (inclusively) for $j=1,...,n$. Since no shorting is allowed, i.e., the remaining constraint vector $\bm{B}_{t}$ must be element-wise non-negative for all $t=1,...,T$. Notice that the true probability distribution $\bm{p}=(p_1,...,p_n)$ is unknown. The counts $n_{j}(t)$'s will be used by the online algorithm to construct empirical estimates for the corresponding probabilities. 

\begin{algorithm}[ht!]
\caption{Adaptive Allocation Algorithm}
\label{alg:DAA}
\begin{algorithmic}[1]
\State Input: $\bm{B}, T, \{(\mu_j, \bm{c}_j)\}_{j=1}^n$
\State Initialize $\bm{B}_1 = \bm{b}$, $\bm{b}_1=\bm{B}_1/T$
\State Set $x_1=1$
\For {$t=2,..., T$}
\State Compute $\bm{B}_t=\bm{B}_{t-1}-\bm{a}_{t-1}x_{t-1}$
\State Compute $\bm{b}_t=\bm{B}_{t}/(T-t+1)$
\State Solve the following linear program where the decision variables are $(y_1,...,y_n)$:
\begin{align}
   \max \ \ & \sum_{j=1}^n \frac{n_{t-1}(j)}{t-1} \mu_j y_j \label{adpt_lp}  \\
    \text{s.t. }\ & \sum_{j=1}^n  \frac{n_{t-1}(j)}{t-1}\bm{c}_j\cdot y_j \le \bm{b}_t \nonumber  \\
    & 0 \le y_j \le 1,  \ \ j=1,...,n \nonumber
\end{align}
\State Denote the optimal solution as $\bm{y}_t^*=(y_{1,t}^*,...,y_{n,t}^*)$
\State Observe $(r_t,\bm{a}_t)$ and identify $(r_t,\bm{a}_t) = (\mu_j,\bm{c}_j)$ for some $j$
\State Set 
\begin{align*}
    x_t = \begin{cases}
    1, & \text{ with probability } y_{j,t}^*\\
    0, & \text{ with probability } 1-y_{j,t}^*
    \end{cases}
\end{align*}
when the constraint permits; otherwise set $x_t=0.$ 
\State Update the counts 
\begin{align*}
    n_{j}(t) = \begin{cases}
    n_{j}(t-1)+1, & \text{ if } (r_t,\bm{a}_t) = (\mu_j,\bm{c}_j)\\
    n_{j}(t-1), & \text{ otherwise }
    \end{cases}
\end{align*}
\EndFor
\State Output: $\bm{x} = (x_1,...,x_T)$
\end{algorithmic}
\end{algorithm}

Now we formally present the LP-based adaptive algorithm as Algorithm \ref{alg:DAA}. At each time $t$, the algorithm solves a sampled linear program \eqref{adpt_lp} to compute the probability of acceptance for each order type $(\mu_j, \bm{c}_j).$ The LP \eqref{adpt_lp} takes a similar form as LP \eqref{eqn:DLP} but differs in two aspects: (i) the probabilities $p_j$'s in \eqref{eqn:DLP} are replaced with their empirical estimates since the underlying distribution is assumed unknown; (ii) the right-hand-side $\bm{b}$ in \eqref{eqn:DLP} is replaced with its adaptive counterpart $\bm{b}_t$. Algorithm \ref{alg:DAA} then uses the LP's optimal solution $\bm{y}^*_t$ to determine the online solution $x_t$ at time $t$. Recall that $y^*_{j,t}$ denotes the optimal proportion of acceptance rate for the $j$-th order type for the deterministic LP associated with time $t$. Thus the probabilistic decision rule in Algorithm 1 aims to follow the prescription of the optimal solution by accepting the $j$-th order type with probability $y^*_{j,t}$. 

The algorithm is not new, and the adaptive design (using a dynamic right-hand-side) is commonly known as the re-solving technique in the network revenue management literature. Essentially, the algorithm has the same structure as the re-solving algorithms in \cite{jasin2012re}, \cite{jasin2015performance}, \cite{bumpensanti2020re}, and \cite{li2021online}. The algorithm here re-solves the problem in each time period, and it uses the sample counts as estimates for the true probabilities. In the following two sections, we analyze the performance of the algorithm for the cases when the underlying DLP \eqref{eqn:DLP} is (i) nondegenerate and (ii) degenerate. 

\section{Regret Analysis for Nondegenerate Case}
\label{Analysis_DLP}

The standard form of the DLP $\eqref{eqn:DLP}$ is as below. 
\begin{equation}
    \begin{aligned}
  \max \ \ & \bm{\mu}^\top \bm{y} \label{eqn:SDDLP}  \\
    \text{s.t. }\ & \bm{C}\bm{y} + \bm{s} = \bm{b} \text{ \ (dual variable: $\bm{\lambda}$)}  \\
    & \bm{y} + \bm{z} = \bm{1}\\
    & \bm{y}, \bm{s}, \bm{z} \geq \bm{0}.
    \end{aligned}
\end{equation}
where the decision variable vector $\bm{y}$ is the same as the one in $\eqref{eqn:DLP}$. With a slight abuse of notation (omitting the effect of the probability vector $\bm{p}$), we use $\bm{\mu}$ to denote the vector $(p_1\mu_1, ..., p_n\mu_n)^\top$ and $\bm{C}$ to denote the matrix $(p_1\bm{c}_1, ..., p_n\bm{c}_n)$. The additional decision variables $\bm{s}\in\mathbb{R}^m$ and $\bm{z}\in\mathbb{R}^n$ represent the slack variables for the corresponding constraints. The dual program for both the DLP \eqref{eqn:DLP} and its standard form \eqref{eqn:SDDLP} is
\begin{equation}\label{eqn:DDLP}
    \begin{aligned}
   \min \ \ & \bm{b}^\top{\bm{\lambda}} + \sum_{j=1}^n \gamma_j \\
    s.t. \ & p_j\bm{c}_j^\top{\bm{\lambda}} + \gamma_j \ge p_j\mu_j, \ \
    j=1,...,n \\ 
    & \bm{\lambda}\ge \bm{0}, \gamma_j \ge 0, \,\,j = 1,...,n
    \end{aligned}
\end{equation}
where the decision variables are $\bm{\lambda}$ and $\gamma_j$'s.

Denote the optimal solution to the LP \eqref{eqn:SDDLP} as $(\bm{y}^*, \bm{s}^*, \bm{z}^*)$.
and the dual optimal solution (also known as dual price) of \eqref{eqn:DDLP} for the resource constraints in \eqref{eqn:SDDLP} as $\bm{\lambda}^*$. Accordingly, we define the sets of \textit{basic} and \textit{non-basic} variables/order types as
$$\mathcal{J}^*\coloneqq \{j:\mu_j \ge\bm{c}_j^\top \bm{\lambda}^*,j=1,...,n\},\ \ \mathcal{J}'\coloneqq \{j:\mu_j<\bm{c}_j^\top \bm{\lambda}^*,j=1,...,n\},$$
and the sets of \textit{binding} and \textit{non-binding} constraints as
$$\mathcal{I}^*\coloneqq \{i:\bm{b}_i=\bm{C}_{i,:}^\top \bm{y}^* ,i=1,...,m\}, \ \ \mathcal{I}'\coloneqq \{i:\bm{b}_i>\bm{C}_{i,:}^\top \bm{y}^*,i=1,...,m\},$$
where $\bm{C}_{i,:}$ denotes the $i$-th row of the constraint coefficient matrix $\bm{C}$ in \eqref{eqn:SDDLP}. Here, the dual optimal solution $\bm{\lambda}^*$ provides a pricing rule for the resource consumption of a certain order.

Throughout this section, we assume a nondegeneracy structure for the standard form LP as below.

\begin{assumption}[Nondegeneracy]
\label{assp_nondeg}
The optimal solution to $\eqref{eqn:SDDLP}$ is unique and nondegenerate, i.e.,
$$|\{j: y_j^* \neq 0, j = 1, \cdots, n\}| + |\{i: s_i^* \neq 0, i= 1, \cdots, m\}| + |\{j: z_j^* \neq 0, j = 1, \cdots, n\}| = m+n.$$
\end{assumption}

The assumption is a standard one in the literature of linear programming, and with an arbitrarily small perturbation any LP can satisfy the assumption \citep{megiddo1989varepsilon}. An implication of Assumption \ref{assp_nondeg} is a stability structure for the underlying LP as follows. 

\begin{lemma}[Stability under nondegeneracy]
\label{lem_gen_stability}
Under Assumption \ref{assp_dist} and \ref{assp_nondeg}, there exists a positive constant $L$ which depends on $\bm{\mu},$ $\bm{C}$ and $\bm{b}$ such that if 
$$\max\{\|\hat{\bm{C}}-\bm{C}\|_\infty, \|\hat{\bm{\mu}}-\bm{\mu}\|_\infty, \max_{i\in\mathcal{I}^*}\{|b_i- \hat{b}_i|\}, \max_{i\in\mathcal{I}'}\{b_i- \hat{b}_i\}\} \le L,$$
then the following LP shares the same optimal basis and set of binding constraints with LP \eqref{eqn:SDDLP},
\begin{equation}
    \begin{aligned}
   \max \ \ & \hat{\bm{\mu}}^\top \bm{y} \label{eqn:SDDLP_1}  \\
    \text{s.t. }\ & \hat{\bm{C}}\bm{y} + \bm{s} = \hat{\bm{b}}  \\
    & \bm{y} + \bm{z} = \bm{1}\\
    & \bm{y}, \bm{s}, \bm{z} \geq \bm{0}.
    \end{aligned}
\end{equation}
\end{lemma}

From a geometric viewpoint, this lemma ensures that the optimal solutions of \eqref{eqn:SDDLP} and \eqref{eqn:SDDLP_1} coincide at the same corner point of the corresponding feasible simplex. We remark that for the non-binding constraints $i\in\mathcal{I}',$ it only needs a lower bound for $\hat{b}_i.$ For the online problem, the parameters $(\bm{\mu},\bm{C})$ are estimated through observations, and the resource level may deviate from the initial $\bm{b}.$ Thus the perturbed LP \eqref{eqn:SDDLP_1} is analogous to the adaptive LP \eqref{adpt_lp} used in the algorithm. In the previous analysis of the network revenue management problem \citep{jasin2012re, jasin2015performance}, the nondegeneracy assumption is used in a similar way to establish a stability for the underlying LP. We defer the proof of Lemma \ref{lem_gen_stability} to Appendix \ref{secLPStab} where we further relate the parameter $L$ explicitly with several key parameters of the underlying LP. The relationship refines the analysis in \cite{mangasarian1987lipschitz} and may be of independent interest. In the rest of the paper, we will express our regret bounds in terms of the parameter $L$ in Lemma \ref{lem_gen_stability}.

\subsection{Regret Decomposition}

The starting point of our analysis is to decompose the regret into three parts: (i) the first two parts concern the ``incorrect'' number of order acceptance; (ii) the third part concerns the remaining resources weighted by the dual price. We summarize the result in the following proposition, and defer the proof to Appendix \ref{ap_p1}.

\begin{proposition}\label{prop_general_regret_form}
Under Assumption \ref{assp_dist} and \ref{assp_nondeg}, the following equality holds 
\begin{equation}
    \begin{aligned}
    \mathrm{Reg}_{T}^\pi  = 
    \sum_{j\in \mathcal{J}^*} (\mu_j-\bm{c}_j^\top \bm{\lambda}^*)\cdot \E\left[n_j(T) -n_{j}^{a}(\tau) \right]+ \sum_{j\in \mathcal{J}'} (\bm{c}_j^\top \bm{\lambda}^*-\mu_j)\cdot \E\left[n_{j}^{a}(\tau) \right] + \bm{\lambda}^{*\top}\cdot \E\left[\bm{B}_{\tau}\right]
    \end{aligned}
    \label{regretUB}
\end{equation}
where $n_{j}^{a}(t)$ denotes the number of accepted orders of the $j$-th type up to time $t$ (inclusively) under policy $\pi$. Here $n_{j}(t)$ denotes the total number of orders of the $j$-th type up to time $t$ (inclusively) as defined earlier.
\end{proposition}

The equality has an intuitive interpretation. Recall that $\tau$ is the stopping time of the algorithm and the first time that some resource is (almost) depleted, and $\bm{B}_\tau$ denotes the remaining resource vector when the algorithm terminates. Thus the last part on the right-hand-side of \eqref{regretUB} penalizes the wasted resources when the process terminates. In particular, only residuals of the binding resources will be penalized. As to the first two parts on the right-hand-side of \eqref{regretUB}, we categorize the order types and elaborate as below:
\begin{itemize}
    \item All-accepted orders: $\bm{\mu}_j-\bm{c}_j^\top \bm{\lambda}^*>0$. For these orders, the optimal decision should be to accept all of them. We will observe $n_j(T)$ such orders throughout the horizon and aim to have the number of acceptance $n_{j}^a(\tau)$ close to that. 
    \item All-rejected orders:  $\bm{\mu}_j-\bm{c}_j^\top \bm{\lambda}^*<0$. On the opposite of the previous case, the optimal decision should be to reject all of these orders. Each acceptance of such order will induce a cost of $\bm{c}_j^\top \bm{\lambda}^*-\mu_j$: resources of value $\bm{c}_j^\top \bm{\lambda}^*$ are spent, but only reward of value $\mu_j$ is received.
    \item Partially-accepted orders: $\bm{\mu}_j-\bm{c}_j^\top \bm{\lambda}^*=0$. The condition may lead to a proportional acceptance of the orders, i.e., $0\le y_j^*\le 1$. Analysis-wise, there is no need to worry about these orders because they do not contribute to the first two terms of the right-hand-side of \eqref{regretUB}.
\end{itemize}

We make the following remarks for the regret decomposition. First, when the underlying probability distribution is unknown, the above categorization is also unknown a priori and should be learned. Second, in general, to identify the right proportion of acceptance for the partially-accepted order types is more challenging than the other two categories of orders. However, with the regret decomposition, we only need to focus on analyzing the constraint consumption $\bm{B}_{\tau}$ and avoid the complication related to the analysis of the partially-accepted orders. This is in contrast with the existing works \citep{jasin2012re, jasin2015performance, bumpensanti2020re} where it needs to carefully chase after (i) the number of the acceptance and (ii) the number of total arrivals of the partially-accepted order types. This paradigm shift from decision-variable-centric to constraint-centric is crucial in tightening the regret bound. 
Third, the nondegeneracy assumption makes the all-accepted orders and the all-rejected orders well separated from the remaining order types. Without this assumption, some order type's categorization (e.g. as either partially-accepted or all-rejected) can only be revealed in hindsight or at the very end of the online procedure \citep{bumpensanti2020re}. 


The following corollary extends Proposition \ref{prop_general_regret_form} to the case of a more general stopping time, and its proof can be found in Appendix \ref{ap_c1}.
\begin{corollary}
\label{cor_general_regret_form}
The following inequality holds 
\begin{equation}
    \begin{aligned}\label{ineq_regretUB}
    \mathrm{Reg}_{T}^\pi  & \leq 
    \sum_{j\in \mathcal{J}^*} (\mu_j-\bm{c}_j^\top \bm{\lambda}^*)\cdot  \E\left[n_j(\tau')-n_{j}^{a}(\tau') \right]+ \sum_{j\in \mathcal{J}'} (\bm{c}_j^\top \bm{\lambda}^*-\mu_j)\cdot \E\left[n_{j}^{a}(\tau') \right] \\ & \ \ \ \  \ \ \ +\left(T -\E[\tau']\right)\cdot \max_{j\in[n]}|\mu_j-\bm{c}_j^\top \bm{\lambda}^*|  + \bm{\lambda}^{*\top}\cdot \E\left[\bm{B}_{\tau'}\right]
    \end{aligned}
\end{equation}
where $\tau'$ is a stopping time adapted to the process $\bm{B}_t$'s and $\tau'\le\tau$ almost surely. 
\end{corollary}

The corollary replaces the algorithm termination time $\tau$ in Proposition \ref{prop_general_regret_form} with a general stopping time $\tau'$ and it includes an additional term $T-\E[\tau']$ measuring the closeness of $\tau'$ to the end of the horizon. This generalization gives much more flexibility in choosing a proper stopping time when analyzing the regret. In the following two subsections, we will analyze the terms in Corollary \ref{cor_general_regret_form} part by part.

\subsection{Order Acceptance -- First Two Terms in Regret Decomposition}

\label{secOrderAcc}

Next, we utilize Lemma \ref{lem_gen_stability} to analyze the first and second term in \eqref{ineq_regretUB}. Recall that in the setting of the online resource allocation problem, one passively collects the observations of orders. The sampled LP \eqref{adpt_lp} solved in Algorithm \ref{alg:DAA} will gradually converge to the DLP \eqref{eqn:DLP} if given the same right-hand-side. So, a concentration argument leads to a bound on the number of time periods it takes until the condition in Lemma \ref{lem_gen_stability} is met. Since then, the algorithm will make no further mistake on the acceptance (or rejection) of all-accepted orders (or all-rejected orders) defined by Proposition \ref{prop_general_regret_form}. One caveat is that part of the condition in Lemma \ref{lem_gen_stability} concerns the right-hand-side of the LP, so we also need to impose some restrictions on $\bm{b}_t$ -- the right-hand-side of the sampled LP \eqref{adpt_lp}. 

Specifically, we define a stopping time based on the constant $L$ in Lemma \ref{lem_gen_stability},
$$\tau_{S} \coloneqq \min \left\{t\le T: |b_{i,t}-b_i| > L  \text{ for some }i\in\mathcal{I}^* \right\} \cup \left\{t\le T: b_{i,t}-b_i > - L \text{ for some }i\in\mathcal{I}' \right\} \cup \{T+1\}.$$

Hence for any time $t<\tau_{S}$, the right-hand-side $\bm{b}_t$ meets the condition in Lemma \ref{lem_gen_stability}. Thus the underlying adaptive LP shares the same structure as the DLP \eqref{eqn:DLP} when $t\le \tau_S$. This property makes the stopping time $\tau_S$ easier to analyze than the original stopping time $\tau.$ Because when $t$ approaches to $\tau,$ the underlying LP's optimality and bindingness structure may already change; but this will not happen for $t\le \tau_S$. 

By the definition of $\tau_{S},$ if we restrict our attention to time periods before $\tau_S$, the numbers of mistakes made on the all-accepted and all-rejected orders are purely caused by the inaccurate estimation on the left-hand-side. To obtain an upper bound, the concentration argument can be applied and it leads to the following proposition. 

\begin{proposition}
\label{prop_second_third_term}
Under Assumptions \ref{assp_dist} and \ref{assp_nondeg}, the output of Algorithm \ref{alg:DAA} satisfies
\begin{equation*}
    \begin{aligned}
    &\hspace{5mm}\sum_{j\in \mathcal{J}^*} (\mu_j-\bm{c}_j^\top \bm{\lambda}^*)\cdot \mathbb{E}\left[n(\tau_S)  -n_{j}^{a}(\tau _S)\right]+ \sum_{j\in \mathcal{J}'} (\bm{c}_j^\top \bm{\lambda}^*-\mu_j)\cdot \mathbb{E}\left[n_{j}^{a}(\tau _S)\right]\leq   \frac{2n \max_j \left|\mu_j - \bm{c}_j^\top\bm{\lambda}^*\right|}{1-\exp(-2L^2)}.
    \end{aligned}
\end{equation*}
\end{proposition}

We defer the proof to Appendix \ref{ap_p2}. Note that the right-hand-side of the above inequality is not dependent on the time horizon $T$. The proposition's proof directly follows an application of Lemma \ref{lem_gen_stability}. First, as mentioned earlier, the definition of $\tau_S$ precludes the possibility that the changing right-hand-size $\bm{b}_t$ in Algorithm \ref{alg:DAA} affects the LP's stability. Then, with a sufficient number of observations, the condition in Lemma \ref{lem_gen_stability} can be satisfied with high probability. From that time on, the sampled LP \eqref{adpt_lp} and the deterministic LP \eqref{eqn:DLP} will share the same optimal basis, and consequently the algorithm will not make any further mistake on the all-accepted and all-rejected orders.

If we put together Proposition \ref{prop_second_third_term} with Corollary \ref{cor_general_regret_form}, the remaining task is to bound the remaining time periods $\E[T-\tau_S]$ and the remaining resources $\E[\bm{B}_{\tau_S}].$ In fact, these two aspects are closely related with each other and they will be the focus of the next subsection.

\subsection{Constraint Consumption Process -- Last Two Terms in Regret Decomposition}
\label{sec_reg_binding}

In this section, we analyze the constraint consumption process $\bm{B}_t$ (equivalently, $\bm{b}_t$) and handle the last two terms in the generic upper bound \eqref{ineq_regretUB}. To better illustrate the proof idea, we focus in this section on the case that all the constraints are binding. We will show in Appendix \ref{sec_reg_general} how the analysis can be adapted to the case where there are both binding and non-binding constraints.

\begin{assumption}
All the resource constraints are binding, i.e., $\bm{s}^*=\bm{0}$ in \eqref{eqn:SDDLP}.
\label{assp_binding_new}
\end{assumption}

From Lemma \ref{lem_gen_stability}, we know that under Assumptions \ref{assp_dist} and \ref{assp_nondeg}, the LP's optimality and bindingness structure remains to hold when $\bm{b}$ is perturbed. Let
$$\mathfrak{B} \coloneqq \bigotimes_{i=1}^m [b_i-L, b_i+L].$$
where $L$ is the constant in Lemma \ref{lem_gen_stability}. The following lemma states that for any $\tilde{\bm{b}}\in \mathfrak{B}$, all the constraints of the corresponding LP are binding.

\begin{lemma}
\label{lem_stability_b}
Under Assumption \ref{assp_dist}, \ref{assp_nondeg}, and \ref{assp_binding_new}, for each $\tilde{\bm{b}}\in \mathfrak{B}$, there exists an optimal solution $\tilde{\bm{y}}^*=(\tilde{y}_1^*,....,\tilde{y}_n^*)^\top$ of the DLP \eqref{eqn:DLP} with the right-hand-side being $\tilde{b}$ that satisfies
$$ \sum_{j=1}^n p_j\bm{c}_j \cdot \tilde{y}_j^* = \tilde{\bm{b}}.$$
\end{lemma}

Our goal of analyzing $\bm{B}_t$ is to bound (i) the remaining time periods when the algorithm terminates, $T-\E[\tau_S]$, and (ii) the resource left-over, $\bm{\lambda}^{*\top}\E[\bm{B}_{\tau_S}]$. Throughout the analysis, we will reserve the notation $\bm{b}=(b_1,...,b_m)^\top$ for the initial average resource capacity and use $\bm{b}'$ to denote an arbitrary value in $\mathbb{R}^m.$ Ideally, the process $\bm{b}_t$ should stay near $\bm{b}$ throughout the horizon, as this would imply that the resource is exhausted only at the very end of the horizon. Let $x_t(\bm{b}')$ denote the online solution output by Algorithm \ref{alg:DAA} at the $t$-th time period as a function of the input $\bm{b}_t = \bm{b}'$. Consider the following event defined in the space of the history observations up to time $t-1,$
$$\mathcal{E}_t \coloneqq \left\{\mathcal{H}_{t-1}\Big \vert \sup_{\bm{b}'\in \mathfrak{B}}\left\|\mathbb{E}[\bm{a}_{t}x_{t}(\bm{b}')|\mathcal{H}_{t-1}, \bm{b}_t = \bm{b}']-\bm{b}'\right\|_\infty \le \epsilon_{t-1} \right\}$$
where the history $\mathcal{H}_{t-1}=(r_1,\Aa_1,...,r_{t-1},\Aa_{t-1})$. Here we choose
\begin{equation}\label{eqn:epsilon_t}
    \begin{aligned}
    \epsilon_t \coloneqq \begin{cases}
1 & t \leq \kappa T, \\
t^{-\frac{1}{4}} & t >  \kappa T,
\end{cases}
\end{aligned}
\end{equation}
where the constant $\kappa\in(0,1)$ is to be specified later and will be roughly on the same order of $L$. Without loss of generality, we assume $\kappa T$ takes an integer value. 

Now, we provide some intuitions for the definition of $\mathcal{E}_t.$ First, recall that
$$\bm{b}_{t+1}=\frac{\bm{B}_{t+1}}{T-t}= \frac{\bm{B}_{t}-\bm{a}_tx_t}{T-t}= \bm{b}_t -\frac{1}{T-t}(\bm{a}_tx_t-\bm{b}_t)$$
for $t=1,...,T-1.$ The definition of $\mathcal{E}_t$ is aligned with the hope that the expected resource consumption at each time $t$ stays close to $\bm{b}_t$. Specifically, the event $\mathcal{E}_t$ controls the expectation of $\bm{a}_tx_t-\bm{b}_t$, and in its definition, the supremum taken over $\bm{b}'\in \mathfrak{B}$ is necessary because $\bm{b}_t$ is random. As to the choice of $\epsilon_t$:
\begin{itemize}
\item $\epsilon_t$ should not be too small so that the events $\mathcal{E}_t$'s will happen with high probability.
\item $\epsilon_t$ should not be too large so that conditional on $\mathcal{E}_t$'s, the process $\bm{b}_t$ is ``stable''.
\end{itemize}



To formalize the intuitions, we define a stopping time to capture the ``bad'' event that is either $\bm{b}_t\notin \mathfrak{B}$ or the complement of the event $\mathcal{E}_t$,
$$\tilde{\tau} \coloneqq \min \{t\le T:\bm{b}_t \notin \mathfrak{B} \text{ or } \mathcal{H}_{t-1} \notin \mathcal{E}_t\} \cup \{T+1\}.$$
By comparing the definitions of $\tau_S$ and $\tilde{\tau},$ we claim that $\tilde{\tau}=\tau_S$ with high probability (validated in Lemma \ref{lem_part2decompose}). The sample paths that render $\tilde{\tau}\neq \tau_S$ fall into the event $\bar{\mathcal{E}}_t=\{\mathcal{H}_{t-1} \notin \mathcal{E}_t\}$, i.e., the expected constraint consumption in a single time period has a large deviation from zero. The stopping time $\tilde{\tau}$ chops off such ``bad'' event and through the lens of $\tau_S$, and the constraint process becomes easier to analyze.

Specifically, with $\tilde{\tau}$, we define an auxiliary process $\tilde{\bm{b}}_t$ as follows 
$$\tilde{\bm{b}}_t = \begin{cases} 
\bm{b}_t, & t<\tilde{\tau},\\
\bm{b}_{\tilde{\tau}}, & t\ge \tilde{\tau}.
\end{cases}$$
By its definition, the process $\tilde{\bm{b}}_t$ freezes its value once $\bm{b}_t$ exits the region $\mathfrak{B}$ or the bad event $\bar{\mathcal{E}}_t$ happens. With $\tilde{\bm{b}}_t$, we have
\begin{align}
    \prob\left(\bm{b}_s\notin\mathfrak{B} \text{ for some }s\le t \right)&= \prob\left(\bm{b}_s\notin\mathfrak{B} \text{ for some }s\le t, \cap_{s=1}^t \mathcal{E}_s \right) + \prob\left(\bm{b}_s\notin\mathfrak{B} \text{ for some }s\le t, \cup_{s=1}^t \bar{\mathcal{E}}_s \right) \nonumber \\
    &\le \prob\left(\tilde{\bm{b}}_s\notin\mathfrak{B} \text{ for some }s\le t \right) + \sum_{s=1}^t \prob(\bar{\mathcal{E}}_s) \label{decompose}
\end{align} 
where $\bar{\mathcal{E}}$ denotes the complement of an event $\mathcal{E}.$
For the first part of the second line, it is because given $\cap_{s=1}^t \mathcal{E}_s$, the event that $\bm{b}_s \notin \mathfrak{B}$ for some $s\le t$ is equivalent to the event $\tilde{\tau}\le t$ and thus it entails $\tilde{\bm{b}}_s \notin \mathfrak{B}$.  For the second part, it is obtained by ignoring the condition on $\bm{b}_s$ and then taking a union bound with respect to $s=1,...,t$.

We justify the decomposition of the left-hand-side of \eqref{decompose} by relating it with the last two terms in Corollary \ref{cor_general_regret_form} as the following lemma.
\begin{lemma}\label{lem_bound_cor_last_two_terms}
We have the following relation between the stopping time and the left-hand-side of \eqref{decompose}
\begin{align*}
    \mathbb{E}[T-\tau_S] \le  \sum_{t=1}^T\prob\left(\bm{b}_s\notin\mathfrak{B} \text{ for some }s\le t \right),
\end{align*}
and for the last two terms in Corollary \ref{cor_general_regret_form}, we have
\begin{align*}
    \left(T -\E[\tau_S]\right)\cdot \max_{j\in[n]}|\mu_j-\bm{c}_j^\top \bm{\lambda}^*|  + \bm{\lambda}^{*\top}\cdot \E\left[\bm{B}_{\tau_S}\right] \leq ||\bm{\lambda^*}||_1\cdot(3\mathbb{E}[T-\tau_S]+4).
\end{align*}
\end{lemma}

Now we discuss the motivation for defining $\tilde{\bm{b}}_t$ and why the inequality (\ref{decompose}) is useful. The inequality (\ref{decompose}) separates our goal -- the probability on the left-hand-side -- into two components. The first component concerns the process $\tilde{\bm{b}}_t$ which is a relatively ``well-behaved'' process in that when $t<\tilde{\tau}$, the process $\tilde{\bm{b}}_t$'s fluctuation is subject to the event $\mathcal{E}_t$; when $t\ge \tilde{\tau}$, the process freezes. The event $\mathcal{E}_t$ further controls the fluctuation of the process $\tilde{\bm{b}}_t$, and thus the process behaves roughly like a martingale. The second component concerns the probability of $\bar{\mathcal{E}}_t$'s, which can be analyzed individually for each $t$. Overall, the inequality \eqref{decompose} disentangles the stability of the process $\bm{b}_t$ from the estimation error. Its first component concerns the process stability given a good estimate (the event $\mathcal{E}_t$), while its second component concerns the probability of obtaining the good estimate for the model parameters. Piecing the two components together, we obtain a bound on the probability that $\bm{b}_t$ exists the region $\mathfrak{B}.$ 

Next, we will analyze the two components in \eqref{decompose} separately and then combine the results to derive the regret bound.

\paragraph{\textbf{Analysis of the first component in (\ref{decompose}).}} \ 

The following theorem states a concentration result for a general martingale difference sequence $X_t$'s. The approach that we analyze the first component can be viewed as a two-step procedure: we first construct a martingale as an approximation of the process $\tilde{\bm{b}}_t$. Specifically, the constructed martingale and $\tilde{\bm{b}}_t$ share the same initial value $\bm{b}$, and the difference between the martingale and $\tilde{\bm{b}}_t$ is controllably small. Second, we apply Theorem \ref{thm_dehling} for the constructed martingale and argue that both the martingale and $\tilde{\bm{b}}_t$ will stay within $\mathfrak{B}$ with high probability.
 
\begin{theorem}{(Hoeffding's inequality for dependent data \citep{van2002hoeffding})} \label{thm_dehling}
Consider a sequence of random variables $\{X_t\}_{t=1}^T$ adapted to the filtration $\mathcal{F}_t$'s and 
$$\E[X_t | \mathcal{F}_{t-1}] = 0 \text{ \ for \ } t=1,...,T$$
where $\mathcal{F}_0=\emptyset.$
Suppose $L_t, U_t$ are $\mathcal{F}_{t-1}$-measurable random variables such that $L_t\le X_t\le U_t$ almost surely for $t=1,...,T$. Let $S_t =\sum_{s=1}^t X_t$ and $V_t=\sum_{s=1}^t(U_s-L_s)^2$. Then, the following inequality holds for any $b>0, c>0$ and $T \in \mathbb{N}_+,$
\begin{align*}
\prob(|S_t|\ge b, V_t \le c^2 \text{ for some } t\in\{1,...,T\}) \leq 2e^{-\frac{2b^2}{c^2}}.
\end{align*}
\label{Hf_dependent}
\end{theorem}

The following lemma utilizes the result in Theorem \ref{Hf_dependent} and provides an upper bound on the first component in \eqref{decompose}. Thus it completes the analysis of the first component in (\ref{decompose}).  The proof can be found in Appendix \ref{ap_l3}.

\begin{lemma}
\label{lem_part1decompose}
For $T\ge T_1$ and $t\le T-2$, the following inequality
\begin{equation*}
    \begin{aligned}
    \prob\left(\tilde{\bm{b}}_s\notin \bigotimes_{i=1}^m \left[b_{i}-\Delta, b_{i}+\Delta\right] \text{ for some }s\le t\right)\le 2me^{-\frac{\Delta^2(T-t)}{8}}
    \end{aligned}
\end{equation*}
holds for any $\Delta>0$. The constant $T_1$ is defined as the minimal integer such that $T_1 \ge \frac{1}{\exp{\left(\frac{\Delta}{8}\right) - 1}}  + 2$ and 
$\frac{\log T_1}{T_1^{1/4}} \le \frac{\kappa^{1/4}\Delta}{4}.$
For the parameter in the definition of $\epsilon_t$ in \eqref{eqn:epsilon_t}, we set $\kappa=1-\exp(-\frac{\Delta}{8})$.
\end{lemma}

From the lemma, if we set $\Delta=L$ and $\kappa=1-\exp(-\frac{L}{8})$, we have
\begin{equation*}
    \begin{aligned}
    \prob\left(\tilde{\bm{b}}_s\notin \mathfrak{B} \text{ for some }s\le t\right)\le 2me^{-\frac{L^2(T-t)}{8}}.
    \end{aligned}
\end{equation*}
Thus we obtain a bound for the first component in (\ref{decompose}).

\paragraph{\textbf{Analysis of the second component in (\ref{decompose}).}} \

Define the events for $t=2,...,T$
\begin{equation*}\label{event_At}
    \begin{aligned}
    \mathcal{A}_t^{(j)} \coloneqq \left\{ \left|\frac{n_j(t-1)}{t-1} - p_j\right| \leq L   \right\}
    \end{aligned}
\end{equation*}
and 
$$\mathcal{B}_t^{(j)} \coloneqq \left\{ \left\vert \frac{n_j(t-1)}{t-1}-p_j\right\vert \leq \frac{1}{n(t-1)^{1/4}} \right\}.$$
In addition, we define $\mathcal{A}_1^{(j)}=\mathcal{B}_1^{(j)}=\Omega.$

Now we argue that 
$$\left(\cap_{j=1}^n \mathcal{A}_{t}^{(j)}\right)\cap\left(\cap_{j=1}^n \mathcal{B}_{t}^{(j)}\right) \subseteq \mathcal{E}_t.$$

To see this, for $t \geq 2$, let $\bm{y}^*_{t}=(y_{1,t}^*,...,y_{n,t}^*)^\top$ be the optimal solution of \eqref{adpt_lp} with $\bm{b}_{t}=\bm{b}'$ for some $\bm{b}'\in \mathfrak{B}$. By the algorithm, we have the expected resource consumption at time $t$
$$\E[\bm{a}_{t}x_{t}(\bm{b}')|\mathcal{H}_{t-1},\bm{b}_{t}=\bm{b}'] = \sum_{j=1}^n p_j\bm{c}_j\cdot y_{j,t}^*.$$
Moreover, we know that given the event $\cap_{j=1}^n \mathcal{A}_{t}^{(j)},$ the perturbation of $\bm{C}$ and $\bm{\mu}$ (described in Lemma \ref{lem_gen_stability}) will be within $L$. From Lemma \ref{lem_gen_stability} and Assumption \ref{assp_binding_new} we know that all the constraints of the LP with right-hand-side being $\bm{b}'$ are binding
$$\bm{b}' = \sum_{j=1}^n \frac{n_j(t-1)}{t-1} \bm{c}_j \cdot y_{j,t}^*.$$
Then, taking the difference, 
\begin{equation}\label{eqn_axminusd}
    \begin{aligned}
    \E[\bm{a}_{t}x_{t}(\bm{b}')|\mathcal{H}_{t-1},\bm{b}_{t}=\bm{b}']-\bm{b}' = \sum_{j=1}^n \left(p_j - \frac{n_j(t-1)}{t-1}\right)\bm{c}_j\cdot y_{j,t}^*.
    \end{aligned}
\end{equation}
Next, given the event $\cap_{j=1}^n\mathcal{B}_{t}^{(j)},$ we have
\begin{align*}
 \left\|\E[\bm{a}_{t}x_{t}(\bm{b}')|\mathcal{H}_{t-1},\bm{b}_{t}=\bm{b}']-\bm{b}'\right\|_\infty & = \left\|\sum_{j=1}^n \left(p_j - \frac{n_j(t-1)}{t-1}\right)\bm{c}_j y_{j,t}^*\right\|_\infty\\
 & \le \sum_{j=1}^n \left \vert p_j - \frac{n_j(t-1)}{t-1}\right\vert \\
 & \le \min\left\{(t-1)^{-\frac{1}{4}}, 1\right\},
\end{align*}
where we use the fact that $\|\bm{c}_j\|_\infty\le 1$ from Assumption \ref{assp_dist}. This meets the definition of the event $\mathcal{E}_t$ and this result is summarized in Lemma \ref{lem_part2decompose}. To analyze the events $\mathcal{A}_{t}^{(j)}$ and $\mathcal{B}_t^{(j)}$, we can simply apply the concentration argument. In this way, we complete our analysis of the second component in \eqref{decompose}.
\begin{lemma}\label{lem_part2decompose}
We have 
$\left(\cap_{j=1}^n \mathcal{A}_{t}^{(j)}\right)\cap\left(\cap_{j=1}^n \mathcal{B}_{t}^{(j)}\right) \subseteq \mathcal{E}_t$
and the following inequality holds for each $t=1,...,T,$
$$\prob\left(\bar{\mathcal{E}_t}\right) \leq 2n \exp\left(-2L^2(t-1)\right) +2n\exp\left(-\frac{2(t-1)^{1/2}}{n^2}\right).$$
\end{lemma}
The proof is left in Appendix \ref{ap_l4}.
\subsection{Final Regret Bound}
We can derive the final regret bound by combining Corollary \ref{cor_general_regret_form}, Proposition \ref{prop_second_third_term}, Lemma \ref{lem_bound_cor_last_two_terms}, \ref{lem_part1decompose} and \ref{lem_part2decompose}. The following theorem states that the regret is uniformly bounded in terms of $T$. We defer the detailed proof to Appendix \ref{ap_t2}. The regret bound is related with parameters such as the number of customer/order types $n$, the number of constraints $m$ and the stability parameter $L$ (in Lemma \ref{lem_gen_stability}). In this sense, the regret should be interpreted as a problem-dependent bound rather than a worst-case bound. The implication is that when the underlying LP is well-posed such that the parameter $L$ can be treated as a constant, the algorithm's regret does not scale up with the number of time horizon $T$. We also remark that the theorem requires Assumption \ref{assp_binding_new} where all the constraints are binding. We make the assumption to better illuminate the analysis of the underlying constraint process. In Appendix \ref{sec_reg_general}, we will remove the assumption and study the general setting where there exist both binding and non-binding constraints. We note that the proof idea and the regret bound for the general setting without Assumption \ref{assp_binding_new} are the same as the case of all binding constraints in this section. 

\begin{theorem}\label{theorem_regret}
Under Assumptions \ref{assp_dist}, \ref{assp_nondeg} and \ref{assp_binding_new}, the regret of Algorithm \ref{alg:DAA} satisfies
\begin{align*}
    \mathrm{Reg}_{T}^\pi
     \le \frac{(48m+4n+12)\cdot \|\bm{\lambda}^*\|_1}{L^2} + o(1).
\end{align*}
where $\pi$ denotes the policy specified by Algorithm \ref{alg:DAA}, $\bm{\lambda}^*$ is the dual optimal solution of LP \eqref{eqn:DDLP}, and the last term $o(1)\rightarrow 0$ as $T\rightarrow \infty.$
\end{theorem}

\section{Regret Analysis for Degenerate Case}\label{sec_degenerate}

Now we analyze the algorithm's performance without the nondegeneracy assumption (Assumption \ref{assp_nondeg}). As noted by \cite{bumpensanti2020re}, the existing results on bounded regret \citep{jasin2012re, wu2015algorithms} require the nondegeneracy assumption. The regret bounds therein, though bear no dependency in $T$, will become arbitrarily large when the LP is nondegenerate but nearly degenerate. In this section, we present a positive result on the performance of Algorithm \ref{alg:DAA} for the degenerate case. Specifically, we show that Algorithm \ref{alg:DAA} achieves a $O(\sqrt{T}\log T)$ regret without the nondegeneracy assumption. The bound is nearly tight in that $\Omega(\sqrt{T})$ lower bound for Algorithm \ref{alg:DAA} is established by \cite{bumpensanti2020re} against both the fluid benchmark OPT$_{\text{D}}$ and hindsight benchmark.

We first provide some intuition of the analysis. Consider the following three LPs:
\begin{align*}
R_t \coloneqq \max \ \ & \bm{\mu}_t^\top \bm{y}  & \bar{R}_t \coloneqq \max \ \ & \bm{\mu}^\top \bm{y} & \text{OPT}_{\text{D}} \coloneqq \max \ \ & \bm{\mu}^\top \bm{y}  \\
\text{s.t. }\ & \bm{C}_t\bm{y} \leq \bm{b}_t & \text{s.t. }\ & \bm{C}\bm{y} \leq \bm{b}_t & \text{s.t. }\ & \bm{C}\bm{y} \leq \bm{b}  \\
& \bm{0} \leq \bm{y} \leq \bm{1}, & & \bm{0} \leq \bm{y} \leq \bm{1}, & & \bm{0} \leq \bm{y} \leq \bm{1}, 
\end{align*}
where  
\begin{equation*}
    \begin{aligned}
     \bm{\mu}_t = \left(\frac{n_1(t-1)}{t-1}\mu_1, \cdots , \frac{n_n(t-1)}{t-1}\mu_n\right)^\top, & \ \ 
     \bm{C}_t = \left(\frac{n_1(t-1)}{t-1}\bm{c}_1, \cdots , \frac{n_n(t-1)}{t-1}\bm{c}_n\right),
     \\
    \bm{\mu} = \left(p_1\mu_1, \cdots , p_n\mu_n\right)^\top, & \ \ 
     \bm{C} = \left(p_1\bm{c}_1, \cdots , p_n\bm{c}_n\right).
    \end{aligned}
\end{equation*}
We note that the left LP in above is the sampled LP \eqref{adpt_lp} used in Algorithm \ref{alg:DAA}, while the right LP is the deterministic LP \eqref{eqn:DLP} with optimal objective value OPT$_{\text{D}}$. It is easy to see that $R_t$ is the expected reward (conditional on the history) that Algorithm \ref{alg:DAA} collects at time $t$. Then the single-period regret at time $t$ follows
$$\text{OPT}_{\text{D}} - R_t = (\text{OPT}_{\text{D}} - \bar{R}_t) + (\bar{R}_t-R_t).$$ 
To analyze the right-hand-side, we define for $2\le t \le T-1,$ 
\begin{equation*}
    \begin{aligned}
    \mathcal{C}_t^{(i)} &= \left\{b_{i,t} > b_i  -\frac{\sqrt{4n\log 2T}}{\sqrt{t}}-\frac{\sqrt{4n\log 2T} + \sqrt{2\log 2T}}{\sqrt{T-t}}\right\}, \\
    \mathcal{D}_t^{(j)} &= \left\{\left\vert\frac{n_j(t-1)}{(t-1)p_j} - 1\right\vert < \frac{\sqrt{\log 2T}}{\sqrt{2\underline{p}^2(t-1)}}\right\}, \\
    \end{aligned}
\end{equation*}
where $\underline{p} = \min\{p_1,....,p_n\}.$

To motivate the definition of these events, first, given the events $\mathcal{C}_t^{(i)}$'s, $\bm{b}_t$ is close to $\bm{b}$, and consequently, $\text{OPT}_{\text{D}} - \bar{R}_t$ is small. Second, given the events $\mathcal{D}_t^{(j)}$'s, the left LP and the middle LP in above are close to each other, and thus, $\bar{R}_t-R_t$ is small.

The following lemma establishes that the above events will happen with high probability. The analysis of the event $\mathcal{C}_{t}^{(i)}$ essentially reduces to the analysis of the constraint process, which is in a similar spirit as that of the previous section. The analysis of the event $\mathcal{D}_{t}^{(j)}$ is simply based on a concentration argument. The proof can be found in Appendix \ref{ap_l5}.

\begin{lemma} \label{lem_deg_events}
We have
$$\prob\left(\cap_{i=1}^{m} \mathcal{C}_{t}^{(i)}\right) \ge 1 - \frac{m}{T}$$
and 
$$\prob\left(\cap_{j=1}^{n} \mathcal{D}_{t}^{(j)}\right) \ge 1 - \frac{n}{T}.$$
\end{lemma}
Denote $\mathcal{C}_t = \cap_{i=1}^m \mathcal{C}_t^{(i)}$ and $\mathcal{D}_t = \cap_{j=1}^n \mathcal{D}_t^{(j)}$. We formalize the intuition above into the following lemma, and leave the proof in Appendix \ref{ap_l6}.
\begin{lemma}
\label{lemma_single_step}
Under event $\mathcal{C}_t\cap\mathcal{D}_t$, we have the single-period regret
\begin{equation*}\label{lem_purturb_inter}
    \begin{aligned}
    \mathrm{OPT}_{\mathrm{D}} - R_t \leq \max\{1, \bar{\lambda}\}\cdot \left(m\left(\frac{\sqrt{4n\log 2T}}{\sqrt{t}}+\frac{\sqrt{4n\log 2T} + \sqrt{2\log 2T}}{\sqrt{T-t}}\right) + n\frac{\sqrt{\log 2T}}{\sqrt{2\underline{p}^2(t-1)}}\right),
    \end{aligned}
\end{equation*}
where $\bar{\lambda} := \max \left\{||\bm{\lambda}||_{\infty} : \bm{\lambda} \in \mathcal{FD}_0 \right\}$ and $\mathcal{FD}_0$ denotes the set of basic feasible solutions for the dual of DLP \eqref{eqn:DLP}.
\end{lemma}

Combining the above two lemmas, we yield the regret bound without the nondegeneracy as follows. The theorem provides a regret bound sublinear in $T$ for Algorithm \ref{alg:DAA} without the nondegeneracy assumption. We remark that the analysis also covers the case when the distribution is known: in that case, the algorithm will use the LP in the middle $\bar{R}_t$ to guide the online decision, and thus the event $\mathcal{D}_t^{(j)}$ will happen with probability 1. In contrast, the previous bounds in \cite{jasin2012re} and \cite{wu2015algorithms} will be arbitrarily large when the underlying LP approaches a nondegenerate one. The key difference between our analysis and the previous analyses is the focus on the constraint process. On one hand, the definition of the events $\mathcal{C}_t^{(i)}$ imposes a condition under which a sublinear regret is achievable. On the other hand, martingale tools enable a careful analysis of the constraint process to meet such condition.

\begin{theorem}\label{thm_degenerate}
Under Assumption \ref{assp_dist}, Algorithm \ref{alg:DAA} gives a regret upper bounded by
\begin{equation*}
    \begin{aligned}
    \mathrm{Reg}_{T}^\pi & \leq  \max\{1, \bar{\lambda}\}\cdot \left(m\left(\sqrt{2}+\sqrt{16n}\right)+\frac{n}{\sqrt{2\underline{p}^2}}\right)\sqrt{T}\sqrt{\log 2T} + 1 + n + m \\
    & = O((m\sqrt{n}+n)\sqrt{T}\log T).\\
    \end{aligned}
\end{equation*}
\end{theorem}

The proof of the theorem is deferred to Appendix \ref{ap_t3}. Apart from the $\Omega(\sqrt{T})$ lower bound result in \cite{bumpensanti2020re}, we can also understand from the analysis the reason why bounded regret cannot be achieved for the degenerate case. In the analysis for the degenerate case, the regret decomposition in Proposition \ref{prop_general_regret_form} no longer holds and thus we have to rely on the single-period regret bound in Lemma \ref{lemma_single_step}. The single-period regret bound is looser in a sense that it treats all the time period separately. Specifically, in the nondegenerate case, if the algorithm performs poorly in one time period, it may recover with better reward later through the adaptive (re-solving) mechanism. But this recovering mechanism cannot be captured when we treat all the time periods separately.

\section{Numerical Experiment and Discussions}

We conclude with numerical experiments to illustrate our analysis. We perform our simulation experiments under both degenerate and nondegenerate cases, and the results are consistent with our theoretical findings. Specifically, we consider the following problem instance. There are two types of resources and three types of customer orders. The unknown probability vector of three order types is $(p_1,p_2,p_3) = (0.3, 0.3, 0.4)$, the reward vector is $(\mu_1,\mu_2,\mu_3) = (1, 1.2, 0.8)$, and the resource consumption is $\bm{c}_1= (1,2)^\top,\bm{c}_2 = (2,1)^\top,\bm{c}_3= (1,1)^\top.$
Thus the underlying LP is 
\begin{equation}
    \begin{aligned}
   \max \ \ & 0.3y_1 + 0.36y_2 + 0.32y_3 \\
    \text{s.t. }\ & 0.3y_1 + 0.6y_2 + 0.4y_3 \leq b_1\\ 
    & 0.6y_1 + 0.3y_2 + 0.4y_3 \leq b_2\\ 
    & 0 \le y_j \le 1,  \ \ j=1,2,3.
    \end{aligned}
\end{equation}
For a nondegenerate problem instance, we set $(b_1, b_2) = (1, 1)$; for a degenerate problem instance, we set $(b_1, b_2) = (1, 1.15)$.

\begin{figure}[ht!]
\begin{subfigure}{.5\textwidth}
  \centering
  \includegraphics[width=.9\linewidth]{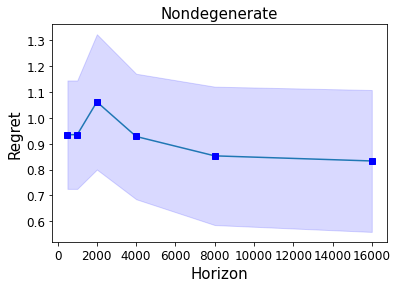}
\end{subfigure}%
\begin{subfigure}{.5\textwidth}
  \centering
  \includegraphics[width=.9\linewidth]{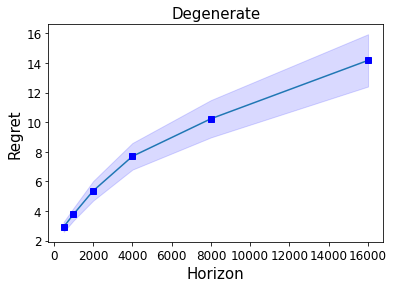}
\end{subfigure}
\caption{Regret of Algorithm \ref{alg:DAA} under different horizon length.}
\label{fig1}
\end{figure}

Figure \ref{fig1} describes the relationship between the regret and the horizon length $T$ under both nondegenerate and degenerate cases. It reports both the sample mean and the 99\%-confidence interval (in a light color). When the underlying LP is nondegenerate, the regret does not scale with $T$, but when the LP is degenerate, the regret grows on the order of $O(\sqrt{T})$.

\begin{figure}[ht!]
\begin{subfigure}{.33\textwidth}
  \centering
  \includegraphics[width=.99\linewidth]{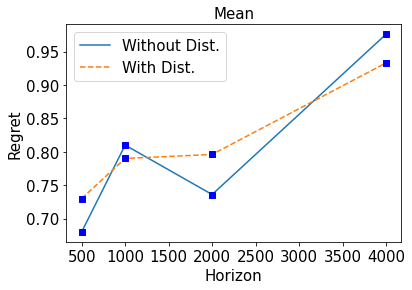}
\end{subfigure}%
\begin{subfigure}{.33\textwidth}
  \centering
  \includegraphics[width=.99\linewidth]{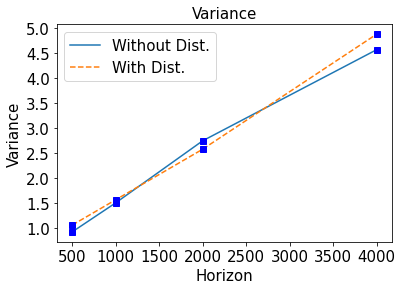}
\end{subfigure}%
\begin{subfigure}{.33\textwidth}
  \centering
  \includegraphics[width=.99\linewidth]{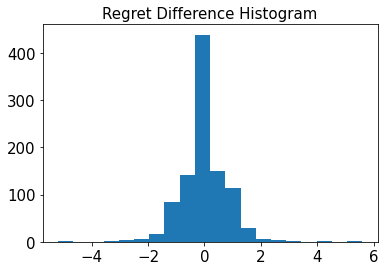}
\end{subfigure}
\caption{Algorithm performance with distribution knowledge v.s. without distribution knowledge.}
\label{fig2}
\end{figure}

Furthermore, we compare the performance of Algorithm \ref{alg:DAA} for the cases of with and without distributional knowledge under a nondegenerate problem instance. The result is presented in Figure \ref{fig2}. The first two plots report the mean and variance computed based on $200$ independent trials. The last plot displays the histogram of the difference of the regrets under the two cases. The histogram is generated under $800$ independent trials with horizon $T=1000$. For the case when we know the underlying distribution, Algorithm \ref{alg:DAA} is implemented based on replacing the estimate with the true distribution when solving the LP \eqref{adpt_lp}.

This result further substantiates the theme of our discussion: the learning of the distribution (under the finite support and nondegeneracy condition) will not affect the regret's dependency on the horizon $T$. The algorithm performances under the cases of known and unknown distribution match in both expectation and variance. In addition, the difference of their performances is also symmetrically distributed as seen from the histogram. 

Lastly, we make the following remarks:
\begin{itemize}
    \item Partial acceptance: In the algorithm and the analysis, we focus on the case of binary decision where the decision maker either accepts or rejects the order, $x_{t}=0$ or $1$. If partial acceptance is allowed, an alternative decision rule is to let $x_{t}=y_{j,t}^*$ in Algorithm \ref{alg:DAA}. Compared to the binary decision rule, this partial rule has one less layer of randomness but no other difference. In that case, all the analyses and regret bounds in this paper still hold.
    \item Known distribution: In this paper, we focus on the case when the underlying distribution is unknown. We remark that our analyses can also be applied to the case when the distribution is known. As noted earlier, the algorithm with known distribution replaces all the estimators in Algorithm \ref{alg:DAA} with their true parameters. For its analysis, all the probability events related to parameter estimation will happen with probability one and all the remaining analyses in our paper still go through. 
    \item Infrequent re-solving: Several existing works consider the case of known distribution and discuss the infrequent re-solving scheme where the adaptive LP \eqref{adpt_lp} in Algorithm \ref{alg:DAA} is solved on an infrequent basis (see, for example, \cite{reiman2008asymptotically,agrawal2014dynamic, jasin2015performance,bumpensanti2020re}). While our paper discusses the case of unknown distribution, we believe the infrequent re-solving scheme is still compatible with our analysis at least for the nondegenerate case. Specifically, in the proof of Proposition \ref{prop_multi_second_third_term}, we identify the number of observations needed to accurately estimate the optimal basis of the underlying LP. After these number of observations, the underlying distribution is estimated accurately enough and we can then adopt the algorithm proposed by \cite{bumpensanti2020re} using the estimated distribution. The re-solving time points can be chosen to meet the condition in $\mathcal{E}_t$, and then the analysis of the constraint process still holds. While the goal of our paper is mainly to provide a thorough analysis for the classic version of the algorithm with a frequent re-solving scheme, we believe such investigation is interesting and deserves future study.
    
    \item More efficient algorithm: Algorithm \ref{alg:DAA} can be computationally costly as it solves a linear program in every time period, so do the infrequent re-solving algorithms when the underlying LP is large-scale. A few recent works study the more efficient algorithms to solve the problem. \cite{li2020simple} and \cite{balseiro2020best} both perform gradient descent in the dual space and use the dual solution to guide the primal decisions. \cite{sun2020near} propose a gradient descent version of the infrequent re-solving algorithm in \cite{bumpensanti2020re}. \cite{gupta2021greedy} adapts the renowned Sum-of-Squares algorithm (originally for bin packing problem) and develops an efficient algorithm that achieves bounded regret with known distribution and nondegeneracy. Though originally motivated from a computational consideration, these works provide new perspectives for both algorithm design and analysis.
\end{itemize}

%
%
%

\bibliographystyle{informs2014} 
\bibliography{main} 

\appendix

\section{Proofs of Section \ref{Analysis_DLP}}\label{append_s3}

We first state three results that will be used in the proof.

\begin{lemma}\label{lem_multinomial}
    Denote $\Xi_d = \{\bm{q}\in \mathbb{R}^d| \sum_{j=1}^d q_j = 1, 0 < q_j < 1, j=1,...,d\}$. Let $\mathrm{Mul}(t,\bm{q})$ to be the multi-nomial distribution. For $\bm{q} \in \Xi_d$ and $\hat{\bm{q}} \sim \frac{1}{t} \cdot \mathrm{Mul}(t,\bm{q})$, the inequality
    \begin{equation}
    \begin{aligned}
    \prob(||\hat{\bm{q}} - \bm{q}||_1 \geq \epsilon) \leq 2^d\exp\left(-\frac{t\epsilon^2}{2}\right)
    \end{aligned}
\end{equation}
holds for any $\epsilon>0.$
\end{lemma}
\textit{Proof}:
See Appedix C1 in \cite{jaksch2010near}.

\begin{lemma}\label{lem_hoef}
(Hoeffding) Let $X_1,\cdots, X_T$ be independent random variables such that $a_t \leq X_t \leq b_t$ almost surely, we have 
$$\prob\left(\left|\sum_{t=1}^T X_t - \sum_{t=1}^T\mathbb{E}[X_t] \right| \geq \epsilon\right) \leq 2\exp\left(-\frac{2\epsilon^2}{\sum_{t=1}^T (b_t - a_t)^2}\right).$$
\end{lemma}

\begin{lemma}\label{lem_azume}
(Azuma-Hoeffding) Let $S_t$ be a martingale such that $S_0 = 0$ and the increment $S_t - S_{t-1}$ is bounded by $\sigma_t$ with probability $1$, then for any $\epsilon>0$ and $T\ge 2,$ we have
$$\prob(S_T \geq \epsilon ) \leq \exp\left(-\frac{\epsilon^2}{2\sum_{t=1}^T \sigma_t^2}\right).$$
\end{lemma}

\subsection{Proof of Proposition \ref{prop_general_regret_form}}
\label{ap_p1}
We have
\begin{equation*}
    \begin{aligned}
    \text{Reg}_{T}^\pi &= \E\left[T\cdot \text{OPT}_{\text{D}} - \sum_{t=1}^\tau r_tx_t\right]\\
    &= \E\left[\bm{\lambda}^{*\top}\bm{B} + \sum_{t=1}^T (r_t - \bm{a}_t^\top\bm{\lambda}^*)_+\right] - \sum_{t=1}^\tau\mathbb{E}\left[r_t x_t\right]\\
    &= \E\left[\bm{\lambda}^{*\top}\bm{B} + \sum_{t=1}^T (r_t - \bm{a}_t^\top\bm{\lambda}^*)_+\right] - \sum_{t=1}^\tau\mathbb{E}\left[(r_t - \bm{a}_t^\top\bm{\lambda}^*)x_t + \bm{a}_t^\top\bm{\lambda}^*x_t\right]\\
    &= \E\left[\bm{\lambda}^{*\top}\left(\bm{B}- \sum_{t=1}^\tau\bm{a}_tx_t\right)\right] + \mathbb{E}\left[ \sum_{t=1}^T(r_t - \bm{a}_t^\top\bm{\lambda}^*)_+ - \sum_{t=1}^\tau(r_t - \bm{a}_t^\top\bm{\lambda}^*)x_t \right]\\
    &= \bm{\lambda}^{*\top}\E\left[\bm{B}_\tau\right] + \mathbb{E}\left[ \sum_{t=1}^T(r_t - \bm{a}_t^\top\bm{\lambda}^*)_+ - \sum_{t=1}^\tau(r_t - \bm{a}_t^\top\bm{\lambda}^*)_+x_t \right] + \mathbb{E}\left[ \sum_{t=1}^\tau(\bm{a}_t^\top\bm{\lambda}^* - r_t)_+x_t \right]\\
    &= \bm{\lambda}^{*\top}\E\left[\bm{B}_\tau\right] + \E\left[\sum_{j\in \mathcal{J}^*} (\mu_j-\bm{c}_j^\top \bm{\lambda}^*)\cdot \left(n_j(T) -n_{j}^{a}(\tau)\right) \right]+ \E\left[\sum_{j\in \mathcal{J}'} (\bm{c}_j^\top \bm{\lambda}^*-\mu_j)\cdot n_{j}^{a}(\tau) \right],\\
    \end{aligned}
\end{equation*}
where the second line comes from duality of the DLP \eqref{eqn:DLP}, and the last line comes from the definition of $\bm{B}_\tau$ and the definition of optimal basis $\mathcal{J}^*$ and its compliment $\mathcal{J}'.$

\subsection{Proof of Corollary \ref{cor_general_regret_form}}\label{ap_c1}
Recall Proposition \ref{prop_general_regret_form} that
\begin{equation*}
    \begin{aligned}
    \mathrm{Reg}_{T}^\pi  = 
    \sum_{j\in \mathcal{J}^*} (\mu_j-\bm{c}_j^\top \bm{\lambda}^*)\cdot \E\left[n_j(T) -n_{j}^{a}(\tau) \right]+ \sum_{j\in \mathcal{J}'} (\bm{c}_j^\top \bm{\lambda}^*-\mu_j)\cdot \E\left[n_{j}^{a}(\tau) \right] + \bm{\lambda}^{*\top}\cdot \E\left[\bm{B}_{\tau}\right].
    \end{aligned}
\end{equation*}
Since $\tau' \leq \tau$, we know $n_{j}^{a}(\tau') \le n_{j}^{a}(\tau)$ for all $j$. If we apply this to the above equality, we obtain
\begin{equation*}
    \begin{aligned}
    \mathrm{Reg}_{T}^\pi &= 
    \sum_{j\in \mathcal{J}^*} (\mu_j-\bm{c}_j^\top \bm{\lambda}^*)\cdot \E\left[n_j(\tau') -n_{j}^{a}(\tau) \right]+ \sum_{j\in \mathcal{J}'} (\bm{c}_j^\top \bm{\lambda}^*-\mu_j)\cdot \E\left[n_{j}^{a}(\tau') \right] + \bm{\lambda}^{*\top}\cdot \E\left[\bm{B}_{\tau'}\right]\\
    &\hspace{5mm}+ 
    \sum_{j\in \mathcal{J}^*} (\mu_j-\bm{c}_j^\top \bm{\lambda}^*)\cdot \E\left[n_j(T) -n_{j}(\tau') \right]+ \sum_{j\in \mathcal{J}'} (\bm{c}_j^\top \bm{\lambda}^*-\mu_j)\cdot \E\left[n_{j}^{a}(\tau)-n_{j}^{a}(\tau') \right]\\
    &\hspace{5mm}+ \bm{\lambda}^{*\top}\cdot \E\left[\bm{B}_{\tau} - \bm{B}_{\tau'}\right]\\
    &\le 
    \sum_{j\in \mathcal{J}^*} (\mu_j-\bm{c}_j^\top \bm{\lambda}^*)\cdot \E\left[n_j(\tau') -n_{j}^{a}(\tau') \right]+ \sum_{j\in \mathcal{J}'} (\bm{c}_j^\top \bm{\lambda}^*-\mu_j)\cdot \E\left[n_{j}^{a}(\tau') \right] + \bm{\lambda}^{*\top}\cdot \E\left[\bm{B}_{\tau'}\right]\\
    &\hspace{5mm}+ 
    \sum_{j\in \mathcal{J}^*} (\mu_j-\bm{c}_j^\top \bm{\lambda}^*)\cdot \E\left[n_j(T) -n_{j}(\tau') \right]+ \sum_{j\in \mathcal{J}'} (\bm{c}_j^\top \bm{\lambda}^*-\mu_j)\cdot \E\left[n_{j}^{a}(\tau)-n_{j}^{a}(\tau') \right]\\
    &\hspace{5mm}+ \bm{\lambda}^{*\top}\cdot \E\left[\bm{B}_{\tau} - \bm{B}_{\tau'}\right].\\
    \end{aligned}
\end{equation*}

To bound the leftover terms compared to the statement of the corollary,
\begin{equation*}
    \begin{aligned}
    &\hspace{5mm}\sum_{j\in \mathcal{J}^*} (\mu_j-\bm{c}_j^\top \bm{\lambda}^*)\cdot \E\left[n_j(T) -n_{j}(\tau') \right]+ \sum_{j\in \mathcal{J}'} (\bm{c}_j^\top \bm{\lambda}^*-\mu_j)\cdot \E\left[n_{j}^{a}(\tau)-n_{j}^{a}(\tau') \right]+ \bm{\lambda}^{*\top}\cdot \E\left[\bm{B}_{\tau} - \bm{B}_{\tau'}\right]\\
    &\leq \max_{j\in[n]}|\mu_j-\bm{c}_j^\top \bm{\lambda}^*|\cdot\E\left[\sum_{j\in \mathcal{J}^*} \left(n_j(T) -n_{j}(\tau')\right) + \sum_{j\in \mathcal{J}'}  \left(n_{j}^{a}(\tau)-n_{j}^{a}(\tau')\right) \right]\\
    &\leq \max_{j\in[n]}|\mu_j-\bm{c}_j^\top \bm{\lambda}^*|\cdot\E\left[\sum_{j\in \mathcal{J}^*} \left(n_j(T) -n_{j}(\tau')\right) + \sum_{j\in \mathcal{J}'}  \left(n_{j}(T)-n_{j}(\tau')\right) \right]\\
    &= \max_{j\in[n]}|\mu_j-\bm{c}_j^\top \bm{\lambda}^*|\cdot\E\left[T-\tau' \right],
    \end{aligned}
\end{equation*}
where the first inequality comes from the fact that since  $\tau'\le \tau$, we have $\bm{B}_{\tau} \le \bm{B}_{\tau'}$ entry-wisely. The second inequality is due to $n_{j}^{a}(\tau)-n_{j}^{a}(\tau') \le n_{j}(\tau)-n_{j}(\tau') \le n_{j}(T)-n_{j}(\tau'),$ and the intuition is that for a certain order type, within a certain time window, the number of occurrence will always be greater than the number of acceptance. The last line is because of the definition of $n_j$ such that $\sum_{j=1}^nn_j(t)=t$ for $t > 0$.
Combining the above inequalities finishes the proof.

\subsection{Proof of Proposition \ref{prop_second_third_term}}\label{ap_p2}
It suffices to bound (i) the term $\mathbb{E}\left[n_{j}^{a}(\tau_S)\right]$ for $j \in \mathcal{J}'$ and (ii) the term $\mathbb{E}\left[n(\tau_S)  -n_{j}^{a}(\tau _S)\right]$ for $j \in \mathcal{J}^*_0 \coloneqq \{j: \mu_{j}>\bm{c}_{j}^\top \bm{\lambda}^*\}$ . 

Consider the event (the same as its definition in the main paper)
\begin{equation*}
    \begin{aligned}
     \mathcal{A}_{t}^{(j)} = \left\{\left|\frac{n_j(t-1)}{t-1} - p_j\right| \leq L\right\},
    \end{aligned}
\end{equation*}
and with the convention $\mathcal{A}_{1}^{(j)} = \Omega$. 

We first bound $\mathbb{E}\left[n_{j}^{a}(\tau_S)\right]$ where $j \in \mathcal{J}'$. Notice that
\begin{equation*}
    \begin{aligned}
     \mathbb{E}\left[n_{j}^{a}(\tau_S)\right] &\leq \mathbb{E}\left[\sum_{t=1}^TI(\{\text{Accept order type $j\in\mathcal{J}'$ at time $t$}\} \cap \{t < \tau_S\})\right]\\
     &\leq \sum_{t=1}^T\prob\left(\left\{\text{Accept order type $j\in\mathcal{J}'$ at time $t$}\right\} \cap \{t < \tau_S\}\right).\\
    \end{aligned}
\end{equation*}
By the definition of $\mathcal{A}_t^{(j)}$ and $L$, for $j \in \mathcal{J}'$, under the event $\{t < \tau_S\} \cap \left\{\cap_{j=1}^n\mathcal{A}_{t}^{(j)}\right\}$, we have that when solving the LP \eqref{adpt_lp}
$$\left|\frac{n_{j}(t-1)}{t-1}\bm{c}_j - p_j\bm{c}_j\right| \leq L\,\,\,\,\,\text{ and } \left|\frac{n_{j}(t-1)}{t-1}\mu_j - p_j\mu_j\right| \leq L.$$
The above inequalities combined with the fact $\{t < \tau_S\}$ ensures the condition in Lemma \ref{lem_gen_stability}.
Thus we know that the perturbed LP \eqref{adpt_lp} share the same optimal basis with DLP \eqref{eqn:DLP}. The probabilistic decision elicited by the algorithm based on \eqref{adpt_lp} will then reject the order type $j\in \mathcal{J}'.$ Therefore we have
\begin{equation*}
    \begin{aligned}
     \left\{\left\{\text{Accept order type $j\in\mathcal{J}'$ at time $t$}\right\} \cap \{t < \tau_S\}\right\} \subseteq \left\{\left\{\cup_{j=1}^n\bar{\mathcal{A}}_{t}^{(j)}\right\} \cap \{t < \tau_S\}\right\}.
    \end{aligned}
\end{equation*}
Thus, we have for $j\in\mathcal{J}'$
\begin{equation*}
    \begin{aligned}
     \mathbb{E}\left[n_{j}^{a}(\tau_S)\right] 
     &\leq \sum_{t=1}^T \prob\left(\left\{\text{Accept order type $j\in\mathcal{J}'$ at time $t$}\right\} \cap \{t < \tau_S\}\right)\\
     &\leq \sum_{t=1}^T\prob\left(\left\{\cup_{j=1}^n\bar{\mathcal{A}}_{t}^{(j)}\right\} \cap \{t < \tau_S\}\right)\\
     &\leq \sum_{t=1}^T\prob\left(\left\{\cup_{j=1}^n\bar{\mathcal{A}}_{t}^{(j)}\right\}\right) \leq \sum_{t=1}^T \sum_{j=1}^n\prob\left(\bar{\mathcal{A}}_{t}^{(j)}\right)\\
     &\leq n +  \sum_{t=1}^{T-1}2n\exp\left(-2tL^2\right)\leq \frac{2n}{1-\exp(-2L^2)}.
    \end{aligned}
\end{equation*}
Here the first two lines come from the above argument, the third line comes from the union bound, and the last line comes from concentration inequality (Lemma \ref{lem_hoef}).

Next, we start to bound $\mathbb{E}\left[n(\tau_S)  -n_{j}^{a}(\tau _S)\right]$ for $j \in \mathcal{J}^*_0 = \{j: \mu_{j}>\bm{c}_{j}^\top \bm{\lambda}^*\}$. In a similar way as before, under the event $\{t < \tau_S\} \cap \left\{\cap_{j=1}^n\mathcal{A}_{t-1}^{(j)}\right\}$, the LP \eqref{adpt_lp} shares the same optimal basis with DLP \eqref{eqn:DLP}. Therefore,
\begin{equation*}
    \begin{aligned}
     \left\{\left\{\text{Reject order type $j\in\mathcal{J}^*_0$ at time $t$}\right\} \cap \{t < \tau_S\}\right\} \subseteq \left\{\left\{\cup_{j=1}^n\bar{\mathcal{A}}_{t}^{(j)}\right\} \cap \{t < \tau_S\}\right\}.
    \end{aligned}
\end{equation*}
Then,  we have for $j\in \mathcal{J}^*_0$,
\begin{equation*}
    \begin{aligned}
     \mathbb{E}\left[n(\tau_S) - n_{j}^{a}(\tau _S)\right] 
     &\leq \mathbb{E}\left[\sum_{t=1}^T I(\{\text{Reject order type $j\in\mathcal{J}^*_0$ at time $t$} \cap \{t < \tau_S\}\})\right]\\
     &\leq \sum_{t=1}^T\prob\left(\left\{\text{Reject order type $j\in\mathcal{J}^*_0$ at time $t$}\right\} \cap \{t < \tau_S\}\right)\\
     &\leq \sum_{t=1}^T\prob\left(\left\{\cup_{j=1}^n\bar{\mathcal{A}}_{t}^{(j)}\right\} \cap \{t < \tau_S\}\right)\\
     &\leq \frac{2n}{1-\exp(-2L^2)},
    \end{aligned}
\end{equation*}
thereby finishing the proof.

\subsection{Proof of Lemma \ref{lem_bound_cor_last_two_terms}}\label{ap_nl3}
From definition we know that 
$$\prob(\tau_S \leq t) = \prob\left(\bm{b}_s\notin\mathfrak{B} \text{ for some }s\le t \right),$$
and since $\tau_S \geq 0$, we know
$$\E[\tau_S] = \sum_{t=1}^{+\infty}\prob(\tau_S \geq t) = \sum_{t=1}^{+\infty}\left(1 - \prob(\tau_S < t)\right) \geq \sum_{t=1}^{+\infty}\left(1 - \prob(\tau_S \leq t)\right)=T - \sum_{t=1}^T\prob\left(\bm{b}_s\notin\mathfrak{B} \text{ for some }s\le t \right).$$
Next, to finish the proof, notice that $\max_{j\in[n]}|\mu_j-\bm{c}_j^\top \bm{\lambda}^*| \leq ||\bm{\lambda}||_1$, and 
\begin{equation*}
    \begin{aligned}
    \bm{\lambda}^{*\top}\E\left[\bm{B}_{\tau_S}\right] \leq \bm{\lambda}^{*\top}\E\left[(T-\tau_S+2)\bm{b}_{\tau_S-1}\right] \leq \bm{\lambda}^{*\top}\left(\bm{b} + L\right)\E\left[(T-\tau_S+2)\right] \leq 2||\bm{\lambda}^*||_1\left(T+2-\E\left[\tau_S\right]\right),
    \end{aligned}
\end{equation*}
where the last inequality comes from the fact that for all the binding resources $i$ we will have initial average inventory $b_i \leq 1$, for all the nonbinding resource $i$ we know $\lambda^*_i = 0$, and lastly we can assume $L \leq 1$.

\subsection{Proof of Lemma \ref{lem_part1decompose}}\label{ap_l3}

We first analyze the process $\tilde{b}_{i,t}$, the $i$-th component of $\tilde{\bm{b}}_t$, and then take union bound with respect to $i$. Define $\mathcal{H}_{t}=\{(r_s, \bm{a}_s)\}_{s=1}^t$ for $t=1,...,T.$ Let $$Y_t\coloneqq\tilde{b}_{i,t+1}-\tilde{b}_{i,t}$$ for $t\ge1$ and $$X_t \coloneqq Y_t - \mathbb{E}[Y_t|\mathcal{H}_{t-1}].$$
In this way, to analyze the process $\tilde{b}_{i,t}$, we can equivalently analyze the summation $\sum_{s=1}^{t-1} Y_t$. From the definition of the process $\tilde{\bm{b}}_t$, we know that when $t\ge\tilde{\tau} - 1$, we have $$\tilde{b}_{i,t+1}=\tilde{b}_{i,t},$$
and when $1\le t<\tilde{\tau},$ we have
$$\tilde{b}_{i,t+1}= \tilde{b}_{i,t} +\frac{1}{T-t}(\tilde{b}_{i,t}-{a}_{i,t}x_{t}).$$
From the fact that $\tilde{b}_{i,t}$ is $\mathcal{H}_{t-1}$-measurable, we can bound the absolute value of $|x_t|$ such that
\begin{equation*}
    \begin{aligned}
        |X_t| &= \left|\frac{1}{T-t}(\tilde{b}_{i,t}-{a}_{i,t}x_{t}) - \mathbb{E}\left[\left. \frac{1}{T-t}(\tilde{b}_{i,t} - {a}_{i,t}x_{t})\right|\mathcal{H}_{t-1}\right]\right|\\
        &=\frac{1}{T-t}\left|\mathbb{E}\left[{a}_{i,t}x_t|\mathcal{H}_{t-1}\right] - {a}_{i,t}x_{t}\right| \leq \frac{1}{T-t}
    \end{aligned}
\end{equation*}
for each $t\le T-1.$
So we can define $L_t$ and $U_t$ as 
\begin{align*}
    L_t & \coloneqq -\frac{1}{T-t},\\
    U_t & \coloneqq \frac{1}{T-t},
\end{align*}
and the conditions of Theorem \ref{thm_dehling} are met for the process $X_t, L_t$ and $U_t$. Then as in Theorem \ref{thm_dehling}, 
$$V_t= \sum_{s=1}^t (U_s-L_s)^2 = \sum_{s=1}^t \frac{4}{(T-s)^2} \le \frac{4}{T-t-1}$$
for $t=1,...,T-2.$
From Theorem \ref{Hf_dependent}, we know that
\begin{equation}
  \prob\left(\left\vert\sum_{j=1}^s X_j \right\vert \ge \Delta \text{ for some } s\le t\right)\le 2e^{-\frac{\Delta^2(T-t-1)}{2}} \label{sum_X}
\end{equation}
holds for all $\Delta>0$ and $t\le T-2.$ 

With this bound on the summation of $X_t$, we return to analyze the summation of $Y_t$ by bounding the difference between these two sequences. By the definition, we have
  \begin{align}
    |X_t-Y_t| & = |\mathbb{E}[Y_t|\mathcal{H}_{t-1}]| \nonumber \\
    & =|\mathbb{E}[\tilde{b}_{i,t+1}-\tilde{b}_{i,t}|\mathcal{H}_{t-1}]|\nonumber \\
    & = \left\vert\frac{1}{T-t}\mathbb{E}[({a}_{i,t}x_{t}-\tilde{b}_{i,t})I(t<\tau)|\mathcal{H}_{t-1}] \right\vert \nonumber \\& \le \frac{\epsilon_t}{T-t} = \frac{1}{T-t}I(t\leq \kappa T) + \frac{1}{(T-t)t^{1/4}}I(t >\kappa T)
    \label{difference_X_Y}
  \end{align}
for $1\le t\le T-1.$ The second line comes from the definition of $Y_t$. The third line comes from splitting the second line with two indicators $I(t<\tau)$ and $I(t\ge \tau)$. The process $\tilde{\bm{b}}_t$ freezes and $\tilde{b}_{i,t+1}=\tilde{b}_{i,t}$ for $t\ge\tau$. The second last line comes from the definition of $\mathcal{E}_t$ and the definition of $\tau$, and the last line comes from the definition of $\epsilon_t$.
  
 By taking summation of (\ref{difference_X_Y}), we have for $s\leq T-2$,
  $$\left\vert\sum_{j=1}^s X_j - \sum_{j=1}^s Y_j \right\vert \le  \sum_{j=1}^{\kappa T} \frac{1}{T-j} + \sum_{j=\kappa T +1 }^{s} \frac{1}{(T-j)j^{1/4}}.$$
The next is to find a proper value of $\kappa$ such that the equation above is bounded by $\frac{\Delta}{2}$. For the first part, we have
 \begin{equation*}
    \begin{aligned}
    \sum_{j=1}^{\kappa T} \frac{1}{T-j} &\leq  \int_{T - \kappa T - 1}^{T-1} \frac{1}{x}dx = \log\left(\frac{T-1}{T -\kappa T - 1}\right)\\
    &\leq\log\left(\frac{T-1}{T - \kappa T- 2 - 2\kappa}\right) = \log\left(\frac{T-1}{T-2}\right) - \log\left(1-\kappa\right).
    \end{aligned}
\end{equation*}
For the second part,
\begin{equation*}
    \begin{aligned}
    \sum_{j=\kappa T+ 1}^{T-1} \frac{1}{(T-j)j^{1/4}} 
    &\leq \frac{1}{(\kappa T)^{1/4}} \sum_{j=\kappa T + 1}^{T-1} \frac{1}{T-j} \leq \frac{\log T}{(\kappa T)^{1/4}}.
    \end{aligned}
\end{equation*}
Henceforth, if we set $\kappa = 1-\exp(-\frac{\Delta}{8})$ and define $T_1$ as the minimal integer such that $T_1 \ge \frac{1}{\exp{\left(\frac{\Delta}{8}\right) - 1}}  + 2$ and 
$\frac{\log T_1}{T_1^{1/4}} \le \frac{\kappa^{1/4}\Delta}{4}$, then the following inequality holds for $T\ge T_1$
\begin{align*}
\sum_{j=1}^{\kappa T} \frac{1}{T-j} + \sum_{j=\kappa T +1 }^{T-1} \frac{1}{(T-j)j^{1/4}} \leq \frac{\Delta}{4} + \frac{\Delta}{4} = \frac{\Delta}{2}.
\end{align*}
With the choice of $\kappa$ and $T\ge T_1$, we have  $$\left\vert\sum_{j=1}^s X_j - \sum_{j=1}^s Y_j \right\vert \leq \frac{\Delta}{2}$$
holds almost surely. Consequently,
\begin{equation*}
    \begin{aligned}
    \left\{|\tilde{b}_{i,s}-b_i| > \Delta \text{ for some } s \leq t\right\} &= \left\{\left\vert\sum_{j=1}^{s-1} Y_j \right\vert > \Delta \text{ for some } s \leq t\right\}\\
    &= \left\{\left\vert\sum_{j=1}^{s} Y_j \right\vert > \Delta \text{ for some } s \leq t-1\right\}\\
    &\subseteq \left\{\left\vert\sum_{j=1}^{s} X_j \right\vert > \Delta/2 \text{ for some } s \leq t-1\right\}.\\
    \end{aligned}
\end{equation*}
Therefore, if we apply union bound with respect to constraint index $i=1,...,m$, we have
\begin{equation*}
    \begin{aligned}
    \prob\left(\tilde{\bm{b}}_s\notin \bigotimes_{i=1}^m \left[b_{i}-\Delta, b_{i}+\Delta\right] \text{ for some }s\le t\right)\le 2me^{-\frac{\Delta^2(T-t)}{8}}
    \end{aligned}
\end{equation*}
for $t\le T-2$ and $T\ge T_1$.

\subsection{Proof of Lemma \ref{lem_part2decompose}}\label{ap_l4}
In the statement before Lemma \ref{lem_part2decompose}, we have already shown $\left(\cap_{j=1}^n \mathcal{A}_{t}^{(j)}\right)\cap \left(\cap_{j=1}^n \mathcal{B}_{t}^{(j)}\right) \subseteq \mathcal{E}_t$ for $t \geq 2$, and we just have to discuss the special case that $\mathcal{E}_1 = \Omega$. From Algorithm \ref{alg:DAA}, we know that when $t=1$, we will always accept the order, and have $\bm{b} = \bm{b}_1$ with probability 1. Therefore, 
$$\E[\bm{a}_{1}x_{1}(\bm{b})|\mathcal{H}_{0},\bm{b}_{1}=\bm{b}]-\bm{b} = \E[\bm{a}_{1}]-\bm{b}.$$
Then notice that we have $-b_i\leq\E[a_{1,i}]-b_i\leq 1 - b_i$ for all $i$, because $||\bm{c}_j||_{\infty} \leq 1$ for all $j$. Moreover, all resource being binding implies that $0\leq b_i \leq 1$. Combining these facts we can have 
$$||\E[\bm{a}_{1}x_{1}(\bm{b})|\mathcal{H}_{0},\bm{b}_{1}=\bm{b}]-\bm{b}||_{\infty} = ||\E[\bm{a}_{1}]-\bm{b}||_{\infty} \leq 1.$$

Next, to show
$$\prob(\mathcal{E}_t) \ge \prob\left(\left(\cap_{j=1}^n \mathcal{A}_{t}^{(j)}\right)\cap \left(\cap_{j=1}^n \mathcal{B}_{t}^{(j)}\right)\right) \ge 1- \sum_{j=1}^n \prob\left(\bar{\mathcal{A}}_t^{(j)}\right) -  \sum_{j=1}^n\prob\left(\bar{\mathcal{B}}_t^{(j)}\right),$$
we analyze each component in the summation with Hoeffding's inequality (Lemma \ref{lem_hoef}). Specifically, $$n_j(t) = \sum_{s=1}^t \mathbb{I}((r_s, \bm{a}_s)=(\mu_j, \bm{c}_j))$$
where $\mathbb{I}(\cdot)$ denotes the indicator function and $\mathbb{I}((r_s, \bm{a}_s)=(\mu_j, \bm{c}_j))$'s are i.i.d. random variables. In addition, $\mathbb{E}[\mathbb{I}((r_s, \bm{a}_s)=(\mu_j, \bm{c}_j))] = p_j$ for each $j=1,...,n$. Therefore, we have
\begin{align*}
    \prob\left(\bar{\mathcal{A}}_t^{(j)}\right) &= \prob\left(\left|\frac{n_j(t-1)}{t-1} - p_j\right| > L\right) \\
    & \le 2\exp\left(-2L^2(t-1)\right)
\end{align*}
and 
\begin{align*}
    \prob\left(\bar{\mathcal{B}}_t^{(j)}\right) & =\prob\left(\left\vert \frac{n_j(t-1)}{t-1}-p_j \right\vert > \frac{1}{n(t-1)^{1/4}} \right)   \\
    & \le 2\exp\left(-\frac{2(t-1)^{1/2}}{n^2}\right).
\end{align*}
Combining these two inequalities, we have 
\begin{align*}
\prob\left(\left(\cap_{j=1}^n \mathcal{A}_{t}^{(j)}\right)\cap \left(\cap_{j=1}^n \mathcal{B}_{t}^{(j)}\right)\right) & \ge 1- 2n \exp\left(-2L^2(t-1)\right) -2n\exp\left(-\frac{2(t-1)^{1/2}}{n^2}\right).
\end{align*}

\subsection{Proof of Theorem \ref{theorem_regret}}\label{ap_t2}

We first state a lemma that takes a summation for both sides of the inequality in Lemma \ref{lem_part2decompose}. 

\begin{lemma}\label{lem_ineq_sum}
The following inequality holds for $T\in \mathbb{N}^+$
\begin{equation}\label{ineq_lemma_sum}
    \begin{aligned}
        \sum_{t=1}^T\sum_{s=\kappa T + 1}^t\mathbb{P}(\bar{\mathcal{E}}_s) \leq \cfrac{nT}{L^2} \exp\left(-2L^2 \kappa T\right) + 4n^3T^{3/2}\exp\left({-\frac{1}{n^2}(\kappa T)^{1/2}}\right).
    \end{aligned}
\end{equation}
\end{lemma}

\textit{Proof}:
Firstly, from Lemma \ref{lem_part2decompose},
    \begin{equation}\label{ineq_sum_part_1}
        \begin{aligned}
        \sum_{s=\kappa T + 1}^t \mathbb{P}(\bar{\mathcal{E}}_s) &\leq \sum_{s=\kappa T }^{t-1}  2n \exp\left(-2L^2s\right) +2n\exp\left(\frac{-s^{1/2}}{n^2}\right)\\
        & \leq \int_{\kappa T}^{T}\left(2n \exp\left(-2L^2s\right) +2n\exp\left(\frac{-s^{1/2}}{n^2}\right)\right)ds.\\
        \end{aligned}
    \end{equation}
For the second term observe that for any $\alpha > 0, \beta,\kappa \in (0,1)$ we can have the following bound
\begin{equation}\label{ineq_simplify}
        \begin{aligned}
        \int_{\kappa T}^T e^{-\alpha x^{\beta}}dx \leq \frac{-T^{1-\beta}}{\alpha\beta}\int_{\kappa T}^T \frac{-\alpha\beta}{x^{1-\beta}}e^{-\alpha x^{\beta}}dx\leq \frac{T^{1-\beta}}{\alpha\beta}e^{-\alpha(\kappa T)^{\beta}}.
        \end{aligned}
    \end{equation}
Combining $\eqref{ineq_sum_part_1}$ and $\eqref{ineq_simplify}$ yields
\begin{equation*}
    \begin{aligned}
        \sum_{t=1}^T\sum_{s=\kappa T + 1}^t\mathbb{P}(\bar{\mathcal{E}}_s) \leq \cfrac{nT}{L^2} \exp\left(-2L^2 \kappa T\right) + 4n^3T^{3/2}\exp\left({-\frac{1}{n^2}(\kappa T)^{1/2}}\right).
    \end{aligned}
\end{equation*}

\paragraph{\textbf{Proof of Theorem \ref{theorem_regret}}.} \ 

First we provide a slightly more careful analysis for \eqref{decompose}. Recall the stopping time
\begin{align*}
\tau_S \coloneqq \min \{t\le T:\bm{b}_t\notin \mathfrak{B}\} \cup \{T+1\}.
\end{align*}
Then, for $t=1,...,T,$
\begin{align*}
    \prob\left(\tau_S \leq t \right) & =\prob\left(\bm{b}_s\notin\mathfrak{B} \text{ for some }s\le t \right)  \\ & \le \prob\left(\tilde{\bm{b}}_s\notin\mathfrak{B} \text{ for some }s\le t \right) + \sum_{s=1}^t \prob((r_1,\Aa_1...,r_{s-1},\Aa_{s-1}) \notin \mathcal{E}_s) \nonumber \\
    &= \prob\left(\tilde{\bm{b}}_s\notin\mathfrak{B} \text{ for some }s\le t \right) + \sum_{s=\kappa T+1}^t \prob((r_1,\Aa_1...,r_{s-1},\Aa_{s-1}) \notin \mathcal{E}_s).
\end{align*} 
where the second line comes from \eqref{decompose} and the third line comes from the definition of $\mathcal{E}_s.$ By a ``more careful'' analysis, it means that the second component takes into account the definition of $\epsilon_t$ and thus removes the first $\kappa T$ summands. Next, we are going to apply Corollary \ref{cor_general_regret_form} to bound the regret.

With $T_1$ and $\kappa$ defined in Lemma \ref{lem_part1decompose} corresponding to $\Delta = L$, if $T \ge T_1$, we have
\begin{equation}\label{ineq_stptime_1}
    \begin{aligned}
    \mathbb{E}[\tau_S] &\geq \sum_{t=1}^T (1 - \mathbb{P}(\tau_S \leq t))\\
    &\geq T - 2 - \sum_{t=1}^{T-2} \prob\left(\tilde{\bm{b}}_s\notin\mathfrak{B} \text{ for some }s \leq t \right) - \sum_{t=1}^T \sum_{s=\kappa T + 1}^t \mathbb{P}(\bar{\mathcal{E}}_s)\\
    &\ge T - 2 - 2m \frac{1-e^{-(T-2)L^2/8}}{1-e^{-L^2/8}} - \sum_{t=1}^T\sum_{s=\kappa T+ 1}^t \mathbb{P}(\bar{\mathcal{E}}_s)
    \end{aligned}
\end{equation}
where the last line applies Lemma \ref{lem_part1decompose} for the first summation. The following result applies Lemma \ref{lem_part2decompose} for the second summation.

By combining $\eqref{ineq_stptime_1}$ with Lemma \ref{lem_ineq_sum}, and noting that
$$2m \frac{1-e^{-(T-2)L^2/8}}{1-e^{-L^2/8}} \le \frac{16m}{L^2},\,\text{ and }\,\,\kappa = 1-\exp\left(-\frac{L}{8}\right)\geq\frac{L}{32}\,\,\text{ when } \frac{L}{8} \leq 1,$$
we have
\begin{equation*}
    \begin{aligned}
    \E[\tau_S] &\geq T - 2 - \frac{16m}{L^2} - \cfrac{nT}{L^2} \exp\left(-\frac{L^3 T}{16}\right)- 4n^3T^{3/2}\exp\left({-\frac{L ^{1/2}}{6n^2}}\cdot T^{1/2}\right).
    \end{aligned}
\end{equation*}

Recall from Corollary \ref{cor_general_regret_form} we have 
\begin{equation*}
    \begin{aligned}
    \text{Reg}_{T}^\pi  & \leq \sum_{j\in \mathcal{J}^*} (\mu_j-\bm{c}_j^\top \bm{\lambda}^*)\cdot  \E\left[n_j(\tau')-n_{j}^{a}(\tau') \right]+ \sum_{j\in \mathcal{J}'} (\bm{c}_j^\top \bm{\lambda}^*-\mu_j)\cdot \E\left[n_{j}^{a}(\tau') \right] \\ & \ \ \ \  \ \ \ +\left(T -\E[\tau']\right)\cdot \max_{j\in[n]}|\mu_j-\bm{c}_j^\top \bm{\lambda}^*|+\bm{\lambda}^{*\top}\cdot \E\left[\bm{B}_{\tau'}\right].
    \end{aligned}
\end{equation*}
If we substitute $\tau_S$ for this $\tau' \leq \tau$ and apply \eqref{ineq_stptime_1}, Lemma \ref{lem_bound_cor_last_two_terms} and \ref{lem_ineq_sum}, we have
$$\left(T -\E[\tau_S]\right)\cdot \max_{j\in[n]}|\mu_j-\bm{c}_j^\top \bm{\lambda}^*| = O\left(\frac{m+n}{L^2} +  n^3T^{3/2}\exp\left(-\frac{T^{1/2}}{n^2}\right)\right),$$
\begin{equation*}
    \begin{aligned}
    \bm{\lambda}^{*\top}\E\left[\bm{B}_{\tau_S}\right] &\leq 2||\bm{\lambda}^*||_1\left(4+\frac{16m}{L^2}+ \cfrac{nT}{L^2} \exp\left(-\frac{L^3 T}{16}\right)+ 4n^3T^{3/2}\exp\left(-\frac{L ^{1/2}}{6n^2}\cdot T^{1/2}\right)\right)\\
    &= O\left(\frac{m+n}{L^2} +  n^3T^{3/2}\exp\left(-\frac{T^{1/2}}{n^2}\right)\right).\\
    \end{aligned}
\end{equation*}
Using the result for the bound of $\left(T -\E[\tau_S]\right)\cdot \max_{j\in[n]}|\mu_j-\bm{c}_j^\top \bm{\lambda}^*|$ and $\bm{\lambda}^{*\top}\E\left[\bm{B}_{\tau_S}\right]$, we have 
\begin{equation*}
    \begin{aligned}
    \text{Reg}_{T}^\pi &\leq \sum_{j\in \mathcal{J}^*} (\mu_j-\bm{c}_j^\top \bm{\lambda}^*)\cdot  \E\left[n_j(\tau_S)-n_{j}^{a}(\tau_S) \right]+ \sum_{j\in \mathcal{J}'} (\bm{c}_j^\top \bm{\lambda}^*-\mu_j)\cdot \E\left[n_{j}^{a}(\tau_S) \right]\\
      &\hspace{-8mm}+ \left(\max_{j\in[n]}|\mu_j-\bm{c}_j^\top \bm{\lambda}^*|+2||\bm{\lambda}^*||_1\right)\left(4+\frac{16m}{L^2}+ \cfrac{nT}{L^2} \exp\left(-\frac{L^3 T}{16}\right)+ 4n^3T^{3/2}\exp\left(-\frac{L ^{1/2}}{6n^2}\cdot T^{1/2}\right)\right)\\
    \end{aligned}
\end{equation*}
Lastly, we finish the proof by applying Proposition \ref{prop_second_third_term} to the first two terms. The binding assumption do not prevent us from getting a problem dependent parameter $L$ for the bound in Proposition \ref{prop_second_third_term}. Therefore we have
\begin{equation*}
    \begin{aligned}
    \text{Reg}_{T}^\pi  &\leq \frac{2n \max_j \left|\mu_j - \bm{c}_j^\top\bm{\lambda}^*\right|}{1-\exp(-2L^2)}\\
    &\hspace{-8mm}+ \left(\max_{j\in[n]}|\mu_j-\bm{c}_j^\top \bm{\lambda}^*|+2||\bm{\lambda}^*||_1\right)\left(4+\frac{16m}{L^2}+ \cfrac{nT}{L^2} \exp\left(-\frac{L^3 T}{16}\right)+ 4n^3T^{3/2}\exp\left(-\frac{L ^{1/2}}{6n^2}\cdot T^{1/2}\right)\right)\\
    \end{aligned}
\end{equation*}
by noticing that $1-\exp(-2L^2) \ge \frac{1}{2}L^2$ for $L\le 1$, we have
\begin{equation*}
    \begin{aligned}
    \text{Reg}_{T}^\pi
    &\leq \frac{(48m+4n+12)\cdot \|\bm{\lambda}^*\|_1}{L^2} + 3\|\bm{\lambda}^*\|_1\cdot \left(\cfrac{nT}{L^2} \exp\left(-\frac{L^3 T}{16}\right)+ 4n^3T^{3/2}\exp\left(-\frac{L ^{1/2}}{6n^2}\cdot T^{1/2}\right)\right)\\
    & = \frac{(48m+4n+12)\cdot \|\bm{\lambda}^*\|_1}{L^2} + o(1).
    \end{aligned}
\end{equation*}

\section{Proof of Section \ref{sec_degenerate}}
\subsection{Proof of Lemma \ref{lem_deg_events}}\label{ap_l5}
For notational simplicity, assume that there is only one constraint, i.e. $m=1$ and both $\bm{b}_t$ and $\bm{a}_t$ are one-dimensional. To obtain the multi-dimensional result, we can simply take a union bound.

Recall the dynamic of the constraint process 
$$b_{t+1} = b_t + \frac{b_t -a_tx_t}{T-t}.$$
Define
$$Y_t = b_{t+1} - b_t, \,\,\,X_t = Y_t - \mathbb{E}[Y_t|\mathcal{F}_{t-1}],$$
then we know $S_t = \sum_{s=1}^tX_s$ is a martingale, and the difference is bounded by
$$|S_t-S_{t-1}| = \left|\frac{b_t -a_tx_t - \mathbb{E}[b_t - a_tx_t |\mathcal{F}_{t-1}] }{T-t} \right| \leq \frac{1}{T-t}.$$
Next, denote the sample mean estimator as $\hat{\bm{p}}_s = \left(\frac{n_1(s)}{s},\cdots,\frac{n_n(s)}{s}\right)$ and define the event
$$\mathcal{G}_t = \cap_{s=1}^{t-1}\left\{||\hat{\bm{p}}_s - \bm{p}||_1 \leq \frac{\sqrt{4n\log 2T}}{\sqrt{s}}\right\}.$$
Under the event $\mathcal{G}_t$ we have
\begin{equation}\label{ineq_deg_boundY}
    \begin{aligned}
    \sum_{s=1}^{t-1} \mathbb{E}[Y_s|\mathcal{F}_{s-1}] &= \sum_{s=1}^{t-1}\frac{{b}_s - \mathbb{E}[{a}_sx_s|\mathcal{H}_{s-1}]}{T-s}\\
    &\geq \sum_{s=1}^{t-1}\frac{\sum_{j=1}^n {c}_jy_{j,s}\hat{p}_{j,s} - \sum_{j=1}^n {c}_jy_{j,s}p_{j}}{T-s}\\
    &\geq \sum_{s=1}^{t-1}\frac{-1}{(T-s)\sqrt{s}}\sqrt{4n\log 2T},
    \end{aligned}
\end{equation}
where the first inequality is because there are resources that might be non-binding, and the second inequality is because $||\bm{c}_j||_{\infty} \leq 1$ for any $j$ and $||\bm{y}_s||_{\infty} \leq 1$. Then, 
\begin{equation}
    \begin{aligned}
    &\hspace{5mm}\prob\left(b_t - b < -\frac{\sqrt{4n\log 2T}}{\sqrt{t}}-\frac{\sqrt{4n\log 2T} + \sqrt{2\log 2T}}{\sqrt{T-t}}\right)\\
    &= \prob\left(\sum_{s=1}^{t-1}X_s + \sum_{s=1}^{t-1} \mathbb{E}[Y_s|\mathcal{F}_{s-1}] < -\frac{\sqrt{4n\log 2T}}{\sqrt{t}}-\frac{\sqrt{4n\log 2T} + \sqrt{2\log 2T}}{\sqrt{T-t}}\right)\\
    &\le \prob\left(\left\{\sum_{s=1}^{t-1}X_s + \sum_{s=1}^{t-1} \mathbb{E}[Y_s|\mathcal{F}_{s-1}] < -\frac{\sqrt{4n\log 2T}}{\sqrt{t}}-\frac{\sqrt{4n\log 2T} + \sqrt{2\log 2T}}{\sqrt{T-t}}\right\}\cap\mathcal{G}_t\right)+\prob(\bar{\mathcal{G}}_t)\\
    &\leq \prob\left(\sum_{s=1}^{t-1}X_s < -\frac{\sqrt{2\log 2T}}{\sqrt{T-t}}\right) + \sum_{s=1}^{t-1}\prob\left(||\hat{\bm{p}}_s - \bm{p}||_1 > \frac{\sqrt{4n\log 2T}}{\sqrt{s}}\right)\\
    &\leq \frac{1}{2T} + \sum_{s=1}^{t-1}\frac{1}{4T^2} \leq \frac{1}{T}.
    \end{aligned}
\end{equation}
Here the first inequality comes from the introduction of the event $\mathcal{G}_t$'s,
the second inequality comes from \eqref{ineq_deg_boundY} and the inequality
$$\sum_{s=1}^{t}\frac{1}{(T-s)\sqrt{s}} \leq \frac{1}{\sqrt{t}}+\frac{1}{\sqrt{T-t}},$$
the third inequality comes from applying Lemma \ref{lem_azume} to get
$$\prob\left(\sum_{s=1}^{t-1}X_s < -\frac{\sqrt{1\log 2T}}{\sqrt{T-t}}\right) \leq \exp\left(-\frac{\frac{2\log 2T}{T-t}}{2\frac{1}{T-t+1}}\right) \leq \frac{1}{2T},$$ and the bound on $\prob(\bar{\mathcal{G}}_t)$ comes from applying Lemma \ref{lem_multinomial} and observe
\begin{equation*}
    \begin{aligned}
    \prob\left(||\hat{\bm{p}}_s - \bm{p}||_1 > \frac{\sqrt{4n\log 2T}}{\sqrt{s}}\right)&\leq 2^n\exp\left(-2n\log 2T\right)\leq\exp\left(n\log 2-2n\log 2T\right)\\
    &\leq\exp\left(-n\log 2T\right)\leq \exp\left(-2\log 2T\right) = \frac{1}{4T^2}.
    \end{aligned}
\end{equation*}
Notice that in the last line we assume the number of order types $n\ge 2$ because otherwise the problem will be trivial.

Having shown the property for $\mathcal{C}_t^{(i)}$, we then go ahead to show the result for $\mathcal{D}_t^{(j)}$. For $j\leq n$, from Hoeffding's inequality (Lemma \ref{lem_hoef}) we know that for $t > 1$
\begin{equation}
    \begin{aligned}
    \prob\left(\left|\frac{n_j(t-1)}{(t-1)p_j} -1\right| > \frac{\sqrt{\log 2T}}{\sqrt{2\underline{p}^2(t-1)}}\right) &= \prob\left(\left|\frac{n_j(t-1)}{(t-1)} - p_j\right| > \frac{\sqrt{\log 2T}p_j}{\sqrt{2\underline{p}^2(t-1)}}\right)\\
    &\leq 2\exp \left( -\frac{p^2}{\underline{p}^2}\log 2T \right) = \frac{1}{T},\\
    \end{aligned}
\end{equation}
thereby finishing the proof.

\subsection{Proof of Lemma \ref{lem_purturb_inter}}\label{ap_l6}
For simplicity, we define the perturbed vector
\begin{equation*}
    \begin{aligned}
     \bm{\xi}_t &:= \left(\frac{n_1(t-1)}{(t-1)p_1}, \cdots , \frac{n_n(t-1)}{(t-1)p_n}\right)^\top.\\
    \end{aligned}
\end{equation*}
Recall the DLP
\begin{equation}\label{eqn:deg_DLP}
    \begin{aligned}
   \text{OPT}_{\text{D}} \coloneqq \max \ \ & \bm{\mu}^\top \bm{y}   \\
    \text{s.t. }\ & \bm{C}\bm{y} \leq \bm{b} \\
    & \bm{0} \leq \bm{y} \leq \bm{1},
    \end{aligned}
\end{equation}
where $\bm{0}$ and $\bm{1}$ are $n$-dimensional vectors with entries being $0$ and $1$, respectively. Notice that at every time $t$, the sample LP \eqref{adpt_lp}
\begin{equation}\label{eqn:deg_SLP}
    \begin{aligned}
  R_t \coloneqq \max \ \ & \bm{\mu}_t^\top \bm{y}   \\
    \text{s.t. }\ & \bm{C}_t\bm{y} \leq \bm{b}_t \\
    & \bm{0} \leq \bm{y} \leq \bm{1},
    \end{aligned}
\end{equation}
is equivalent to 
\begin{equation}\label{eqn:deg_ESLP}
    \begin{aligned} 
   R_t \coloneqq \max \ \ & \bm{\mu}^\top \bm{y}'  \\
    \text{s.t. }\ & \bm{C}\bm{y}' \leq \bm{b}_t   \\
    & \bm{0} \leq \bm{y}' \leq \bm{\xi_t}
    \end{aligned}
\end{equation}
where $y_j' = y_j\frac{n_j(t-1)}{(t-1)p_j}$. The reason for the formulation \eqref{eqn:deg_ESLP} is that we can transform the randomness in $\bm{\mu}$, $\bm{C}$ to the randomness in the inventory process $\bm{b}_t$ and the constraint for $y_j$, and this will facilitate the analysis of the reward accumulated at each time period $t$. More specifically, we can view the objective value as a function of $\bm{b}_t$ and $\bm{\xi}_t$, and we can bound the difference of the objective function to the optimal value by bounding $\bm{b}_t - \bm{b}$ and  $\bm{\xi}_t - \bm{1}$. To analyze the dynamics of $\bm{b}_t$ and $\bm{\xi}_t$, we define events that will give us a ``right'' deviation to ensure a $\tilde{O}(\sqrt{T})$ regret.

Define $\text{OPT}(\bm{b}_t,\bm{\xi}_t)$ as the objective value of \eqref{eqn:deg_ESLP} with the right hand side constraint being $(\bm{b}_t,\bm{\xi}_t)$. Clearly, we have $$\text{OPT}(\bm{b}_t,\bm{\xi}_t) = R_t \text{ and } \text{OPT}(\bm{b},\bm{1}) = \text{OPT}_D.$$ The plan is to show the difference $\text{OPT}(\bm{b},\bm{1}) - \text{OPT}(\bm{b}+\Delta \bm{b}, \bm{1})$ and $\text{OPT}(\bm{b}+\Delta \bm{b}, \bm{1}) -\text{OPT}(\bm{b},\bm{1}+\Delta\bm{\xi})$ are bounded, and
then derive a bound for $\text{OPT}(\bm{b},\bm{1}) -\text{OPT}(\bm{b}+\Delta\bm{b},\bm{1}+\Delta\bm{\xi})$.

We start to show
$$\text{OPT}(\bm{b},\bm{1}) - \text{OPT}(\bm{b}+\Delta \bm{b}, \bm{1}) \leq  \bar{\lambda}||(- \Delta\bm{b})_+||_1.$$
We consider the dual program
\begin{equation}\label{eqn:standard_dual}
    \begin{aligned}
   \text{OPT}_{\text{Dual}}:=\min \ \ & \bm{b}^\top{\bm{\lambda}} + \sum_{j=1}^n \gamma_j \\
    s.t. \ & p_j\bm{c}_j^\top{\bm{\lambda}} + \gamma_j \ge p_j\mu_j, \ \
    j=1,...,n \\ 
    & {\lambda_i} \ge 0, i = 1,\cdots,m\\
    & \gamma_j \ge 0, j = 1,\cdots,n
    \end{aligned}
\end{equation}
which is equivalent to 
\begin{equation*}
    \begin{aligned}
   \text{OPT}_{\text{Dual}}:=\min \ \ & \bm{b}^\top{\bm{\lambda}} + \sum_{j=1}^n (\mu_j - \bm{c}_j^\top\bm{\lambda})_+ \\
    s.t.\,\,\,
    & {\lambda_i} \ge 0, i = 1,\cdots,m.
    \end{aligned}
\end{equation*}
The above program could be understood as a LP with finite dimension, where there exists finite many simplex solutions. Recall that $$\bar{\lambda} := \max \left\{||\bm{\lambda}||_{\infty} : \bm{\lambda} \in \mathcal{FD}_0 \right\}$$ and $\mathcal{FD}_0$ denotes the set of basic solutions for the dual of DLP \eqref{eqn:DLP}. We know for any $\Delta \bm{b}\in \mathbb{R}^m$, 
$$\text{OPT}(\bm{b},\bm{1}) - \text{OPT}(\bm{b}+\Delta \bm{b}, \bm{1}) \leq  \bar{\lambda}||(- \Delta\bm{b})_+||_1.$$
This is because there will always exists a basic solution $\bm{\lambda}$ that is the optimal solution, and the rate of change of the optimal value with respect to $\bm{b}$ will always be bounded by $\bar{\lambda}$.

Next, we begin to show
$$\text{OPT}(\bm{b} + \Delta \bm{b}, \bm{1}) - \text{OPT}(\bm{b} + \Delta \bm{b}, \bm{1}+\Delta\bm{\xi}) \leq ||\Delta\bm{\xi}||_1.$$
This is obvious by taking a look at \eqref{eqn:deg_ESLP}. If we increase/decrease the constraint of $y_j'$, the optimal value can at most increase proportionally to $\mu_j$. Combining the pieces together we have
\begin{equation*}
    \begin{aligned}
    &\hspace{5mm}\text{OPT}(\bm{b},\bm{1}) - \text{OPT}(\bm{b}+\Delta\bm{b},\bm{1}+\Delta\bm{\xi})\\ 
    &= \text{OPT}(\bm{b},\bm{1}) - \text{OPT}(\bm{b}+\Delta \bm{b},\bm{1}) + \text{OPT}(\bm{b}+\Delta \bm{b},\bm{1}) - \text{OPT}(\bm{b}+\Delta\bm{b},\bm{1}+\Delta\bm{\xi})\\
    &\leq \bar{\lambda}||(-\Delta\bm{b})_+||_1 + ||\Delta\bm{\xi}||_1
    \end{aligned}
\end{equation*}
Then, under the condition that $\bm{b}_t = \bm{b} + \Delta \bm{b}$, $\bm{\xi}_t = \bm{1}+\Delta \bm{\xi}$ and $\mathcal{C}_t\cap\mathcal{D}_t$, from equations above we have
\begin{equation*}
    \begin{aligned}
    \mathrm{OPT}_{\mathrm{D}} - R_t  &=\text{OPT}(\bm{b},\bm{1}) - \text{OPT}(\bm{b}_t,\bm{\xi}_t)\leq \bar{\lambda}||(\bm{b} - \bm{b}_t)_+||_1 + ||\bm{\xi}_t - \bm{1}||_1\\
    &\leq m\bar{\lambda}\left(\frac{\sqrt{4n\log 2T}}{\sqrt{t}}+\frac{\sqrt{4n\log 2T}+ \sqrt{2\log 2T}}{\sqrt{T-t}}\right)+n\frac{\sqrt{\log 2T}}{\sqrt{2\underline{p}^2(t-1)}}.
    \end{aligned}
\end{equation*}

\subsection{Proof of Theorem \ref{thm_degenerate}}\label{ap_t3}
With Lemma \ref{lemma_single_step}, we can have an upper bound for the single-period regret. By taking the summation, we can obtain an upper bound for the cumulative regret. More specifically, 
\begin{equation}
    \begin{aligned}
       \text{Reg}^{\pi}_T &= \sum_{t=1}^T\left(\bm{\mu}^\top\bm{y}^* - \mathbb{E}[\bm{a}_tx_t]\right)\\
       &= \sum_{t=1}^T\mathbb{E}\left[\text{OPT}_{\text{D}} - R_t\right]\\
       &\leq \sum_{t=1}^T \left(\mathbb{E}\left[\text{OPT}_{\text{D}} - R_t | \mathcal{C}_t\cap\mathcal{D}_t\right]\prob(\mathcal{C}_t\cap\mathcal{D}_t) + 1\cdot \prob(\bar{\mathcal{C}}_t\cup\bar{\mathcal{D}}_t)\right)\\
       &\leq 1 + \sum_{t=2}^T \left(\mathbb{E}\left[\text{OPT}_{\text{D}} - R_t| \mathcal{C}_t\cap\mathcal{D}_t\right]\prob(\mathcal{C}_t\cap\mathcal{D}_t) + \frac{n+m}{T}\right)\\
       &\leq 1 + \sum_{t=2}^T \left(\max\{1, \bar{\lambda}\}\left(m\left(\frac{\sqrt{4n\log 2T}}{\sqrt{t}}+\frac{\sqrt{4n\log 2T} + \sqrt{2\log 2T}}{\sqrt{T-t}}\right) + n\frac{\sqrt{\log 2T}}{\sqrt{2\underline{p}^2(t-1)}}\right) + \frac{n+m}{T}\right)\\
       &\leq \left(m\left(\sqrt{2}+\sqrt{16n}\right)+\frac{n}{\sqrt{2\underline{p}^2}}\right)\max\{1, \bar{\lambda}\}\sqrt{T}\sqrt{\log 2T} + 1 + n + m,
    \end{aligned}
\end{equation}
where in the first inequality we use the fact that $||\bm{\mu}||_{\infty} \leq 1$, therefore $\text{OPT}_{\text{D}} \leq 1$, and in the third inequality we apply Lemma \ref{lemma_single_step} to get the regret bound.

\section{Regret Analysis for Nondegenerate Case with Both Binding and Nonbinding Constraint}
\label{sec_reg_general}

In Section \ref{Analysis_DLP}, we analyze the nondegenerate case under the assumption that all the constraints are binding (Assumption \ref{assp_binding_new}). In this section we remove the assumption and discuss the more general case where for the underlying DLP $\eqref{eqn:DLP}$, both binding and non-binding constraints exist. We present the final result in Theorem \ref{thm_general}. We first note that Corollary \ref{cor_general_regret_form} and Proposition \ref{prop_second_third_term} hold without dependency on the bindingness of the underlying LP. So the remaining task is to reproduce the constraint process analysis under the general case to deal with the last two terms in Corollary \ref{cor_general_regret_form}.


\subsection{Constraint process under general case}

Recall the stopping time that the LP's structure changes
$$\tau_{S} \coloneqq \min \left\{t\le T: |b_{i,t}-b_i| > L  \text{ for some }i\in\mathcal{I}^* \right\} \cup \left\{t\le T: b_{i,t}-b_i < - L \text{ for some }i\in\mathcal{I}' \right\} \cup \{T+1\}.$$
As the arguments in Section \ref{secOrderAcc}, the parameter $L$ and stopping time $\tau_S$ are critical in bounding the order acceptance in Proposition \ref{prop_second_third_term}.

With a slight overload of the notation, we adjust the previous definition of the region $\mathfrak{B}$ as 
\begin{equation*}
    \begin{aligned}
   \mathfrak{B} &\coloneqq \left(\bigotimes_{i\in\mathcal{I}^*} [b_{i}-L, b_{i}+L]\right)\bigotimes \left(\bigotimes_{i\in\mathcal{I}'} \left[\left.b_{i}-L, +\infty\right)\right.\right).
    \end{aligned}
\end{equation*}
In this way, the stopping time $\tau_S$ can be expressed by
\begin{equation*}
    \begin{aligned}
    \tau_S = \min \left\{t \leq  T: \bm{b}_t \notin \mathfrak{B} \right\}\cup \{T+1\}.
    \end{aligned}
\end{equation*}
Unlike the previous section, the definition of $\mathfrak{B}$ here differentiates between binding and non-binding constraints. For binding constraints, the definition is the same as before, and it aims to capture the first time that the remaining average resource capacity $b_{i,t}$ deviates from $b_i$ by $L$. For non-binding constraints, we only concern a downward deviation -- a lower bound for the deviation because that is sufficient to guarantee that the non-binding constraints remain non-binding. 

Accordingly, we need to slightly adjust the previous definition of event $\mathcal{E}_t$ as follows. Denote $\underline{b} \coloneqq \min\{b_1,...,b_m\}$ and
let
$$\epsilon_t^* \coloneqq \begin{cases}
1 & t \leq \kappa T, \\
\frac{1}{t^{1/4}} & t >  \kappa T,
\end{cases}\,\,\,\,\,\epsilon_t' \coloneqq \begin{cases}
1 +L - \underline{b} & t \leq \kappa T, \\
\frac{1}{t^{1/4}} & t >  \kappa T,
\end{cases}$$
with $\kappa$ to be specified and 
$$\mathcal{E}_t \coloneqq \left\{\mathcal{H}_{t-1}\Big \vert \sup_{\bm{b}'\in \mathfrak{B} }\left\|\mathbb{E}[\bm{a}_{\mathcal{I}^*, t}x_{t}(\bm{b}')|\mathcal{H}_{t-1}]-\bm{b}'_{\mathcal{I}^*}\right\|_\infty \le \epsilon^*_{t-1} \text{ and } \sup_{\bm{b}'\in \mathfrak{B} }\mathbb{E}[\bm{a}_{\mathcal{I}', t}x_{t}(\bm{b}')|\mathcal{H}_{t-1}]-\bm{b}'_{\mathcal{I}'}\le \epsilon'_{t-1} \right\},$$
where the subscripts $\mathcal{I}^*$ and $\mathcal{I}'$ denote the corresponding dimensions of the vectors. Here we define different tolerance levels $\epsilon^*$ and $\epsilon'$ for binding and non-binding dimensions, respectively. The intuition is that for non-binding dimensions, we can tolerate larger deviation for the resource consumption as long as it does not sabotage the non-bindingness. 

According to the new event $\mathcal{E}_t$, we can define the stopping time 
$$\tilde{\tau} \coloneqq \min \{t\le T:\bm{b}_t\notin \mathfrak{B} \text{ or } \mathcal{H}_{t-1} \notin \mathcal{E}_t\} \cup \{T+1\}.$$
and the auxiliary process
$$\tilde{\bm{b}}_t = \begin{cases} 
\bm{b}_t, & t<\tilde{\tau},\\
\bm{b}_{\tilde{\tau}}, & t\ge \tilde{\tau}.
\end{cases}$$
These two definitions lead to the same decomposition as before.
\begin{equation}
    \begin{aligned}\label{eqn:general_term1_decomposition}
    \prob\left(\tau_S \leq t\right) = \prob\left(\bm{b}_s\notin \mathfrak{B} \text{ for some }s\le t \right) \le \prob\left(\tilde{\bm{b}}_s\notin\mathfrak{B} \text{ for some }s\le t \right) + \sum_{s=1}^t \prob((\bm{r}_1,\bm{a}_1...,\bm{r}_{s-1},\bm{a}_{s-1}) \notin \mathcal{E}_s).
    \end{aligned}
\end{equation}
The following lemma generalizes Lemma \ref{lem_part1decompose} to the general case (without Assumption \ref{assp_binding_new}). The proof follows the exact same arguments as Lemma \ref{lem_part1decompose}. As to the statement, for definition of constants there is an additional term compared to Lemma \ref{lem_part1decompose} and it arises from the non-binding dimensions. The proof is left in later subsections.
\begin{lemma}
\label{lem_general_term1_part1}
The following inequality holds for all $T\ge T_1$ and $t\le T-2$,
\begin{equation*}
    \begin{aligned}
    \prob\left(\tilde{\bm{b}}_s\notin\mathfrak{B} \text{ for some }s\le t\right)\le 2me^{-\frac{L^2(T-t)}{8}}
    \end{aligned}
\end{equation*}
where the constant $T_1$ is defined as the minimal integer such that 
$$T_1 \ge \left(\frac{1}{\exp{\left(\frac{L}{8}\right) - 1}}  + 2\right) \vee \left(\frac{1}{\exp{\left(\frac{L}{8(1+L-\underline{b})}\right) - 1}}  + 2\right)$$ and 
$\frac{\log T_1}{T_1^{1/4}} \le \frac{\kappa^{1/4} L}{4}$, where $\kappa$ is set by $\kappa = \left(1-\exp(-\frac{L}{8})\right) \wedge \left(1-\exp(-\frac{L}{8(1+L-\underline{b})})\right)$.
\end{lemma}

Next, for the second term in \eqref{eqn:general_term1_decomposition}, we keep the same definition of events $\mathcal{A}_t^{(j)}$ and $\mathcal{B}_t^{(j)}$ in the previous section such that
$$\mathcal{A}_t^{(j)} \coloneqq \left\{ \left\vert \frac{n_j(t-1)}{t-1}-p_j\right\vert \leq L \right\},$$
$$\mathcal{B}_t^{(j)} \coloneqq \left\{\left\vert \frac{n_j(t-1)}{t-1}-p_j\right\vert \leq \frac{1}{n(t-1)^{1/4}} \right\}.$$

Their definitions convey the same intuition as before. The event $\cap_{j=1}^n\mathcal{A}_t^{(j)}$ ensures that the sampled LP's bindingness structure aligns with that of the DLP, while the event $\cap_{j=1}^n\mathcal{B}_t^{(j)}$ ensures that given $\cap_{j=1}^n\mathcal{A}_t^{(j)}$, the expected resource consumption at time $t$ of the algorithm stays close to $\bm{b}_t$. The following lemma is analogous to Lemma \ref{lem_part2decompose} and provides an lower bound for the ``good'' event $\mathcal{E}_t.$ In addition, we can obtain a probability bound on the stopping time $\tau_S.$

\begin{lemma}
\label{lem_general_ABE}
We have $ \left(\cap_{j=1}^n\mathcal{A}_{t}^{(j)}\right) \cap \left(\cap_{j=1}^n\mathcal{B}_{t}^{(j)}\right) \subset \mathcal{E}_t$ for $t = 1,...,T$. Under Assumption \ref{assp_dist} and \ref{assp_nondeg}, we have for $t \leq T-2$,
$$\prob\left(\mathcal{E}_t\right)
\begin{cases}
 = 1 \,\,\,&\text{ for $t \leq \kappa T$} \\
\ge 1- 2n\exp\left(-2L^2(t-1)\right) -2n\exp\left(-\frac{2(t-1)^{1/2}}{n^2}\right) \,\,\,&\text{ for $t > \kappa T$}
\end{cases}
$$
Consequently,
\begin{equation*}
    \begin{aligned}
    \E[\tau_S] &\geq T - 2 - \frac{16m}{L^2} - \cfrac{nT}{L^2} \exp\left(-2L^2\kappa T\right)- 4n^3T^{3/2}\exp\left({-\frac{\kappa ^{1/2}}{n^2}\cdot T^{1/2}}\right).
    \end{aligned}
\end{equation*}

\end{lemma}

By putting the analysis of the stopping time $\tau_S$ together with the analysis in Corollary \ref{cor_general_regret_form} and Proposition \ref{prop_second_third_term}, we can obtain the final regret bound as in the theorem below. As the bound for the all-binding case in the previous section, the regret bound bears no dependency in terms of the time horizon $T$. Our analysis mainly focuses on removing the dependency on time $T$, and the result thus indicates that the adaptive design of the algorithm can significantly mitigate the effect of the parameter learning/estimation error on the regret.

\begin{theorem}\label{thm_general}
Under Assumption \ref{assp_dist} and \ref{assp_nondeg}, Algorithm \ref{alg:DAA} give a regret in the order of 
\begin{equation*}
    \begin{aligned}
    \mathrm{Reg}_{T}^\pi
     \le \frac{(48m+4n+12)\cdot \|\bm{\lambda}^*\|_1}{L^2} + o(1).
    \end{aligned}
\end{equation*}
\end{theorem}
In the following subsections, we elaborate the proofs for the results in above.

\subsection{Proof of Lemma \ref{lem_general_term1_part1}}

We treat the case for binding index and non-binding index separately. For binding index $i \in \mathcal{I}^*$ and $t=1,...,T$, let $$Y_t\coloneqq\tilde{b}_{i,t+1}-\tilde{b}_{i,t},\,\,\,\,\,\,X_t \coloneqq Y_t - \mathbb{E}[Y_t|\mathcal{H}_{t-1}].$$
For binding resources, the setup is completely the same as the proof in Lemma \ref{lem_part1decompose}. Therefore from the same definition that $\kappa \leq 1-\exp(-\frac{L}{8})$, and that $T_1$ is the minimal integer such that $T_1 \ge \frac{1}{\exp{\left(\frac{L}{8}\right) - 1}}  + 2$ and $\frac{\log T_1}{T_1^{1/4}} \le \frac{\kappa^{1/4} L}{4}$, we have 
\begin{equation*}
    \begin{aligned}
    \prob\left(\tilde{b}_{i,s} \notin \mathfrak{B} \text{ for some } s\le t\right)\le 2e^{-\frac{L^2(T-t)}{8 }}
    \end{aligned}
\end{equation*}
for $i \in \mathcal{I}^*$, $t\leq T-2$, and $T \geq T_1$.

Next, for non-binding index $i\in \mathcal{I}'$, we use the same definition on $Y_t, X_t$ such that
$$Y_t\coloneqq\tilde{b}_{i,t+1}-\tilde{b}_{i,t},\,\,\,\,\,\,X_t \coloneqq Y_t - \mathbb{E}[Y_t|\mathcal{H}_{t-1}],$$
and want to show
$$\sum_{j=1}^s X_j - \sum_{j=1}^s Y_j \leq \frac{L}{2}.$$
From the fact that 
\begin{equation*}
    \begin{aligned}
     b_{i,t}  - a_{i,t}x_t(\bm{b}_t) \geq \inf_{\bm{b}'\in \mathfrak{B}} \bm{b}_{\mathcal{I}'} - \mathbb{E}[\bm{a}_{\mathcal{I}', t}x_{t}(\bm{b}')|\mathcal{H}_{t-1}]\ge 
\begin{cases}
\underline{b} - L - 1 & t \leq \kappa T, \\
-\frac{1}{t^{1/4}} & t >  \kappa T,
\end{cases}
    \end{aligned}
\end{equation*}
We know 
\begin{equation*}
    \begin{aligned}
     X_t - Y_t &= -\mathbb{E}\left[Y_t|\mathcal{H}_{t-1}\right]\\
     &= -\mathbb{E}\left[\left.\frac{1}{T-t}(\bm{b}_t - \bm{a}_tx_t(\bm{b}_t)) I(\tilde{\tau} > t)\right|\mathcal{H}_{t-1}\right]\\
     &\leq \frac{\epsilon'_{t-1}}{T-t}\\
     &= \frac{1+L-\underline{b}}{T-t}I(t \leq \kappa T) + \frac{1}{t^{1/4}(T-t)}I(t > \kappa T)
    \end{aligned}
\end{equation*}
From the same approach in the proof of Lemma \ref{lem_part1decompose}, by defining $\kappa \leq 1-\exp(-\frac{L}{8(1+L-\underline{b})})$, and $T_1$ to be the minimal integer such that $T_1 \ge \frac{1}{\exp{\left(\frac{L}{8(1+L-\underline{b})}\right) - 1}}  + 2$ and $\frac{\log T_1}{T_1^{1/4}} \leq \frac{\kappa^{1/4} L}{4}$, we know that for $T > T_1$,
\begin{equation*}
    \begin{aligned}
    \sum_{j=1}^{\kappa T} \frac{1+L-\underline{b}}{T-j} &\leq  \left(1+L-\underline{b}\right)\left(\log\left(\frac{T-1}{T-2}\right) - \log\left(1-\kappa\right)\right)\leq \frac{L}{4},\\
    \sum_{j=\kappa T +1}^{ T-1}\frac{1}{(T-j)j^{1/4}}
    &\leq \frac{1}{(\kappa T)^{1/4}} \sum_{j=\kappa T + 1}^{T-1} \frac{1}{T-j} \leq \frac{\log T}{(\kappa T)^{1/4}} \leq \frac{L}{4}.
    \end{aligned}
\end{equation*}
Therefore, with the choice of $\kappa$ and $T\ge T_1$, we have
$$\sum_{j=1}^s X_j - \sum_{j=1}^s Y_j =  \sum_{j=1}^{\kappa T} \frac{1+L-\underline{b}}{T-j} + \sum_{j=\kappa T + 1}^{s} \frac{1}{(T-j)j^{1/4}} \leq \frac{L}{2}.$$
Next, similar to Lemma \ref{lem_part1decompose}, for $t\le T-2.$ we have that  
\begin{equation*}
  \prob\left(\sum_{j=1}^s X_j  \le -\frac{L}{2} \text{ for some } s\le t\right)\le e^{-\frac{L^2(T-t-1)}{8}}.
\end{equation*}
Then, we have for $i \in \mathcal{I}'$
\begin{equation*}
    \begin{aligned}
    \left\{\tilde{b}_{i,s}-b_i \leq -L \text{ for some } s \leq t\right\} &= \left\{\sum_{j=1}^{s-1} Y_j \leq -L \text{ for some } s \leq t\right\}\\
    &= \left\{\sum_{j=1}^{s} Y_j \leq -L \text{ for some } s \leq t-1\right\}\\
    &\subseteq \left\{\sum_{j=1}^{s} X_j \leq -\frac{L}{2} \text{ for some } s \leq t-1\right\}.\\
    \end{aligned}
\end{equation*}
Therefore, for $i \in \mathcal{I}'$ we have that
\begin{equation*}
    \begin{aligned}
    \prob\left(\tilde{b}_{i,s} \notin \mathfrak{B} \text{ for some } s\le t\right)\le 2e^{-\frac{L^2(T-t)}{8}}.
    \end{aligned}
\end{equation*}
Summing up and taking the union bound, we know that for all $T\ge T_1$ and $t\le T-2$, we have
\begin{equation*}
    \begin{aligned}
    \prob\left(\tilde{\bm{b}}_s\notin\mathfrak{B} \text{ for some }s\le t\right)\le 2me^{-\frac{L^2(T-t)}{8}}.
    \end{aligned}
\end{equation*}

\subsection{Proof of Lemma \ref{lem_general_ABE}}

To show that $\left(\cap_{j=1}^n\mathcal{A}_{t}^{(j)}\right) \cap \left(\cap_{j=1}^n\mathcal{B}_{t}^{(j)}\right)\subset \mathcal{E}_t$, we firstly denote that it suffices to show such property for $t \geq \kappa T$, because the requirement for $t \leq \kappa T$ holds with probability $1$ for both the binding resource in $\mathcal{I}^*$ and the non-binding resource in $\mathcal{I}'$. To see this, notice that for any $\bm{b}'\in \mathfrak{B}$ and $i\in\mathcal{I}^*$, we have 
$$- b'_i \leq \mathbb{E}[\bm{a}_{i, t}] - b'_i \leq 1-b'_i.$$
Because $\bm{b}'\in \mathfrak{B}$ and $i\in\mathcal{I}^*$, without loss of generality we can take $L \leq 1-b_i'$ then assume that $b'_i < 1$ (otherwise if there is no slackness for $b_i$, the problem will be degenerate). Therefore, we know that
$$\left\|\mathbb{E}[\bm{a}_{\mathcal{I}^*, t}x_{t}(\bm{b}')|\mathcal{H}_{t-1}]-\bm{b}'_{\mathcal{I}^*}\right\|_\infty \le 1.$$
For $\bm{b}'\in \mathfrak{B}$ and $i\in\mathcal{I}'$, from $b'_i > \underline{b}-L$ we know that 
$$\mathbb{E}[\bm{a}_{\mathcal{I}', t}x_{t}(\bm{b}')|\mathcal{H}_{t-1}]-\bm{b}'_{\mathcal{I}'}\le 1 + L - \underline{b}.$$
From above we can see that $P(\mathcal{E}_t) = 1$ for $t\leq \kappa T$.

For $t > \kappa T$, under $\left(\cap_{j=1}^n\mathcal{A}_{t}^{(j)}\right)$ and $\bm{b}_t = \bm{b}' \in \mathfrak{B}$, we know that the sampled LP \eqref{adpt_lp} is stable such that the optimal basis and the resource bindinness is the same as the DLP \eqref{eqn:DLP}. Then, Lemma \ref{lem_gen_stability} tell us that there exists $\hat{\bm{b}}$ such that $\hat{\bm{b}}_{\mathcal{I}^*} = \bm{b}'_{\mathcal{I}^*}$, and the same binding structure also implies that we can define $\hat{\bm{b}}_{\mathcal{I}'}$ in the non-binding dimension to be
$$\hat{\bm{b}} = \sum_{j=1}^n \bm{c}_{j} \frac{n_k(t-1)}{t-1}y_j^*(\bm{b}'),$$
where $\bm{y}^*$ is the solution for \eqref{adpt_lp}.
More specifically, from above we know that
\begin{equation}
    \begin{aligned}\label{eqn:general_feasibility}
        \hat{\bm{b}}_{\mathcal{I}^*} &=  \bm{b}'_{\mathcal{I}^*} = \sum_{j=1}^n \bm{c}_{\mathcal{I}^*, j} \frac{n_j(t-1)}{t-1}y_j^*(\bm{b}'),\\
        \hat{\bm{b}}_{\mathcal{I}'} &=  \sum_{j=1}^n \bm{c}_{\mathcal{I}', j} \frac{n_j(t-1)}{t-1}y_j^*(\bm{b}') \leq \bm{b}'_{\mathcal{I}'}.
    \end{aligned}
\end{equation}
Therefore, from equation $\eqref{eqn:general_feasibility}$ and Algorithm \ref{alg:DAA} we know that
\begin{equation*}
    \begin{aligned}
    \bm{b}'_{\mathcal{I}'} - \mathbb{E}[\bm{a}_{\mathcal{I}', t}x_{t}(\bm{b}')|\mathcal{H}_{t-1}] &= \bm{b}'_{\mathcal{I}'} -  \sum_{j=1}^n \bm{c}_{j} p_j y_j^*(\bm{b}')\\
    &= \bm{b}'_{\mathcal{I}'} - \hat{\bm{b}}_{\mathcal{I}'}  + \sum_{j=1}^n \bm{c}_{j} \left(\frac{n_j(t-1)}{t-1} - p_j\right) y_j^*(\bm{b}')\\
&\ge \left(\bm{b}_{\mathcal{I}'} - L\right) + \left(L - \bm{b}_{\mathcal{I}'}\right)  - \frac{1}{(t-1)^{1/4}} =  - \frac{1}{(t-1)^{1/4}},
    \end{aligned}
\end{equation*}
where in the last line, $\bm{b}'_{\mathcal{I}'} \geq \bm{b}_{\mathcal{I}'} - L$ is from the fact that $\bm{b}' \in \mathfrak{B}$, $-\hat{\bm{b}}_{\mathcal{I}'} \geq L - \bm{b}_{\mathcal{I}'}$ is from  Lemma \ref{lem_gen_stability} (more specifically, the positivity of slack variable and the stability of the optimal basic index), and the last one is from the inequality ensured by $\left(\cap_{j=1}^n\mathcal{B}_t^{(j)}\right)$ that $ \sum_{j=1}^n\left\vert \frac{n_j(t-1)}{t-1}-p_j\right\vert \leq \frac{1}{(t-1)^{1/4}} $. Therefore, under $\left(\cap_{j=1}^n\mathcal{A}_{t}^{(j)}\right) \cap \left(\cap_{j=1}^n\mathcal{B}_t^{(j)}\right)$ we have shown that
$$\sup_{\bm{b}'\in \mathfrak{B}}\mathbb{E}[\bm{a}_{\mathcal{I}', t}x_{t}(\bm{b}')|\mathcal{H}_{t-1}]-\bm{b}'_{\mathcal{I}'}\le \frac{1}{(t-1)^{1/4}}.$$
Then, in order to show $ \left(\cap_{j=1}^n\mathcal{A}_{t}^{(j)}\right) \cap \left(\cap_{j=1}^n\mathcal{B}_t^{(j)}\right) \subset \mathcal{E}_t$, it suffices to show that
$$\sup_{\bm{b}'\in \mathfrak{B}}\left\|\mathbb{E}[\bm{a}_{\mathcal{I}^*, t}x_{t}(\bm{b}')|\mathcal{H}_{t-1}]-\bm{b}'_{\mathcal{I}^*}\right\|_\infty \le \frac{1}{(t-1)^{1/4}}.$$
Under event $\left(\cap_{j=1}^n\mathcal{A}_{t}^{(j)}\right) \cap \left(\cap_{j=1}^n\mathcal{B}_{t}^{(j)}\right)$, from equation $\eqref{eqn:general_feasibility}$ and Algorithm \ref{alg:DAA} we know that
\begin{equation*}
    \begin{aligned}
    \left\|\mathbb{E}[\bm{a}_{\mathcal{I}^*, t}x_{t}(\bm{b}')|\mathcal{H}_{t-1}]-\bm{b}'_{\mathcal{I}^*}\right\|_{\infty} & = \left\|\sum_{j=1}^n \bm{c}_{\mathcal{I}^*, j} \left(p_j- \frac{n_j(t-1)}{t-1}\right)y_j^*(\bm{b}')\right\|_{\infty} \leq \frac{1}{(t-1)^{1/4}},
    \end{aligned}
\end{equation*}
where the last equation follows from the definition of $\cap_{j=1}^n\mathcal{B}_t^{(j)}$. Therefore, we are done with showing $\left(\cap_{j=1}^n\mathcal{A}_{t}^{(j)}\right) \cap \left(\cap_{j=1}^n\mathcal{B}_t^{(j)}\right) \subset \mathcal{E}_t$.

Next, for $t \leq T-2$, we have to show
$$\prob\left(\mathcal{E}_t\right)
\begin{cases}
 = 1 \,\,\,&\text{ for $t \leq \kappa T$} \\
\ge 1- 2n\exp\left(2L^2(t-1)\right) -2n\exp\left(-\frac{2(t-1)^{1/2}}{n^2}\right) \,\,\,&\text{ for $t > \kappa T$}
\end{cases}
$$
This is trivial since by definition of $\mathcal{E}_t$, when $t \leq \kappa T$, $\mathcal{E}_t$ is an event with probability 1. For $t > \kappa T$, the proof follows from applying Lemma \ref{lem_multinomial}.

Finally, the bound 
\begin{equation*}
    \begin{aligned}
    \E[\tau_S] &\geq T - 2 - \frac{16m}{L^2} - \cfrac{nT}{L^2} \exp\left(-2L^2\kappa T\right)- 4n^3T^{3/2}\exp\left({-\frac{\kappa ^{1/2}}{n^2}\cdot T^{1/2}}\right)
    \end{aligned}
\end{equation*}
follows from the same approach in \eqref{ineq_stptime_1}, Lemma \ref{lem_ineq_sum} and \ref{lem_general_term1_part1}, because all the terms in \eqref{ineq_stptime_1} have the same bound as their counterparts in Section \ref{sec_reg_binding}.

\subsection{Proof of Theorem \ref{thm_general}}
Recall the Corollary \ref{cor_general_regret_form} also holds without Assumption \ref{assp_binding_new}, and we have
\begin{equation*}
    \begin{aligned}
    \text{Reg}_{T}^\pi  &\leq \sum_{j\in \mathcal{J}^*} (\mu_j-\bm{c}_j^\top \bm{\lambda}^*)\cdot  \E\left[n_j(\tau')-n_{j}^{a}(\tau') \right]+ \sum_{j\in \mathcal{J}'} (\bm{c}_j^\top \bm{\lambda}^*-\mu_j)\cdot \E\left[n_{j}^{a}(\tau') \right] \\ & \ \ \ \ +\left(T -\E[\tau']\right)\cdot \max_{j\in[n]}|\mu_j-\bm{c}_j^\top \bm{\lambda}^*| + \bm{\lambda}^{*\top}\cdot \E\left[\bm{B}_{\tau'}\right].
    \end{aligned}
\end{equation*}
By substituting $\tau_S$ for $\tau'$, we can again bound the first two terms using Proposition \ref{prop_second_third_term}, and thereby having 
\begin{equation*}
    \begin{aligned}
    \text{Reg}_{T}^\pi  &\leq \frac{2n \max_j \left|\mu_j - \bm{c}_j^\top\bm{\lambda}^*\right|}{1-\exp(-2L^2)}  +\left(T -\E[\tau_S]\right)\cdot \max_{j\in[n]}|\mu_j-\bm{c}_j^\top \bm{\lambda}^*| + \bm{\lambda}^{*\top}\cdot \E\left[\bm{B}_{\tau_S}\right].
    \end{aligned}
\end{equation*}
As for the last two terms, from Lemma \ref{lem_general_ABE}
we know that
\begin{equation*}
    \begin{aligned}
    \E[\tau_S] &\geq T - 2 - \frac{16m}{L^2} - \cfrac{nT}{L^2} \exp\left(-2L^2\kappa T\right)- 4n^3T^{3/2}\exp\left({-\frac{\kappa ^{1/2}}{n^2}\cdot T^{1/2}}\right),
    \end{aligned}
\end{equation*}
and without loss of generality let us assume $1+L-\underline{b} \leq 1$ and $\frac{L}{8} \leq 1$, thereby having $\kappa = 1-\exp\left(-\frac{L}{8}\right)\geq\frac{L}{32}$, and 
\begin{equation*}
    \begin{aligned}
    \E[\tau_S] &\geq T - 2 - \frac{16m}{L^2} - \cfrac{nT}{L^2} \exp\left(-\frac{L^3 T}{16}\right)- 4n^3T^{3/2}\exp\left({-\frac{L^{1/2}}{6n^2}\cdot T^{1/2}}\right).
    \end{aligned}
\end{equation*}
Then, we have
\begin{equation*}
    \begin{aligned}
    &\hspace{5mm}\left(T -\E[\tau_S]\right)\cdot \max_{j\in[n]}|\mu_j-\bm{c}_j^\top \bm{\lambda}^*|\\ &\leq  \max_{j\in[n]}|\mu_j-\bm{c}_j^\top \bm{\lambda}^*|\left(4 + \frac{16m}{L^2} + \cfrac{nT}{L^2} \exp\left(-\frac{L^3 T}{16}\right)+ 4n^3T^{3/2}\exp\left({-\frac{L^{1/2}}{6n^2}\cdot T^{1/2}}\right)\right),
    \end{aligned}
\end{equation*}
and from Lemma \ref{lem_bound_cor_last_two_terms}, we have 
\begin{equation*}
    \begin{aligned}
    \bm{\lambda}^{*\top}\E\left[\bm{B}_{\tau_S}\right]
    &\leq 2||\bm{\lambda}^*||_1\left(4+\frac{16m}{L^2}+ \cfrac{nT}{L^2} \exp\left(-\frac{L^3 T}{16}\right)+ 4n^3T^{3/2}\exp\left({-\frac{L^{1/2}}{6n^2}\cdot T^{1/2}}\right)\right).\\
    \end{aligned}
\end{equation*}
Lastly noticing $1-\exp(-2L^2) \ge \frac{1}{2}L^2$ for $L\le 1$ and combining the results above we have
\begin{equation*}
    \begin{aligned}
    \text{Reg}_{T}^\pi  &\leq \frac{2n \max_j \left|\mu_j - \bm{c}_j^\top\bm{\lambda}^*\right|}{1-\exp(-2L^2)}  +\left(T -\E[\tau_S]\right)\cdot \max_{j\in[n]}|\mu_j-\bm{c}_j^\top \bm{\lambda}^*| + \bm{\lambda}^{*\top}\cdot \E\left[\bm{B}_{\tau_S}\right]\\
    & \leq \frac{4n \max_j \left|\mu_j - \bm{c}_j^\top\bm{\lambda}^*\right|}{L^2} + \left(\max_{j\in[n]}|\mu_j-\bm{c}_j^\top \bm{\lambda}^*|+2||\bm{\lambda}^*||_1\right)\cdot \\
    &\hspace{8mm}\left(4 + \frac{16m}{L^2} + \cfrac{nT}{L^2} \exp\left(-\frac{L^3 T}{16}\right)+ 4n^3T^{3/2}\exp\left({-\frac{L^{1/2}}{6n^2}\cdot T^{1/2}}\right)\right)\\
    &\leq \frac{(48m+4n+12)\cdot \|\bm{\lambda}^*\|_1}{L^2} + 3\|\bm{\lambda}^*\|_1\cdot \left(\cfrac{nT}{L^2} \exp\left(-\frac{L^3 T}{16}\right)+ 4n^3T^{3/2}\exp\left(-\frac{L ^{1/2}}{6n^2}\cdot T^{1/2}\right)\right)\\
    & = \frac{(48m+4n+12)\cdot \|\bm{\lambda}^*\|_1}{L^2} + o(1).
    \end{aligned}
\end{equation*}

\section{Multi-dimensional Case}
\label{sec_multi}

Now, we return to the general multi-dimensional formulation in LP $\eqref{eqn:multiILP}$. We note that the analysis of the multi-dimensional problem is essentially the same as the one-dimensional case discussed previously and the presentation here is mainly for completeness. So this section will mainly emphasize the different parts between the multi-dimensional case and the one-dimensional case, and the rest will follow. 

Firstly, for the multi-dimensional case, the distribution assumption becomes the following. The assumption and its interpretation here are parallel to Assumption \ref{assp_dist}.
\begin{assumption}\label{assp_multi} We assume
\begin{itemize}
    \item[(a)]  The column-coefficient pair $(\bm{r}_t,\bm{A}_t)$'s are i.i.d. sampled from a distribution $\mathcal{P}.$ The distribution $\mathcal{P}$ takes a finite and known support $\{(\bm{\mu}_j, \bm{c}_j)\}_{j=1}^n$ where $\bm{\mu}_j \in \mathbb{R}^k$ and $\bm{c}_j \in \mathbb{R}^{m\times k}$. Specifically, $\prob((\bm{r}_t, \bm{A}_t) = (\bm{\mu}_j, \bm{c}_j)) = p_j$ for $j =1,...,n$ and the parameters $\bm{p} = (p_1,...,p_n)^\top$ are unknown.
    \item[(b)] Positiveness and Boundedness: $\bm{\mu_j}, \bm{c}_j\ge \bm{0}$ and $\|\bm{\mu}_j\|_\infty, \|\bm{c}_j\|_\infty\le 1$ for $j=1,...,n.$
    \item[(c)] The right-hand-side $\bm{B}=T\bm{b}\ge 0$ where $\bm{b}=(b_1,...,b_m)^\top$.
\end{itemize}
\end{assumption}

The multi-dimensional DLP will become the following
\begin{equation}\label{eqn:multi_DLP}
    \begin{aligned}
   \max \ \ & \sum_{j=1}^n \bm{\mu}_j^\top\bm{y}_j p_j   \\
    \text{s.t. }\ & \sum_{j=1}^n \bm{c}_j \bm{y}_j p_j \le \bm{b} \\ 
    &\bm{1}^\top \bm{y}_j \le 1, \ \ \bm{y}_j \ge \bm{0}, \ \ j=1,...,n,\\
    \end{aligned}
\end{equation}
where $\bm{y}_j \in \mathbb{R}^k$.
Accordingly, the sampled LP that we solve in Algorithm \ref{alg:DAA} at time $t$ will be in the form
\begin{equation}\label{eqn:multi_SLP}
    \begin{aligned}
   \max \ \ & \sum_{j=1}^n \bm{\mu}_j^\top\bm{y}_j \frac{n_j(t-1)}{t-1}   \\
    \text{s.t. }\ & \sum_{j=1}^n \bm{c}_j \bm{y}_j \frac{n_j(t-1)}{t-1} \le \bm{b} \\ 
    &\bm{1}^\top \bm{y}_j \le 1, \ \ \bm{y}_j \ge \bm{0}. \ \ j=1,...,n.\\
    \end{aligned}
\end{equation}
Intuitively, we can understand the above LPs in the following way. At each round $t$, we are presented with an order bundle of type $j$ (from $n$ different order bundles), and each order bundle $j$ is drawn from an i.i.d distribution. Inside the order bundle type $j$, there are $k$ different orders, and we can accept up to one order out of the $k$ orders in the bundle or reject the bundle. We denote that the theoretical framework for the binding case and the general case (in Section \ref{sec_reg_binding} and Appendix \ref{sec_reg_general}) still applies to the proof of multi-dimensional case. Therefore, in the rest of the appendix we will just outline the procedure on how to get similar results. 

We note that there is a slight difference for the regret decomposition. Since there are multiple orders in each order bundle $j$, and we denote the set that contains all the order in order bundle $j$ to be $\mathcal{K}_j$. Moreover, we overload the notation such that $\bm{\lambda}^*$ is the optimal dual solution of \eqref{eqn:multi_DLP}, and $\mathcal{J}^*$ and $\mathcal{J}'$ is defined as
$$\mathcal{J}^*\coloneqq \{j: \exists \,l \in \mathcal{K}_j \text{ such that } \mu_l > \bm{c}_l^\top \bm{\lambda}^*,j=1,...,n\},$$
$$\mathcal{J}'\coloneqq \{j:\mu_l < \bm{c}_l^\top \bm{\lambda}^* \text{ for all } l \in \mathcal{K}_j,j=1,...,n\}.$$
The intuition is that in the multi-dimensional case, we will always accept one order in the bundle $j$ if there exists an order $l$ that features $\mu_l > \bm{c}_l^\top\bm{\lambda}^*$. On the opposite, we will reject the whole bundle if all the orders $l$ in the bundle feature $\mu_l < \bm{c}_l^\top\bm{\lambda}^*$.
Next, we denote $n_l^a(t)$ to be the accepted orders of type $l$ up to time $t$, and $n_j(t)$ to be the occurrence of order group $j$ up to time $t$.

\begin{proposition}\label{prop_multi_general_regret_form}
Under Assumption \ref{assp_multi} and assume the DLP \eqref{eqn:multi_DLP} is nondegenerate, the following equality holds 
\begin{equation*}
    \begin{aligned}
    \text{Reg}_{T}^\pi & = \bm{\lambda}^{*\top}\cdot \E\left[\bm{B}_{\tau}\right]
   + \sum_{j\in \mathcal{J}^*} \E\left[\max_{l\in \mathcal{K}_j}(\mu_l-\bm{c}_l^\top \bm{\lambda}^*)\cdot n_j(T) -\sum_{l\in \mathcal{K}_j}(\mu_l - \bm{c}_l^\top\bm{\lambda}^*)n_l^a(\tau) \right]\\
   &\hspace{20mm} + \sum_{j\in \mathcal{J}'}\sum_{l\in \mathcal{K}_j} (\bm{c}_l^\top \bm{\lambda}^*-\mu_l)\cdot \E\left[n_{l}^{a}(\tau) \right].
    \end{aligned}
\end{equation*}
\end{proposition}

To bound the second and third term, an important observation is that if we are able to always choose the right $l \in \mathcal{K}_j$ such that $\mu_l - \bm{c}_l^\top\bm{\lambda}^* = \max_{l\in \mathcal{K}_j}(\mu_l-\bm{c}_l^\top \bm{\lambda}^*)$, the second term will be $0$. Moreover, if we are able to reject all the orders in the order group $j \in \mathcal{J}'$, the third term will be $0$.

A similar corollary as Corollary \ref{cor_general_regret_form} transforms the above result to the case that involves a general stopping time $\tau'.$

\begin{corollary}
\label{coro_multi_general_regret_form}
The following inequality holds 
\begin{equation*}
    \begin{aligned}
    \text{Reg}_{T}^\pi  & \leq \bm{\lambda}^{*\top}\cdot \E\left[\bm{B}_{\tau'}\right]
   + \sum_{j\in \mathcal{J}^*} \E\left[\max_{l\in \mathcal{K}_j}(\mu_l-\bm{c}_l^\top \bm{\lambda}^*)\cdot n_j(\tau') -\sum_{l\in \mathcal{K}_j}(\mu_l - \bm{c}_l^\top\bm{\lambda}^*)n_l^a(\tau') \right]\\  & \hspace{10mm}
   + \sum_{j\in \mathcal{J}'}\sum_{l\in \mathcal{K}_j} (\bm{c}_l^\top \bm{\lambda}^*-\mu_l)\cdot \E\left[n_{l}^{a}(\tau') \right] +\left(T -\E[\tau']\right)\cdot \max_{j\in[n],l\in \mathcal{K}_j}|\mu_l-\bm{c}_l^\top \bm{\lambda}^*| 
    \end{aligned}
\end{equation*}
where $\tau'$ is a stopping time adapted to the process $\bm{B}_t$'s and $\tau'\le\tau$ almost surely. 
\end{corollary}

The following proposition mimics Proposition \ref{prop_second_third_term} to provide a bound for the second and third term for the above corollary. Its proof idea is similar to that of Proposition \ref{prop_second_third_term}, both of which utilize the stability result in Lemma \ref{lem_gen_stability}. 

\begin{proposition}
\label{prop_multi_second_third_term}
Under Assumption \ref{assp_multi} and assume the DLP \eqref{eqn:multi_DLP} is nondegenerate, the output of Algorithm \ref{alg:DAA} satisfies
\begin{equation*}
    \begin{aligned}
    \hspace{5mm}&\sum_{j\in \mathcal{J}^*} \E\left[\max_{l\in \mathcal{K}_j}(\mu_l-\bm{c}_l^\top \bm{\lambda}^*)\cdot n_j(\tau') -\sum_{l\in \mathcal{K}_j}(\mu_l - \bm{c}_l^\top\bm{\lambda}^*)n_l^a(\tau') \right]
   + \sum_{j\in \mathcal{J}'}\sum_{l\in \mathcal{K}_j} (\bm{c}_l^\top \bm{\lambda}^*-\mu_l)\cdot \E\left[n_{l}^{a}(\tau') \right]\\  
   &\leq   \frac{2n \max_j \left|\mu_j - \bm{c}_j^\top\bm{\lambda}^*\right|}{1-\exp(-2L^2)}  
    \end{aligned}
\end{equation*}
\end{proposition}

In the following subsection, we lay out the details of the proof.

\subsection{Analysis of the constraint process for the multi-dimensional case}
The dimension of the constraint process for the multi-dimensional case remains the same as before, and its dynamics can also be analyzed in a similar way. To complete the analysis, a few minor changes need to be made and are elaborated as follows.
We can carry out the same analysis as before. The only caveat is that we need to adjust the definitions for vectors $\bm{r}_t$'s and matrices $\bm{A}_t$'s as below. For simplicity, we assume the case that all the resource is binding again and one will see that the non-binding case follows by the same analysis. Define
$$\mathfrak{B}\coloneqq \bigotimes_{i=1}^m [b_{i}-L, b_{i}+L],$$
$$\mathcal{E}_t \coloneqq \left\{\mathcal{H}_{t-1}\Big \vert \sup_{\bm{b}'\in \mathfrak{B}}\left\|\mathbb{E}[\bm{A}_{t}\bm{x}_{t}(\bm{b}')|\mathcal{H}_{t-1}]-\bm{b}'\right\|_\infty \le \epsilon_{t-1} \right\},$$
$$
\epsilon_t \coloneqq \begin{cases}
1 & t \leq \kappa T, \\
\frac{1}{t^{1/4}} & t >  \kappa T,
\end{cases}
$$
Again by defining 
\begin{align*}
\tau_S \coloneqq \min \{t\le T:\bm{b}_t\notin \mathfrak{B}\} \cup \{T+1\},
\end{align*}
$$\tilde{\tau} \coloneqq \min \{t\le T:\bm{b}_t \notin \mathfrak{B} \text{ or } \mathcal{H}_{t-1} \notin \mathcal{E}_t\} \cup \{T+1\},$$
and
$$\tilde{\bm{b}}_t = \begin{cases} 
\bm{b}_t, & t<\tilde{\tau},\\
\bm{b}_{\tilde{\tau}}, & t\ge \tilde{\tau},
\end{cases}$$
we can have the same decomposition
\begin{align}
\prob\left(\tau_S\le t \right) &= \prob\left(\bm{b}_s\notin\mathfrak{B} \text{ for some }s\le t \right)\nonumber\\
&\le \prob\left(\tilde{\bm{b}}_s\notin\mathfrak{B} \text{ for some }s\le t \right) + \sum_{s=1}^t \prob((\bm{r}_1,\bm{A}_1...,\bm{r}_{s-1},\bm{A}_{s-1}) \notin \mathcal{E}_s)\label{ineq_multi_twoterms}
\end{align}
For the first term in \eqref{ineq_multi_twoterms}, we find that Lemma \ref{lem_part1decompose} follows easily since the analysis therein does not involve the extra dimension $k$. For the second term, we can just change the decision variables $y_j$'s from scalars to vectors and everything will follow.

From the definition of the events
\begin{equation*}
    \begin{aligned}
    \mathcal{A}_t^{(j)} \coloneqq \left\{ \left|\frac{n_j(t-1)}{t-1} - p_j\right| \leq L \right\} \text{ and } \mathcal{B}_t^{(j)} \coloneqq \left\{\left\vert \frac{n_j(t-1)}{t-1}-p_j\right\vert \leq \frac{1}{n(t-1)^{1/4}} \right\},
    \end{aligned}
\end{equation*}
we remark that the derivation of
$$ \left(\cap_{j=1}^n \mathcal{A}_{t}^{(j)}\right)\cap \left(\cap_{j=1}^n \mathcal{B}_{t}^{(j)}\right) \subset \mathcal{E}_t$$
is independent of extra dimension $k$ because only the norms on $\bm{p}$, $\bm{A}$ and $\bm{b}$ are required. Therefore, the result in Theorem \ref{theorem_regret} also holds for the multi-dimensional case under Assumption \ref{assp_multi} and the assumption of nondegeneracy.
We have laid out the same pathways as in the one-dimensional case, and from now on, the analysis in Section \ref{sec_reg_general} in terms of the handling of the non-binding constraints will be the same.

\subsection{Proof of Proposition \ref{prop_multi_general_regret_form}}
Denote $\mathcal{T}_t$ to be the arriving order group at time $t$, also notice that here $(r_t, \bm{a}_t)$ is the reward/consumption pair \textbf{accepted} by our algorithm at time $t$. 
\begin{equation*}
    \begin{aligned}
    \text{Reg}_{T}^\pi &= \E\left[T\cdot \text{OPT}_D - \sum_{t=1}^\tau r_tx_t\right]\\
    &= \E\left[\bm{\lambda}^{*\top}\bm{B} + \sum_{t=1}^T \max_{l\in \mathcal{T}_t}(\mu_l - \bm{c}_l^\top\bm{\lambda}^*)_+\right] - \sum_{t=1}^\tau\mathbb{E}\left[r_t x_t\right]\\
    &= \E\left[\bm{\lambda}^{*\top}\bm{B} + \sum_{t=1}^T \max_{l\in \mathcal{T}_t}(\mu_l - \bm{c}_l^\top\bm{\lambda}^*)_+\right] - \sum_{t=1}^\tau\mathbb{E}\left[(r_t - \bm{a}_t^\top\bm{\lambda}^*)x_t + \bm{a}_t^\top\bm{\lambda}^*x_t\right]\\
    &= \E\left[\bm{\lambda}^{*\top}\left(\bm{B}- \sum_{t=1}^\tau\bm{a}_tx_t\right)\right] + \mathbb{E}\left[ \sum_{t=1}^T\max_{l\in \mathcal{T}_t}(\mu_l - \bm{c}_l^\top\bm{\lambda}^*)_+ - \sum_{t=1}^\tau(r_t - \bm{a}_t^\top\bm{\lambda}^*)x_t \right]\\
    &= \bm{\lambda}^{*\top}\E\left[\bm{B}_\tau\right] + \mathbb{E}\left[ \sum_{t=1}^T\max_{l\in \mathcal{T}_t}(\mu_l - \bm{c}_l^\top\bm{\lambda}^*)_+ - \sum_{t=1}^\tau(r_t - \bm{a}_t^\top\bm{\lambda}^*)_+x_t \right] + \mathbb{E}\left[ \sum_{t=1}^\tau(\bm{a}_t^\top\bm{\lambda}^* - r_t)_+x_t \right]\\
    &= \bm{\lambda}^{*\top}\E\left[\bm{B}_\tau\right] +  \sum_{j\in \mathcal{J}^*} \E\left[\max_{l\in \mathcal{K}_j}(\mu_l-\bm{c}_l^\top \bm{\lambda}^*)\cdot n_j(T) -\sum_{l\in \mathcal{K}_j}(\mu_l - \bm{c}_l^\top\bm{\lambda}^*)n_l^a(\tau) \right]\\
    &\hspace{20mm}+ \sum_{j\in \mathcal{J}'}\sum_{l\in \mathcal{K}_j} (\bm{c}_l^\top \bm{\lambda}^*-\mu_j)\cdot \E\left[n_{l}^{a}(\tau) \right],\\
    \end{aligned}
\end{equation*}
where the second inequality comes from the fact that we are dealing with the DLP \eqref{eqn:multi_DLP} and the corresponding dual value is different, and the rest just follows.

\subsection{Proof of Proposition \ref{prop_multi_second_third_term}}
We note that the proof will be very similar to the proof of Proposition \ref{prop_second_third_term}, therefore we will only point out the differences.
Firstly, we define events such that
\begin{equation*}
    \begin{aligned}
     \mathcal{A}_{t}^{(j)} = \left\{\left|\frac{n_j(t-1)}{t-1} - p_j\right| \leq L\right\},
    \end{aligned}
\end{equation*}
and start to bound $\sum_{l\in \mathcal{K}_j}\E\left[n_{l}^{a}(\tau_S) \right]$ for $j \in \mathcal{J}'$. With the same decomposition we have
\begin{equation*}
    \begin{aligned}
     \sum_{l\in \mathcal{K}_j}\mathbb{E}\left[n_{l}^{a}(\tau_S)\right] &\leq \mathbb{E}\left[\sum_{t=1}^TI(\{\text{Accept any order in order group $j \in \mathcal{J}'$ at time $t$}\} \cap \{t < \tau_S\})\right]\\
     &\leq \sum_{t=1}^T\prob\left(\left\{\text{Accept any order in order group $j \in \mathcal{J}'$ at time $t$}\right\} \cap \{t < \tau_S\}\right).\\
    \end{aligned}
\end{equation*}
Notice that under the event $\{t < \tau_S\} \cap \left\{\cap_{j=1}^n\mathcal{A}_{t}^{(j)}\right\}$, the perturbed LP \eqref{eqn:multi_SLP} has the same optimal basis as \eqref{eqn:multi_DLP}. Hence for $j \in \mathcal{J}'$, from complementary slackness condition we know $x_l\cdot(\bm{c}_l^\top\bm{\lambda} - \mu_l) = 0$ for all $l \in \mathcal{K}_j$. Therefore we will reject all the order in order group $j$ with probability $1$, and
\begin{equation*}
    \begin{aligned}
     \left\{\left\{\text{Accept any order in order group $j \in \mathcal{J}'$ at time $t$}\right\} \cap \{t < \tau_S\}\right\} \subseteq \left\{\left\{\cup_{j=1}^n\bar{\mathcal{A}}_{t}^{(j)}\right\} \cap \{t < \tau_S\}\right\}.
    \end{aligned}
\end{equation*}
The rest of bounding $\sum_{l\in \mathcal{K}_j}\mathbb{E}\left[n_{l}^{a}(\tau_S)\right]$ follows the same approach in Proposition \ref{prop_second_third_term}.

Next, to bound the term $\E\left[\max_{l\in \mathcal{K}_j}(\mu_l-\bm{c}_l^\top \bm{\lambda}^*)\cdot n_j(\tau') -\sum_{l\in \mathcal{K}_j}(\mu_l - \bm{c}_l^\top\bm{\lambda}^*)n_l^a(\tau') \right]$ for $j \in \mathcal{J}^*$. We define the \textbf{suboptimal} order $l$ in order group $j$ to be the order such that $\mu_l - \bm{c}_l^\top\bm{\lambda}^* < \max_{l\in \mathcal{K}_j}(\mu_l - \bm{c}_l^\top\bm{\lambda}^*)$, and define the \textbf{optimal} order $l$ in order group $j$ to be such that $\mu_l - \bm{c}_l^\top\bm{\lambda}^* = \max_{l\in \mathcal{K}_j}(\mu_l - \bm{c}_l^\top\bm{\lambda}^*)$. Notice that 
\begin{equation*}
    \begin{aligned}
     &\hspace{5mm}\E\left[\max_{l\in \mathcal{K}_j}(\mu_l-\bm{c}_l^\top \bm{\lambda}^*)\cdot n_j(\tau') -\sum_{l\in \mathcal{K}_j}(\mu_l - \bm{c}_l^\top\bm{\lambda}^*)n_l^a(\tau') \right]\\
     &\leq 2\max_{l\in \mathcal{K}_j, j\in[n]}|\mu_l -\bm{c}_l^\top\bm{\lambda}^*|\\
     &\hspace{10mm}\cdot\mathbb{E}\left[\sum_{t=1}^TI(\{\text{Not accepting anything, or accepting the suboptimal $l \in \mathcal{K}_j$ and $j\in\mathcal{J}^*$}\} \cap \{t < \tau_S\})\right].\\
    \end{aligned}
\end{equation*}
Next, we analyze the case that happens with probability $1$ under the event $\{t < \tau_S\} \cap \left\{\cap_{j=1}^n\mathcal{A}_{t}^{(j)}\right\}$, where the optimal basis for the perturbed LP \eqref{eqn:multi_SLP} and DLP \eqref{eqn:multi_DLP} are the same. Observe that for $j \in \mathcal{J}^*$, no matter if there are one or multiple \textbf{optimal} order $l$ in group $j$, we will always have $\bm{1}^\top\bm{y}_j = 1$. This is again because of the complementary slackness condition : firstly, we have $(1 - \bm{1}^\top\bm{y}_j)(\max_{l\in \mathcal{K}_j}(\mu_l - \bm{c}_l^\top\bm{\lambda}^*)_+) = 0$, and this ensures $ \bm{1}^\top\bm{y}_j = 1$; secondly, we have 
$$y_l\cdot(\bm{c}_l^\top\bm{\lambda}^* - \mu_l + \max_{l\in \mathcal{K}_j}(\bm{c}_l^\top\bm{\lambda}^* - \mu_l)_+ ) = 0,$$
and this ensures that for $j \in \mathcal{J}^*$, $y_l > 0$ if and only if it is the \textbf{optimal} order. Henceforth, we have
\begin{equation*}
    \begin{aligned}
     \hspace{5mm}&\left\{\left\{\text{Not accepting anything, or accepting the suboptimal $l \in \mathcal{K}_j$ and $j\in\mathcal{J}^*$}\right\} \cap \{t < \tau_S\}\right\}\\
     &\subseteq \left\{\left\{\cup_{j=1}^n\bar{\mathcal{A}}_{t}^{(j)}\right\} \cap \{t < \tau_S\}\right\}.
    \end{aligned}
\end{equation*}
Lastly, the rest of the proof follows the proof of Proposition \ref{prop_second_third_term}.

\section{Proof of Lemma 1 and LP's Stability}
\label{secLPStab}

Here we explore a stability property for LPs by identifying conditions under which a small perturbation of LP's input will not change the optimal basis and the bindingness of the constraints. As a side product, we prove Lemma \ref{lem_gen_stability} and relate the constant $L$ therein with a number of parameters of the underlying LP. In the context of our paper, the perturbation of LP can be viewed as the estimation error of using \eqref{adpt_lp} as a proxy for \eqref{eqn:DLP}. We emphasize that the property is not pertaining to the resource allocation LP, so we present the result under the general standard form,
\begin{align}
    \min \ \ & \bm{c}^\top \bm{x},\nonumber \\
     \text{s.t.} \ \  & \bm{A} \bm{x} = \bm{b}, \label{LP1}\\
     & \bm{x} \ge \bm{0}.\nonumber
\end{align}
where $\bm{c}\in \mathbb{R}^n$, $\bm{A}\in \mathbb{R}^{m\times n},$ and $\bm{b}\in \mathbb{R}^n.$ We overload the notations a bit: within this subsection, the LP's input $\bm{c}$, $\bm{A}$ and $\bm{b}$ (and their dimensions $m,n$) all refer to a general vector or matrix, different from their contextual meanings in the previous sections. 
As the convention, we define the basic/non-basic variable set
$$\mathcal{B}^*=\{j:\bm{x}^*_j>0,j=1,...,n\},$$
$$\mathcal{B}'=\{j:\bm{x}^*_j=0,j=1,...,n\}$$
where $\mathcal{B}^*$ is also known as the optimal basis.  Consider a second LP of the same size as \eqref{LP1},
\begin{align}
        \min \ \ & \hat{\bm{c}}^{\top} \bm{x},\nonumber \\
     \text{s.t.} \ \  & \hat{\bm{A}} \bm{x} = \hat{\bm{b}}, \label{LP2}\\
     & \bm{x} \ge \bm{0},\nonumber
\end{align}
where $\hat{\bm{c}}\in \mathbb{R}^n$, $\hat{\bm{A}}\in \mathbb{R}^{m\times n},$ and $\hat{\bm{b}}\in \mathbb{R}^n.$ We can interpret $\hat{\bm{c}}, \hat{\bm{A}}, \hat{\bm{b}}$ as sample-based estimates or perturbations of their counterparts in \eqref{LP1}. The following proposition states the conditions under which the two LPs \eqref{LP1} and \eqref{LP2} share the same index set of basic variables. 

\begin{proposition}
\label{prop_stability_general}
Suppose the optimal solution of $\eqref{LP1}$ is unique and nondegenerate. Define
\begin{align*}
    \chi &\coloneqq \min \{{x}_j^*: {x}_j^* > 0\}, \\
    \sigma& \coloneqq \sigma_{\min}(\bm{A}_{\mathcal{B}^*}), \\
    \delta & \coloneqq \min \{\bm{A}_{j}^\top \bm{\lambda}^*-c_j: \bm{A}_{j}^\top\bm{\lambda}^*-c_j>0\}, 
\end{align*}
where $\bm{A}_{\mathcal{B}^*}$ is the sub-matrix of $\bm{A}$ containing the columns in $\mathcal{B}^*$, $\bm{A}_j$ is the $j$-th column of the matrix $\bm{A},$ and $\bm{\lambda}^*$ is the dual optimal solution of \eqref{LP1}.

If the following conditions hold
\begin{equation} \label{eqn:condition_stability}
    \begin{split}
        \left\|\hat{\bm{A}}_{j} - \bm{A}_{j}\right\|_{\infty}&
        \leq
        \left\{
        \begin{matrix}
            \frac{\min\{1,\sigma,\sigma^2\}\cdot \min\{\chi,\delta\}}{12m^2\sqrt{m}}, & \text{\ for $j\in\mathcal{B}^*$},\\
            \frac{\sigma\delta}{12m}, & \text{\ for $j\in \mathcal{B}'$},
        \end{matrix}
        \right.\\
        |\hat{c}_j - c_j| &
        \leq
        \left\{
        \begin{matrix}
            \frac{\sigma\delta}{12m}, & \text{\ for $j\in\mathcal{B}^*$},\\
            \frac{\delta}{6}, & \text{\ for $j\in \mathcal{B}'$},
        \end{matrix}
        \right.\\    
        |\hat{b}_i - b_i|&
        \leq \frac{\sigma\chi}{8\sqrt{m}}, \,\,\, \text{\ for $i = 1,...,m$},
    \end{split}
\end{equation}
then the LP \eqref{LP2} has the same index sets of basic and non-basic variables as the LP \eqref{LP1}. In addition, the optimal solution of \eqref{LP1} is also unique and nondegenerate. 
\end{proposition}

The condition in the proposition is expressed by the deviation of the perturbed LP's inputs from the original LP's input. It relates the stability of $\mathcal{B}^*$ with a few quantities of the LP: $\chi$ captures the stability of the original primal optimal solution; $\sigma$ describes the singularity of the constraint matrix restricted to the columns in $\mathcal{B}^*$; $\delta$ is a sub-optimality measure for the non-basic variables and it is computed based on the reduced costs. The proof of the proposition is based on standard linear algebra analysis and it formalizes the intuition that the optimal basis $\mathcal{B}^*$ should exhibit some continuity with respect to the LP's input. 

The nondegeneracy assumption (Assumption \ref{assp_nondeg}) can be further illustrated from Proposition \ref{prop_stability_general}. The role of the assumption is to ensure $\chi, \sigma$ and $\delta$ defined in the proposition to be positive, and consequently, the stability property holds for the DLP \eqref{eqn:DLP}. To achieve a bounded regret for the online resource allocation problem, the nondegeneracy assumption is indeed necessary so that the parameters $\chi, \sigma$ and $\delta$ can be treated as constant and not dependent on $T$. On one hand, when the nondegeneracy assumption is violated or when the parameters such as $\chi, \sigma$ and $\delta$ may scale with  $\frac{1}{\sqrt{T}}$ or $\frac{1}{T}$, there can be examples for which constant regret is not achievable \citep{arlotto2019uniformly, bumpensanti2020re}. On the other hand, the constant condition is arguably reasonable in the nondegeneracy context in that these parameters are computed based on the DLP \eqref{eqn:DLP} and its standard form both of which bear no dependence of $T$.

\subsection{Proof of Proposition \ref{prop_stability_general}}

Recall that $\mathcal{B}^*$ and $\mathcal{B}'$ denote the optimal and non-optimal basis for \eqref{LP1}, respectively. The idea here is to show that the perturbed LP \eqref{LP2} has the same optimal and non-optimal basis under condition \eqref{eqn:condition_stability}.
Consider a basic solution $\hat{\bm{x}}$ of the perturbed LP \eqref{LP2} defined as follows. If $\hat{\bm{A}}_{\mathcal{B}^*}$ is invertible, we can define $\hat{\bm{x}}$ as
\begin{align*}
(\hat{\bm{x}})_{\mathcal{B}^*}  & = (\hat{\bm{A}}_{\mathcal{B}^*})^{-1}\hat{\bm{b}},\\
(\hat{\bm{x}})_{\mathcal{B}'}  & = \bm{0}
\end{align*}

Next, we prove the following results:
\begin{itemize}
    \item[(a)] The matrix $\hat{\bm{A}}_{\mathcal{B}^*}$ is non-singular and thus $(\hat{\bm{x}})_{\mathcal{B}^*}$ is a well-defined basic solution.
    \item[(b)] $(\hat{\bm{x}})_{\mathcal{B}^*} > \bm{0}$, and thus $\hat{\bm{x}}$ is a basic feasible solution of $\eqref{LP2}$.
    \item[(c)] The reduced costs associated with non-basic variables in $\mathcal{B}'$ are all negative and hence $\hat{\bm{x}}$ is the unique optimal solution of the perturbed LP $\eqref{LP2}$.
\end{itemize}

Throughout the proof, we use $\bm{A}^\epsilon$ to denote $\hat{\bm{A}} - \bm{A}$, and similarly for $\bm{c}^\epsilon$ and $\bm{b}^\epsilon$. To show part (a), we prove that the smallest singular value of the matrix is positive. We use $\sigma_{\min}(\bm{M})$ and $\sigma_{\max}(\bm{M})$ to denote the smallest and the largest singular value of a matrix $\bm{M}$. Then we have
\begin{align*}
   \sigma_{\min}\left(\hat{\bm{A}}_{\mathcal{B}^*}\right)
    &\geq
    \sigma_{\min}\left(\bm{A}_{\mathcal{B}^*}\right) -\sigma_{\max}\left(-\bm{A}^{\epsilon}_{\mathcal{B}^*}\right)\\
    &=
    \sigma - \sigma_{\max}\left(-\bm{A}^{\epsilon}_{\mathcal{B}^*}\right)\\
    &\geq
    \sigma - \sqrt{m}\|\bm{A}^{\epsilon}_{\mathcal{B}^*}\|_{\infty},
\end{align*}
where the first line comes from Weyl's inequality on matrix eigenvalues/singular values, the second line comes from the definition of $\sigma$, and the third line is obtained from the relation between the spectral norm and the infinity norm of a matrix. From condition $\eqref{eqn:condition_stability}$ and $||\bm{A}^{\epsilon}_{\mathcal{B}^*}||_{\infty} \leq m\max\limits_{i\in\mathcal{B}^*}\left\|\bm{A}^{\epsilon}_{i}\right\|_{\infty}$, we have
$$\sigma_{\min}\left(\hat{\bm{A}}_{\mathcal{B}^*}\right) \geq \frac{\sigma}{2},$$ 
and consequently,
\begin{align}
\label{bdd_purturb_matrix_sigma}
\sigma_{\max}\left(\left(\hat{\bm{A}}_{\mathcal{B}^*}\right)^{-1}\right)\leq\frac{2}{\sigma}.
\end{align}

For part (b), we show $\hat{\bm{x}}$ is a feasible basic solution. From Assumption \ref{assp_nondeg} we know that
$$(\bm{x}^*)_{\mathcal{B}^*} \geq \chi >0,$$
where the inequality holds element-wise. To ensure that $\left(\hat{\bm{x}}\right)_{\mathcal{B}^*}$ is strict positive, it suffices to show
$$
\left\|(\hat{\bm{x}})_{\mathcal{B}^*}-(\bm{x}^*)_{\mathcal{B}^*}\right\|_{\infty}
=
\left\|(\bm{A} + \bm{A}^{\epsilon})_{\mathcal{B}^*}^{-1}(\bm{b} + \bm{b}^{\epsilon})-(\bm{A}_{\mathcal{B}^*})^{-1}\bm{b}\right\|_{\infty}
\leq \frac{\chi}{2}.
$$
From condition $\eqref{eqn:condition_stability}$, if $\max\limits_{i\in{\mathcal{B}^*}}\|\bm{A}_i^{\epsilon}\|_{\infty}\leq \frac{\sigma^2\chi}{8m\sqrt{m}}$ and $\|\bm{b}^{\epsilon}\|_{\infty}\leq \frac{\sigma\chi}{8\sqrt{m}}$ , then

\begin{align*}
    &\hspace{4mm}\left\|(\bm{A} + \bm{A}^{\epsilon})_{\mathcal{B}^*}^{-1}(\bm{b} + \bm{b}^{\epsilon})-(\bm{A}_{\mathcal{B}^*})^{-1}\bm{b}\right\|_{\infty}\\
    &\leq
     \left\|(\bm{A} + \bm{A}^{\epsilon})_{\mathcal{B}^*}^{-1}\bm{b} -(\bm{A}_{\mathcal{B}^*})^{-1}\bm{b}\right\|_{\infty} +  \left\|(\bm{A} + \bm{A}^{\epsilon})_{\mathcal{B}^*}^{-1}\bm{b}^\epsilon\right\|_{\infty}\\
    &\leq
    \left\|(\bm{A} + \bm{A}^{\epsilon})_{\mathcal{B}^*}^{-1} -(\bm{A}_{\mathcal{B}^*})^{-1}\right\|_{\infty} \left\|\bm{b}\right\|_{\infty} +  \left\|(\bm{A} + \bm{A}^{\epsilon})_{\mathcal{B}^*}^{-1}\bm{b}^\epsilon\right\|_{\infty}\\
    &\leq
    \left\|(\bm{A}_{\mathcal{B}^*})^{-1}\left(\bm{A}_{\mathcal{B}^*}(\bm{A} + \bm{A}^{\epsilon})_{\mathcal{B}^*}^{-1} -\bm{I}\right)\right\|_{\infty}  +  \left\|(\bm{A} + \bm{A}^{\epsilon})_{\mathcal{B}^*}^{-1}\bm{b}^\epsilon\right\|_{\infty}\\
    &=\left\|(\bm{A}_{\mathcal{B}^*})^{-1}\left((\bm{A} +\bm{A}^\epsilon - \bm{A}^\epsilon )_{\mathcal{B}^*}(\bm{A} + \bm{A}^{\epsilon})_{\mathcal{B}^*}^{-1} -\bm{I}\right)\right\|_{\infty}  +  \left\|(\bm{A} + \bm{A}^{\epsilon})_{\mathcal{B}^*}^{-1}\bm{b}^\epsilon\right\|_{\infty}\\
    &=\left\|(\bm{A}_{\mathcal{B}^*})^{-1}\bm{A}^\epsilon _{\mathcal{B}^*}(\bm{A} + \bm{A}^{\epsilon})_{\mathcal{B}^*}^{-1}\right\|_{\infty}  +  \left\|(\bm{A} + \bm{A}^{\epsilon})_{\mathcal{B}^*}^{-1}\bm{b}^\epsilon\right\|_{\infty}\\
    &\leq \left\|(\bm{A}_{\mathcal{B}^*})^{-1}(\bm{A} + \bm{A}^{\epsilon})_{\mathcal{B}^*}^{-1}\right\|_{\infty}\|\bm{A}^\epsilon _{\mathcal{B}^*}\|_{\infty}  +  \left\|(\bm{A} + \bm{A}^{\epsilon})_{\mathcal{B}^*}^{-1}\bm{b}^\epsilon\right\|_{\infty}\\
    &\leq \sqrt{m}\sigma_{\max}\left((\bm{A}_{\mathcal{B}^*})^{-1}(\bm{A} + \bm{A}^{\epsilon})_{\mathcal{B}^*}^{-1}\right)\|\bm{A}^\epsilon _{\mathcal{B}^*}\|_{\infty}  +  \sqrt{m}\sigma_{\max}(\bm{A} + \bm{A}^{\epsilon})_{\mathcal{B}^*}^{-1}\left\|\bm{b}^\epsilon\right\|_{\infty}\\
    &\leq\frac{2\sqrt{m}}{\sigma^2}\|\bm{A}^\epsilon _{\mathcal{B}^*}\|_{\infty} + \frac{2\sqrt{m}}{\sigma}\|\bm{b}^\epsilon\|_{\infty}\\
    &\leq \frac{2m\sqrt{m}}{\sigma^2}\max_{i\in \mathcal{B}^*}\|\bm{A}_i^\epsilon\|_\infty   +  \frac{2\sqrt{m}}{\sigma}\left\|\bm{b}^\epsilon\right\|_{\infty} \leq \frac{\chi}{2}.
\end{align*}
The third and seventh line come from the sub-multiplicativity of matrix L$_\infty$ norm. The eighth line come from the definition of $\sigma$ following Assumption \ref{assp_nondeg} and the relation between the spectral norm $\sigma_{max}$ and L$_\infty$ norm. The last line is from the inequality  $||\bm{A}^{\epsilon}_{\mathcal{B}^*}||_{\infty} \leq m\max\limits_{i\in\mathcal{B}^*}\left\|\bm{A}^{\epsilon}_{i}\right\|_{\infty}$ and condition $\eqref{eqn:condition_stability}$. Thus we finish the part on the feasibility.

For part (c), we prove that reduced costs of non-basic variables in $\mathcal{B}'$ are all strictly negative. For a non-basic variable $x_i$, the reduced cost of the perturbed LP $\eqref{LP2}$ (denoted by $\hat{\Psi}_i$) can be expressed as follows,
\begin{equation}
\label{reduced_cost_decompose}
    \begin{aligned}
    \hat{\Psi}_i
    &\coloneqq
    c_i + c_i^\epsilon - ({\bm{c}}+\bm{c}^{\epsilon})_{\mathcal{B}^*}^{\top}(\bm{A}+\bm{A}^{\epsilon})_{\mathcal{B}^*}^{-1} (\bm{A}_i+\bm{A}^{\epsilon}_{i}) \\
    &=c_i - \bm{c}_{\mathcal{B}^*}^{\top}(\bm{A}_{\mathcal{B}^*})^{-1}\bm{A}_i + \bm{c}_{\mathcal{B}^*}^{\top}(\bm{A}_{\mathcal{B}^*})^{-1}\bm{A}_i + c_i^\epsilon - (\bm{c}+\bm{c}^{\epsilon})_{\mathcal{B}^*}^{\top}(\bm{A}+\bm{A}^{\epsilon})_{\mathcal{B}^*}^{-1} (\bm{A}_i+\bm{A}^{\epsilon}_{i})\\
    &= \Psi_i + c_i^\epsilon + \bm{c}_{\mathcal{B}^*}^{\top}\left((\bm{A}_{\mathcal{B}^*})^{-1}- (\bm{A}+\bm{A}^{\epsilon})_{\mathcal{B}^*}^{-1}\right)\bm{A}_i - \bm{c}^{\epsilon \top}_{\mathcal{B}^*}(\bm{A}+\bm{A}^{\epsilon})_{\mathcal{B}^*}^{-1} \bm{A}_i \\ 
    & \hspace{5mm} - \bm{c}_{\mathcal{B}^*}^{\top}(\bm{A}+\bm{A}^{\epsilon})_{\mathcal{B}^*}^{-1} \bm{A}_i^\epsilon - \bm{c}_{\mathcal{B}^*}^{\epsilon \top}(\bm{A}+\bm{A}^{\epsilon})_{\mathcal{B}^*}^{-1} \bm{A}_i^\epsilon.\\
    \end{aligned}
\end{equation}
Since the reduced cost of LP \eqref{LP1} $$\Psi_i\coloneqq c_i-\bm{A}_i^\top \bm{\lambda}^*>0$$
for $i\in\mathcal{B}'$, a sufficient condition for $\hat{\Psi}_i>0$ is the absolution values of all the rest five terms are no greater than $\frac{|\Psi_i|}{6}$. Next, we are going to bound each component in \eqref{reduced_cost_decompose}.

For $c_i^\epsilon$, the inequality $|c_i^\epsilon|\leq\frac{|\Psi_i|}{6}$ is directly implied from the condition $\eqref{eqn:condition_stability}$. For the second term in \eqref{reduced_cost_decompose}, we have
\begin{equation*}
    \begin{aligned}
    \left|\bm{c}_{\mathcal{B}^*}^{\top}\left((\bm{A}_{\mathcal{B}^*})^{-1}- (\bm{A}+\bm{A}^{\epsilon})_{\mathcal{B}^*}^{-1}\right)\bm{A}_i\right|
    &\leq \left\|\bm{c}_{\mathcal{B}^*}\right\|_1\cdot \left\|\left((\bm{A}_{\mathcal{B}^*})^{-1}- (\bm{A}+\bm{A}^{\epsilon})_{\mathcal{B}^*}^{-1}\right)\bm{A}_i \right\|_{\infty}\\
    &\leq m \left\|(\bm{A}_{\mathcal{B}^*})^{-1}\left(\bm{I}- \bm{A}_{\mathcal{B}^*}(\bm{A}+\bm{A}^{\epsilon})_{\mathcal{B}^*}^{-1}\right)\bm{A}_i \right\|_{\infty}\\
    &= m \left\|(\bm{A}_{\mathcal{B}^*})^{-1} \bm{A}^{\epsilon}_{\mathcal{B}^*}(\bm{A}+\bm{A}^{\epsilon})_{\mathcal{B}^*}^{-1}\bm{A}_i \right\|_{\infty} \\
    &\leq  m \left\|(\bm{A}_{\mathcal{B}^*})^{-1} (\bm{A}+\bm{A}^{\epsilon})_{\mathcal{B}^*}^{-1}\right\|_{\infty}\left\|\bm{A}^{\epsilon}_{\mathcal{B}^*}\right\|_{\infty}\left\|\bm{A}_i \right\|_{\infty} \\
    &\leq  m \left\|(\bm{A}_{\mathcal{B}^*})^{-1} (\bm{A}+\bm{A}^{\epsilon})_{\mathcal{B}^*}^{-1}\right\|_{\infty}\left\|\bm{A}^{\epsilon}_{\mathcal{B}^*}\right\|_{\infty} \\
    &\leq m\sqrt{m} \sigma_{\max}\left((\bm{A}_{\mathcal{B}^*})^{-1} (\bm{A}+\bm{A}^{\epsilon})_{\mathcal{B}^*}^{-1}\right)\left\|\bm{A}^{\epsilon}_{\mathcal{B}^*}\right\|_{\infty} \\
    &\leq \frac{2m\sqrt{m}}{\sigma^2} \left\|\bm{A}^{\epsilon}_{\mathcal{B}^*}\right\|_{\infty}\\
    &\leq \frac{2m^{2}\sqrt{m}}{\sigma^2} \max\limits_{i\in\mathcal{B}^*}\left\|\bm{A}^{\epsilon}_{i}\right\|_{\infty},\\
    \end{aligned}
\end{equation*}
where the first line is obtained by Holder's inequality, the fourth line is obtained by the sub-multiplicativity and \eqref{bdd_purturb_matrix_sigma}, and the sixth line comes from the relation between the spectral norm and $L_{\infty}$ norm, and the last line is again from $||\bm{A}^{\epsilon}_{\mathcal{B}^*}||_{\infty} \leq m
\max\limits_{i\in\mathcal{B}^*}\left\|\bm{A}^{\epsilon}_{i}\right\|_{\infty}$. Thus, from $\eqref{eqn:condition_stability}$ we have
$$\left|\bm{c}_{\mathcal{B}^*}^{\top}\left((\bm{A}_{\mathcal{B}^*})^{-1}- (\bm{A}+\bm{A}^{\epsilon})_{\mathcal{B}^*}^{-1}\right)\bm{A}_i\right|\leq\frac{\delta}{6}.$$
For the third term in \eqref{reduced_cost_decompose}, we have
\begin{equation*}
    \begin{aligned}
    \left|\bm{c}^{\epsilon \top}_{\mathcal{B}^*}(\bm{A}+\bm{A}^{\epsilon})_{\mathcal{B}^*}^{-1} \bm{A}_i\right|
    &\leq \left\|\bm{c}^{\epsilon}_{\mathcal{B}^*}\right\|_2\left\|(\bm{A}+\bm{A}^{\epsilon})_{\mathcal{B}^*}^{-1} \bm{A}_i\right\|_2\\
    &\leq \sqrt{m} \left\|\bm{c}^{\epsilon}_{\mathcal{B}^*}\right\|_{\infty}\sigma_{\max}\left((\bm{A}+\bm{A}^{\epsilon})_{\mathcal{B}^*}^{-1}\right)\|\bm{A}_i\|_2\\
    &\leq m \left\|\bm{c}^{\epsilon}_{\mathcal{B}^*}\right\|_{\infty}\sigma_{\max}\left((\bm{A}+\bm{A}^{\epsilon})_{\mathcal{B}^*}^{-1}\right)\\
    &\leq \frac{2m}{\sigma} \left\|\bm{c}^{\epsilon}_{\mathcal{B}^*}\right\|_{\infty}.\\
    \end{aligned}
\end{equation*}
 Thus, from $\eqref{eqn:condition_stability}$ again we have
$$\left|\bm{c}^{\epsilon \top}_{\mathcal{B}^*}(\bm{A}+\bm{A}^{\epsilon})_{\mathcal{B}^*}^{-1} \bm{A}_i\right| \leq\frac{\delta}{6}.$$

In a similar manner, for the last two terms we have
\begin{equation*}
    \begin{aligned}
        \left|\bar{\bm{c}}_{\mathcal{B}^*}^{\top}(\bm{A}+\bm{A}^{\epsilon})_{\mathcal{B}^*}^{-1} \bm{A}_i^\epsilon\right|
        &\leq
        \frac{2m}{\sigma}\max\limits_{i\in\mathcal{B}'} \left\|\bm{A}^{\epsilon}_{i}\right\|_{\infty},\\
        \left|\bm{c}_{\mathcal{B}^*}^{\epsilon \top}(\bm{A}+\bm{A}^{\epsilon})_{\mathcal{B}^*}^{-1} \bm{A}_i^\epsilon\right| 
        &\leq
        \frac{2m}{\sigma} \max\limits_{i\in\mathcal{B}'} \left\|\bm{A}^{\epsilon}_{i}\right\|_{\infty}\left\|\bm{c}^{\epsilon}_{\mathcal{B}^*}\right\|_{\infty}.
    \end{aligned}
\end{equation*}
Both of them are no larger than  $\frac{\delta}{6}$  condition because of $\eqref{eqn:condition_stability}$. 

Therefore, for any non-basic variable $x_i$ with $i\in\mathcal{B}'$, we conclude that its reduced cost in the perturbed LP $\hat{\Psi}_{i} \ge \frac{\delta}{6} > 0.$ Thus we establish the optimality of the solution $\hat{\bm{x}}.$ Lastly, given the stability of the sign for the reduced-cost, we know the optimal solution is unique.

\subsection{Proof of Lemma \ref{lem_gen_stability}}\label{ap_l1}
Firstly, we apply the result of Proposition \ref{prop_stability_general} to the DLP \eqref{eqn:DLP} (also the equivalent form \eqref{eqn:SDDLP}). Notice that in Proposition \ref{prop_stability_general} the matrix dimension is $m\times n$, while in the canonical form of the DLP \eqref{eqn:SDDLP}, the dimension becomes $(n+m) \times (2n+m)$. However, instead of plugging in the dimension to the result of Proposition \ref{prop_stability_general} directly, one can find that in \eqref{eqn:SDDLP}, only $\bm{C}$, $\bm{\mu}$ and $\bm{b}$ component will have random perturbation, and other entries will always be $1$. Having this observation, we know that under Assumption \ref{assp_dist}, \ref{assp_nondeg}, and the condition that
\begin{equation*}
    \begin{split}
        \left\|\hat{\bm{C}}_j - \bm{C}_j\right\|_{\infty}&
        \leq
        \left\{
        \begin{matrix}
            \frac{\min\{1,\sigma,\sigma^2\}\cdot \min\{\chi,\delta\}}{12n^2\sqrt{n+m}}, & \text{\ for $j\in\mathcal{J}^*$},\\
            \frac{\sigma\delta}{12\sqrt{n(n+m)}}, & \text{\ for $j\in \mathcal{J}'$},
        \end{matrix}
        \right.\\
        |\hat{\bm{\mu}}_j - \bm{\mu}_{j}| &
        \leq
        \left\{
        \begin{matrix}
            \frac{\sigma\delta}{12\sqrt{n(n+m)}}, & \text{\ for $j\in\mathcal{J}^*$},\\
            \frac{\delta}{6}, & \text{\ for $j\in \mathcal{J}'$},
        \end{matrix}
        \right.\\    
        |\hat{\bm{b}}_i - \bm{b}_i|&
        \leq \frac{\sigma\chi}{8\sqrt{n+m}}, \,\,\, \text{\ for $i\in [m]$}\\
    \end{split}
\end{equation*}
the optimal solution to LP \eqref{eqn:SDDLP_1} is unique and it shares the same optimal basis and binding/non-binding structure with the solution of LP \eqref{eqn:SDDLP} (which equivalent to \eqref{eqn:DLP}). The last step is to show that the condition above could be extended to
\begin{equation*}
    \begin{split}
        \left\|\hat{\bm{C}}_j - \bm{C}_j\right\|_{\infty}&
        \leq
        \left\{
        \begin{matrix}
            \frac{\min\{1,\sigma,\sigma^2\}\cdot \min\{\chi,\delta\}}{12n^2\sqrt{n+m}}, & \text{\ for $j\in\mathcal{J}^*$},\\
            \frac{\sigma\delta}{12\sqrt{n(n+m)}}, & \text{\ for $j\in \mathcal{J}'$},
        \end{matrix}
        \right.\\
        |\hat{\bm{\mu}}_j - \bm{\mu}_{j}| &
        \leq
        \left\{
        \begin{matrix}
            \frac{\sigma\delta}{12\sqrt{n(n+m)}}, & \text{\ for $j\in\mathcal{J}^*$},\\
            \frac{\delta}{6}, & \text{\ for $j\in \mathcal{J}'$},
        \end{matrix}
        \right.\\    
        |\hat{\bm{b}}_i - \bm{b}_i|&
        \leq \frac{\sigma\chi}{8\sqrt{n+m}}, \,\,\, \text{\ for $i\in \mathcal{I}^*$}\\
        \hat{\bm{b}}_i - \bm{b}_i&
        \geq -\frac{\sigma\chi}{8\sqrt{n+m}}, \,\,\, \text{\ for $i\in \mathcal{I}'$}.\\
    \end{split}
\end{equation*}
To show this statement, we only have to look at the dual solution of the LP \eqref{eqn:SDDLP_1}. Since the primal is nondegenerate, its dual solution must be unique. It is easy to construct a primal solution for a larger $\hat{b}_i$ where $i \in \mathcal{I}'$: we just add the corresponding increased value to the primal slack variable, and it remains to show that the constructed primal solution is unique. The increase for $\hat{b}_i$ where $i \in \mathcal{I}'$ will not improve objective value because it corresponds to the dual variable that $\lambda_i^* = 0$, and moreover, it will not bring any new optimal dual solutions (otherwise the original optimal primal solution set would be different). Therefore, the optimal primal solution remains unique for $\hat{b}_i > b_i-\frac{\sigma\chi}{8\sqrt{n+m}}$. Therefore if we take 
$$L = \min\left\{\frac{\min\{1,\sigma,\sigma^2\}\cdot \min\{\chi,\delta\}}{12n^2\sqrt{n+m}}, \frac{\sigma\delta}{12\sqrt{n(n+m)}}, \frac{\sigma\delta}{12\sqrt{n(n+m)}}, \frac{\delta}{6}, \frac{\sigma\chi}{8\sqrt{n+m}}\right\},$$
we finish proving the statement of the lemma.
Here the parameters $\sigma, \chi, \delta, $ follow their definitions in Proposition \ref{prop_stability_general}, but are based on the standard form DLP \eqref{eqn:SDDLP}.

\end{document}